\documentclass[12pt]{article} 

\usepackage[utf8]{inputenc}
\usepackage[T1]{fontenc}

\usepackage{balance}
\usepackage{graphics}
\usepackage{graphicx}
\usepackage{amsfonts}
\usepackage{amssymb}
\usepackage{verbatim}     
\usepackage{amsxtra}%
\usepackage{epsfig}
\usepackage{latexsym}
\usepackage{ifthen}
\usepackage{psfrag}
\usepackage{xcolor}
\usepackage{enumitem}
\usepackage{textcomp}
\usepackage{epstopdf}
\usepackage{algorithm}
\usepackage{algorithmic}
\usepackage{esdiff}
\usepackage{subfig}
\usepackage{multirow}

\usepackage{appendix}

\usepackage{hslTR}

\usepackage{rgsEnvironments}

\usepackage{rgsMacros}

\newcommand{\eoe}
{\hspace*{\fill}{$\vcenter{\hrule height1pt 
			\hbox{\vrule width1pt height3pt 
				\kern3pt \vrule width1pt} \hrule height1pt}$} }

\newcommand{\T}{\mathcal{T}}

\newcommand{\xp}{z}

\newcommand{\up}{u}

\newcommand{\fp}{F_P}

\newcommand{\Cp}{C_P}

\newcommand{\hp}{h}

\newcommand{\hcq}{\kappa_{\xlogic}}
\newcommand{\hczero}{\kappa_0}

\newcommand{\hcone}{\kappa_1}

\newcommand{\xlogic}{q}

\newcommand{\xlogicSpace}{Q}

\newcommand{\f}{F}
\newcommand{\g}{G}

\def\ba{\begin{array}}
	\def\ea{\end{array}}

\renewcommand{\scriptsize}{\tiny}

\newcommand{\minSet}{z_1^*} 

\allowdisplaybreaks

\begin{document}
	
\IfConf{
\begin{frontmatter}
	\runtitle{Uniting Nesterov and Heavy Ball}  
	
	\title{Uniting Nesterov and Heavy Ball Methods for Uniform Global Asymptotic Stability of the Set of Minimizers\thanksref{footnoteinfo}} 
	
	\thanks[footnoteinfo]{This paper was not presented at any IFAC 
		meeting. Corresponding author D.~M.~Hustig-Schultz.\\
		This research has been partially supported by the National Science Foundation under Grant no. ECS-1710621, Grant no. CNS-2039054, and Grant no. CNS-2111688, by the Air Force Office of Scientific Research under Grant nos. FA9550-19-1-0169, FA9550-20-1-0238, and FA9550-23-1-0145, by the Air Force Research Laboratory under Grant nos. FA8651-22-1-0017 and FA8651-23-1-0004, and by the Army Research Office under Grant no. W911NF-20-1-0253.}
	
	\author[UCSC]{Dawn M. Hustig-Schultz}\ead{dhustigs@ucsc.edu},    
	\author[UCSC]{Ricardo G. Sanfelice}\ead{ricardo@ucsc.edu}                
	
	\address[UCSC]{Department of Electrical and Computer Engineering, University of California, Santa Cruz.}  


\begin{abstract} 
	We propose a hybrid control algorithm that guarantees fast convergence and uniform global asymptotic stability of the unique minimizer of a $\mathcal{C}^1$, convex objective function. The
	algorithm, developed using hybrid system tools, employs a uniting control strategy, in which Nesterov's accelerated gradient descent is used ``globally'' and the heavy ball method is used
	``locally,'' relative to the minimizer. Without knowledge of its location, the proposed hybrid control strategy switches between these accelerated methods to ensure convergence to the minimizer without oscillations, with a (hybrid) convergence rate that preserves the convergence rates of the individual optimization algorithms. We analyze key properties of the resulting closed-loop system including existence of solutions, uniform global asymptotic stability, and convergence rate. Additionally, stability properties of Nesterov’s method are analyzed, and extensions on convergence rate results in the existing literature are presented. Numerical results validate the findings.
\end{abstract}
	
\end{frontmatter}
}{
		\ititle{Uniting Nesterov and Heavy Ball Methods for Uniform Global Asymptotic Stability of the Set of Minimizers}
		\iauthor{
			Dawn M. Hustig-Schultz \\
			{\normalsize dhustigs@ucsc.edu} \\
			Ricardo G. Sanfelice \\
			{\normalsize ricardo@ucsc.edu}}
		\idate{\today{}} 
		\iyear{2021}
		\irefnr{04}
		\makeititle
		
		\begin{abstract} 
			We propose a hybrid control algorithm that guarantees fast convergence and uniform global asymptotic stability of the unique minimizer of a $\mathcal{C}^1$, convex objective function. The
			algorithm, developed using hybrid system tools, employs a uniting control strategy, in which Nesterov's accelerated gradient descent is used ``globally'' and the heavy ball method is used
			``locally,'' relative to the minimizer. Without knowledge of its location, the proposed hybrid control strategy switches between these accelerated methods to ensure convergence to the minimizer without oscillations, with a (hybrid) convergence rate that preserves the convergence rates of the individual optimization algorithms. We analyze key properties of the resulting closed-loop system including existence of solutions, uniform global asymptotic stability, and convergence rate. Additionally, stability properties of Nesterov’s method are analyzed, and extensions on convergence rate results in the existing literature are presented. Numerical results validate the findings and demonstrate the robustness of the uniting algorithm.
		\end{abstract}
	}

\section{Introduction}
\label{sec:Intro}

\subsection{Background and Motivation}
\label{sec:Background}
\IfConf{We propose an algorithm that solves optimization problems of the form $\min_{\xi \in \reals^n} L(\xi)$ with accelerated gradient methods.}{There has been growing interest in analyzing accelerated gradient methods from a dynamical systems perspective 
\cite{muehlebach2019dynamical}, which permits the use of well established analysis tools, such as Lyapunov theory, to study convergence and stability properties of accelerated algorithms \cite{polyak2017lyapunov}, \cite{attouch2000heavy}\IfConf{}{, \cite{su2016differential}, \cite{krichene2015accelerated}, \cite{wibisono2016variational}, \cite{wilson2016lyapunov}}.} The {\em heavy ball method} is an accelerated gradient method that guarantees convergence to the minimizer \IfConf{$\xi^*$}{} of a convex function $L$ \cite{polyak1964some}, and that achieves a faster convergence rate than classical gradient descent by adding a ``velocity'' term to \IfConf{$\nabla L$}{the gradient}. The dynamical system characterization for this method is
\begin{equation}\label{eqn:HBF}
	\ddot{\xi} + \lambda \dot{\xi} + \gamma \nabla L(\xi) = 0
\end{equation}
where $\lambda$ and $\gamma$ are positive tunable parameters that represent friction and gravity, respectively; see \cite{attouch2000heavy}, \cite{polyak2017lyapunov}. 
\IfConf{In \cite{ghadimi2015global} and \cite{fazlyab2018design} it is shown that the discrete-time version of the heavy ball method converges exponentially when $L$ is strongly convex with a Lipschitz continuous gradient, and \cite{ghadimi2015global} shows convergence with rate $\frac{1}{k}$ when $L$ is convex. It is shown in \cite{lessard2016analysis} that for strongly convex $L$ with Lipschitz continuous $\nabla L$ global convergence of the discrete-time heavy ball method can only be guaranteed for condition numbers of about 18 or less, and it is found in \cite{badithela2019analysis} that the exact condition number of $9 + 5\sqrt{14} \approx 17.94$ denotes such a boundary between global convergence and non-convergence, for such objective functions.}{In \cite{ghadimi2015global} and \cite{fazlyab2018design} it is shown that the discrete-time version of the heavy ball method converges exponentially when $L$ is strongly convex with a Lipschitz continuous gradient, and \cite{ghadimi2015global} shows convergence with rate $\frac{1}{k}$ when $L$ is convex. It is shown in \cite{lessard2016analysis} that for strongly convex $L$ with Lipschitz continuous $\nabla L$ global convergence of the discrete-time heavy ball method can only be guaranteed for condition numbers of about 18 or less, and it is found in \cite{badithela2019analysis} that the exact condition number of $9 + 5\sqrt{14} \approx 17.94$ denotes such a boundary between global convergence and non-convergence, for such objective functions.} \IfConf{W}{For the case w}hen $L$ is strongly convex, and inspired by the heavy ball algorithm, two algorithms with a resettable velocity term are proposed in \cite{le2021hybrid} and shown to guarantee exponential convergence. In \cite{sebbouh2020convergence}, however, it was \IfConf{shown}{demonstrated} that the heavy ball algorithm converges exponentially for convex $L$ when \IfConf{$L$}{such an objective function} also has the property of quadratic growth away from \IfConf{$\xi^*$.}{its minimizer.} Global asymptotic stability of \IfConf{$\xi^*$,}{the minimizer,} which is the property that all solutions that start close to \IfConf{$\xi^*$}{the minimizer} stay close, and solutions from all initial conditions converge to \IfConf{$\xi^*$,}{the minimizer,} is demonstrated in \IfConf{}{\cite{michalowsky2014multidimensional},} \cite{michalowsky2016extremum}, \IfConf{for convex and smooth $L$.}{when $L$ is convex and smooth.}
\IfConf{}{The work in \cite{polyak2017lyapunov} provides several Lyapunov functions to establish global asymptotic stability of the minimizer and convergence rates for the heavy ball method, both when $L$ is strongly convex and when $L$ is convex.} \IfConf{\vspace{-0.3cm}}{}

Another powerful accelerated method is {\em Nesterov's accelerated gradient descent}.
One characterization of the dynamical system for Nesterov's method, for convex $L$, proposed in \cite{muehlebach2019dynamical}, is 
\IfConf{
	
	\vspace{-0.55cm}
	
}{}
\begin{equation} \label{eqn:MJODENCVX}
	\ddot{\xi} + 2\bar{d}(t)\dot{\xi} + \frac{1}{M\zeta^2}\nabla L(\xi + \bar{\beta}(t) \dot{\xi}) = 0,
\end{equation}
%
%
where $M > 0$ is the Lipschitz constant of \IfConf{$\nabla L$}{the gradient of $L$} and where the constant $\zeta > 0$ rescales time in solutions to \eqref{eqn:MJODENCVX}. The dynamical system in \eqref{eqn:MJODENCVX} resembles the model of a mass-spring-damper, with a curvature-dependent damping term where the total damping is a linear combination of $\bar{d}(t)$ and $\bar{\beta}(t)$. In \cite{muehlebach2019dynamical}, the convergence rate of Nesterov's method is characterized as $\frac{1}{(t+2)^2}$ for \eqref{eqn:MJODENCVX} (for $t \geq 1$), when $\zeta = 1$, \IfConf{$\xi^* = 0$, and $L(\xi^*) = 0$. The stability properties of \eqref{eqn:MJODENCVX} are not revealed in \cite{muehlebach2019dynamical}, however.}{and when the minimizer is the origin, at which $L$ is zero.} 
\IfConf{
	
	\vspace{-0.3cm}
	
}{}

\IfConf{}{Exponential convergence of the discrete-time analogue of Nesterov's\\ method is established for strongly convex $L$ with Lipschitz continous $\nabla L$ in \cite{fazlyab2018design} and \cite{lessard2016analysis}.
The earliest dynamical system characterization for Nesterov's algorithm was proposed in \cite{su2016differential}, including a variation for higher friction, and the proposed characterizations were shown to have a  convergence rate of $\frac{1}{t^2}$. In \cite{krichene2015accelerated}, the analysis of the dynamical system in \cite{su2016differential} is extended to include optimization of objective functions $L$ with non-Euclidean geometries, 
and this dynamical system is combined with mirror descent to design an accelerated mirror descent ODE, with a convergence rate of $\frac{1}{t^2}$. In \cite{wibisono2016variational}, a dynamical system, consisting of an Euler-Lagrange equation, is derived for Nesterov's algorithm via a Bregman Lagrangian, with an exponential rate of convergence under ideal scaling and a rate of convergence of $\frac{1}{t^p}$ with $p \geq 2$ for a polynomial subclass of such a dynamical system.}

\IfConf{}{In \cite{wibisono2016variational} an exponential rate of convergence for such a system under ideal scaling is provided, and, for a polynomial class of dynamical systems, a convergence rate of $\frac{1}{t^p}$ with $p \geq 2$ is shown. 
In \cite{kolarijani2019continuous} and \cite{kolarijani2018fast}, two hybrid algorithms based on the ODE in \cite{su2016differential} are presented: one with a state-dependent, time-invariant damping input and another with an input that controls the magnitude of the gradient term. The algorithms require the objective function to satisfy the Polyak-\L ojasiewicz inequality, which includes a subclass of nonconvex functions in which all stationary points are global minimizers. 
The authors in \cite{poveda2019inducing} propose two hybrid reset algorithms based on the ODE in \cite{su2016differential}, HAND-1 and HAND-2, which yield an exponential convergence rate for strongly convex $L$ and a rate of $\frac{1}{t^2}$ for convex $L$, with the latter rate only assured until the first reset.}

\IfConf{While the results in \cite{fazlyab2018design}, \cite{lessard2016analysis}, and \cite{muehlebach2019dynamical} characterize the convergence properties of Nesterov's method (or a variation of) the stability properties of the method are not revealed.}{While the results in \cite{muehlebach2019dynamical}, \cite{su2016differential}, \cite{krichene2015accelerated}, \cite{wibisono2016variational}, \cite{fazlyab2018design}, and \cite{lessard2016analysis} characterize the convergence properties of Nesterov's method (or a variation of) the stability properties of the method are not revealed.} A particularly useful property for optimization algorithms, called {\em uniform global asymptotic stability} (UGAS), requires that solutions reach a neigborhood of \IfConf{$\xi^*$}{the minimizer} in time that is uniform on the set of initial conditions \cite{poveda2019inducing}, \cite{poveda2021robust}, \cite{teel2019first}. After finite time, the error of such solutions becomes smaller than a given threshold \cite{220}. Due to such a guarantee\IfConf{,}{ for solutions,} UGAS is typically useful for certifying robustness to small perturbations in time-varying dynamical and hybrid systems \cite{65}, \cite{220}. Remarkably, the algorithm with resets in the velocity term proposed in \cite{le2021hybrid} can be shown to induce UGAS of \IfConf{$\xi^*$}{the minimizer} (with zero velocity term) and reduced oscillations, for the particular case when L is strongly convex. The algorithm with resets in \cite{teel2019first} can be shown to induce UGAS of \IfConf{$\xi^*$}{the minimizer} when $L$ is invex, has an exponential convergence rate when $L$ satisfies the Polyak-\L ojasiewicz inequality, and uniform global exponential stability (UGES) when $L$ is strongly convex. Unfortunately, as shown in \cite{poveda2019inducing}, via a counterexample, Nesterov-like algorithms do not necessarily assure UGAS of \IfConf{$\xi^*$}{the minimizer} when $L$ is convex. In response to this, \cite{poveda2019inducing} proposes the HAND-1 and HAND-2 reset algorithms\IfConf{,}{, based on the ODE in \cite{su2016differential},} and prove UGAS of \IfConf{$\xi^*$}{the minimizer} for both algorithms. The exponential convergence rate of HAND-2, however, only applies to strongly convex $L$, and the convergence rate of $\frac{1}{t^2}$ for HAND-1, for convex $L$, only holds up until the first reset. 
\IfConf{
	
	\vspace{-0.2cm}
	
	\begin{figure}[thpb] 
		\centering
		\setlength{\unitlength}{1.0pc} 
		
		\begin{picture}(20,15)(0,0)
			\footnotesize
%
%
			\put(0.5,0){\includegraphics[scale=0.3,trim={0.8cm 0.4cm 0.6cm 0.5cm},clip,width=20\unitlength]{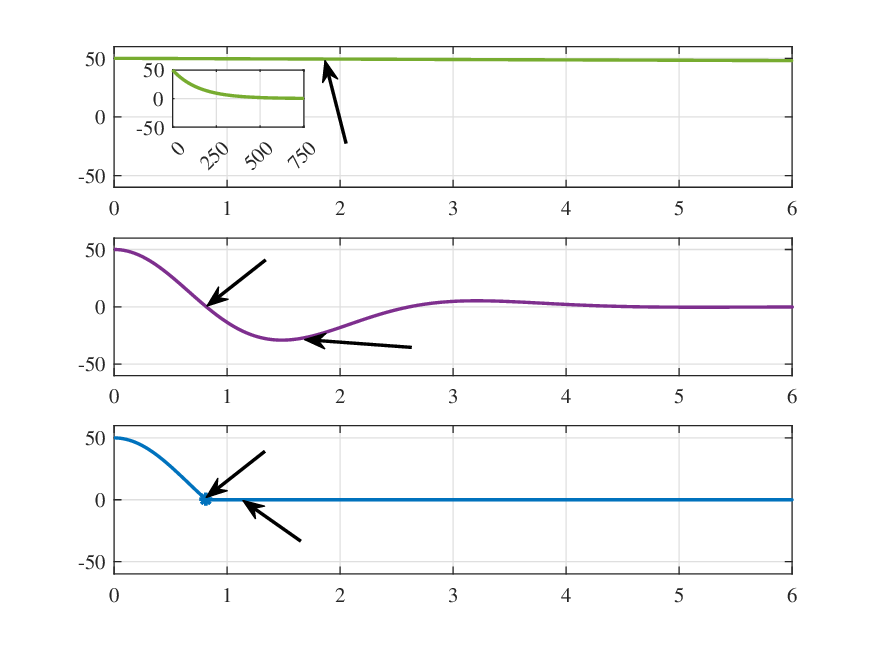}}
			\put(0.8,3.1){$\xi$}
			\put(0.8,8){$\xi$}
			\put(0.8,12.8){$\xi$}
			\put(10.3,0.3){$t [s]$}
			\put(8.2,12.5){slow convergence}
			\put(8.2,11.7){without oscillations}
			\put(9.9,7){oscillations}
			\put(6.2,9.2){fast convergence}
			\put(7,2){no oscillations}
			\put(6.2,4.5){fast convergence}	
		\end{picture}
	
		\vspace{-0.3cm}
	
		\caption{Comparison of the performance of heavy ball, with large $\lambda$, Nesterov's method, and the proposed logic-based algorithm. The objective function is $L(\xi) = \xi^2$. Top: the heavy ball algorithm, with large $\lambda$, converges very slowly. Top inset: zoomed out view of heavy ball. Second from top: Nesterov's method converges quickly, but with oscillations. Third from top: our proposed logic-based algorithm yields fast convergence, with no oscillations.}
		\label{fig:MotivationalPlot}
	\end{figure}}{
	\begin{figure}[thpb] 
		\centering
		\setlength{\unitlength}{1.0pc} 
		
		\begin{picture}(30.8,12)(0,0)
			\footnotesize
			%
			%
			\put(0.7,0.5){\includegraphics[scale=0.3,trim={0.9cm 0.25cm 1.25cm 0.5cm},clip,width=14\unitlength]{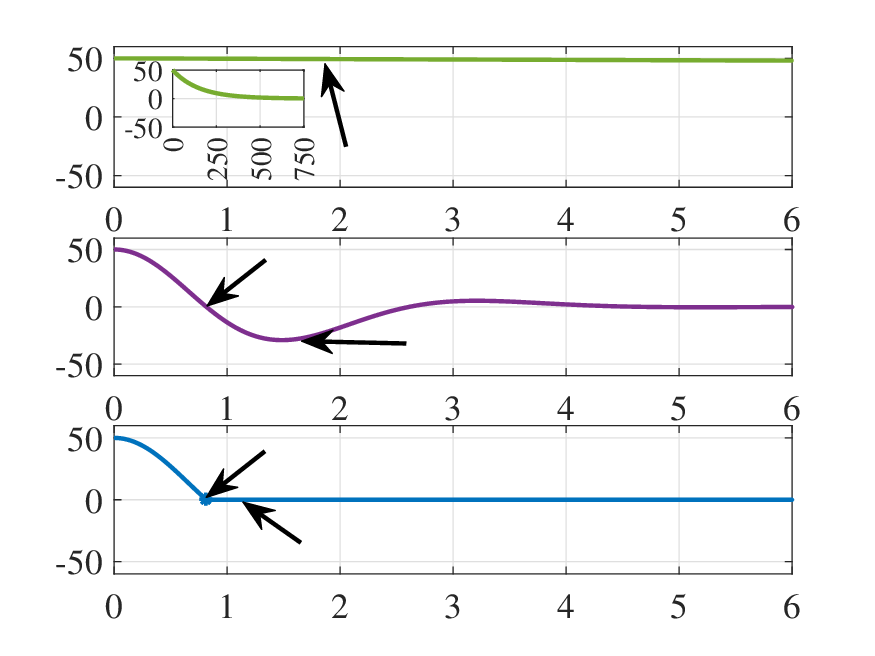}}
			\put(16.5,0.5){\includegraphics[scale=0.3,trim={0.6cm 0cm 1.25cm 0.5cm},clip,width=14\unitlength]{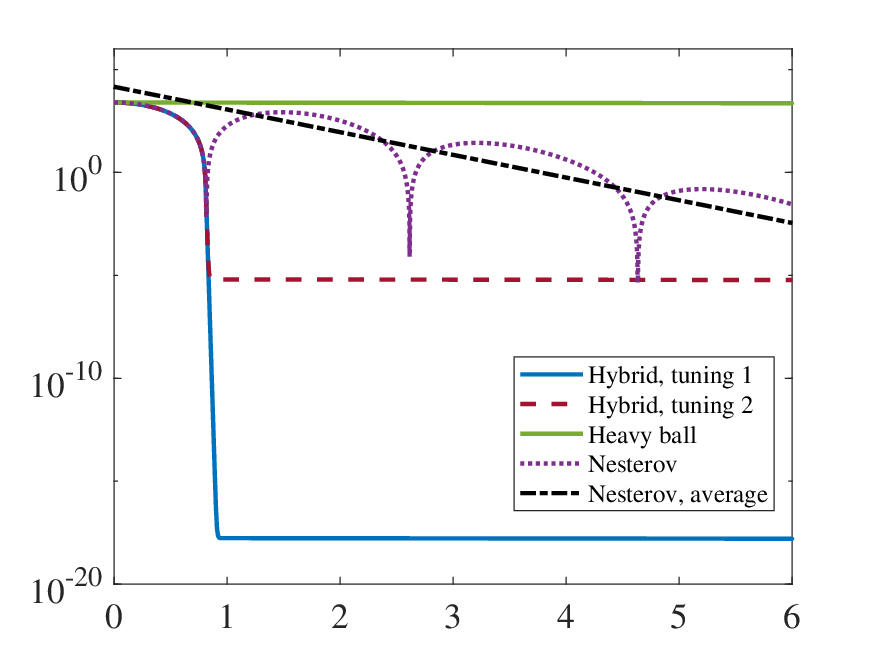}}
			\put(0,3){$\xi$}
			\put(0,6.6){$\xi$}
			\put(0,10.1){$\xi$}
			\put(7.7,0.2){$t [s]$}
			\put(23.7,0.3){$t [s]$}
			\put(15.5,4.8){\rotatebox{90}{$L(\xi) - L^*$}}
			\put(6.4,9.9){slow convergence}
			\put(6.4,9.2){without oscillations}
			\put(7.4,5.9){oscillations}
			\put(4.8,7.4){fast convergence}
			\put(5.5,2){no oscillations}
			\put(4.8,3.8){fast convergence}
		\end{picture}
		\caption{Comparison of the performance of heavy ball, with large $\lambda$, Nesterov's method, and the proposed logic-based algorithm. The objective function is $L(\xi) = \xi^2$. Top left: the heavy ball algorithm, with large $\lambda$, converges very slowly. Top inset: zoomed out view of heavy ball. Middle left: Nesterov's method converges quickly, but with oscillations. Bottom left: our proposed logic-based algorithm yields fast convergence, with no oscillations. Right: comparison of the value of $L(\xi) - L^*$ (in log scale) versus time for each algorithm. Different tunings of the logic-based algorithm's parameters leads to modifications of the solution's profile.}
		\label{fig:MotivationalPlotTR}
\end{figure}}
%
%
The work in this paper is motivated by the lack of an accelerated gradient algorithm assuring UGAS, with a convergence rate that holds for all time and that resembles that of Nesterov’s method (at least far from \IfConf{$\xi^*$}{the minimizer}), when \IfConf{$L$}{the objective function} is convex. However, attaining such a rate is expected to lead to oscillations, which are typically seen in accelerated gradient methods. The performance of the heavy ball method\IfConf{}{, for instance,} depends highly on the choice of $\lambda$\IfConf{.}{ and $\gamma$.}
\IfConf{For rather simple choices of}{In particular, for a fixed value of $\gamma$,
the choice of the friction parameter $\lambda$
significantly affects the asymptotic behavior of the solutions
to \eqref{eqn:HBF}. For rather simple choices of} \IfConf{$L$}{the function $L$},\IfConf{}{ 
the literature on this method indicates that}
large values of $\lambda$\IfConf{}{ are seen to} give rise to slowly converging 
solutions \IfConf{\cite{attouch2000heavy}.}{
resembling solutions yielded by steepest descent \cite{attouch2000heavy}.}
The top plot\footnote{Code at\IfConf{\\}{} \texttt{gitHub.com/HybridSystemsLab/UnitingMotivation}}\IfConf{}{ on the left} in \IfConf{Fig. \ref{fig:MotivationalPlot}}{Figure \ref{fig:MotivationalPlotTR}} demonstrates such behavior. In contrast, smaller values of $\lambda$ give rise 
to fast solutions with oscillations getting wilder as $\lambda$ decreases 
\cite{attouch2000heavy}. Nesterov's method converges quickly but also suffers from oscillations \IfConf{\cite{muehlebach2019dynamical}}{\cite{su2016differential}}, as the coefficient of the velocity term starts small and tends toward zero (but being always positive) as $t$ tends to infinity. \IfConf{Such behavior of}{The oscillatory behavior of} Nesterov's method, with $\zeta = 2$, is shown in the \IfConf{second plot from the top}{middle plot on the left} in \IfConf{Fig. \ref{fig:MotivationalPlot}}{Figure \ref{fig:MotivationalPlotTR}}. \IfConf{
	
	\vspace{-0.2cm}
	
}{}

Due to its implications on robustness, we are particularly interested in an algorithm that assures UGAS of \IfConf{$\xi^*$}{the minimizer of $L$} with a rate of convergence that holds for all time, and without\IfConf{}{ the undesired} oscillations. As pointed out in Section \ref{sec:Background}, these properties are not guaranteed by Nesterov's method. The behavior shown in the \IfConf{first and second}{top and middle} plots in \IfConf{Fig. \ref{fig:MotivationalPlot}}{Figure \ref{fig:MotivationalPlotTR}} motivates the logic-based algorithm proposed in this paper. The proposed algorithm exploits the main features of heavy ball and Nesterov's method to achieve fast convergence and UGAS of \IfConf{$\xi^*$.}{the minimizer.} More precisely, without knowledge of the location of \IfConf{$\xi^*$,}{the minimizer,} it selects Nesterov's method to converge quickly to nearby \IfConf{$\xi^*$}{the minimizer} and, once solutions reach a neighborhood of \IfConf{$\xi^*$,}{the minimizer,} switches to the heavy ball method with large $\lambda$ to avoid oscillations. \IfConf{Such logic-based algorithms, or {\em uniting algorithms}, were first proposed in \cite{teel1997uniting} and \cite{teel1997uniting2}. General uniting algorithms, with examples, are discussed in \cite{65} and \cite{220}.}{Such logic-based algorithms, or {\em uniting algorithms}, were first proposed in \cite{teel1997uniting} and \cite{teel1997uniting2}. General uniting algorithms, with examples, are discussed in \cite{65} and \cite{220}.} We use the hybrid systems framework for our proposed algorithm, as hybrid systems utilize hysteresis to avoid chattering at the switching boundary; see \cite{le2021hybrid}, \cite{220}, \cite{65}, \cite{teel2019first}.
An example solution to our proposed\IfConf{}{ logic-based} algorithm, shown in the \IfConf{third plot from the top in}{bottom plot on the left of} \IfConf{Fig. \ref{fig:MotivationalPlot}}{Figure \ref{fig:MotivationalPlotTR}}, demonstrates the improvement obtained\IfConf{,}{
by using Nesterov's method globally and the heavy ball method locally,}
under relatively mild assumptions on \IfConf{$L$}{the objective function $L$}. The proposed algorithm guarantees UGAS and a (hybrid) convergence rate that holds for all $t \geq 0$. 
\IfConf{
	
\vspace{-0.2cm}

}{}
\subsection{Contributions}
\label{sec:Contributions}
\IfConf{
	
\vspace{-0.2cm}

}{}
The main contributions of this paper are as follows. \IfConf{
	
	\vspace{-0.45cm}
	
\begin{enumerate}[label={\arabic*)},leftmargin=*]
	\item {\em A uniting algorithm for fast convergence and UGAS of $\xi^*$.} 
	
	\item {\em Well-posedness, existence of solutions, and robustness to small perturbations in measurements of $\nabla L$.}
	Nesterov's method can suffer from error accumulation, due to its velocity term \cite{flammarion2015averaging}. To overcome this issue, in Section \ref{sec:NonstronglyConvexNest} we prove well-posedness and existence of solutions for the proposed hybrid closed-loop algorithm. 
	Due to such well-posedness, the established UGAS property is robust to small perturbations in measurements of \IfConf{$\nabla L$}{the gradient of $L$} \cite[Theorem~7.21]{65}. 
	\begin{figure}[thpb]
		\centering
		\setlength{\unitlength}{1.0pc} 
				
				\vspace{-0.2cm}
		
		\begin{picture}(20,15)(0,0)
			\footnotesize
%
%
			\put(0,0.5){\includegraphics[trim={0.5cm 0.3cm 1cm 0.6cm},clip,width=20\unitlength]{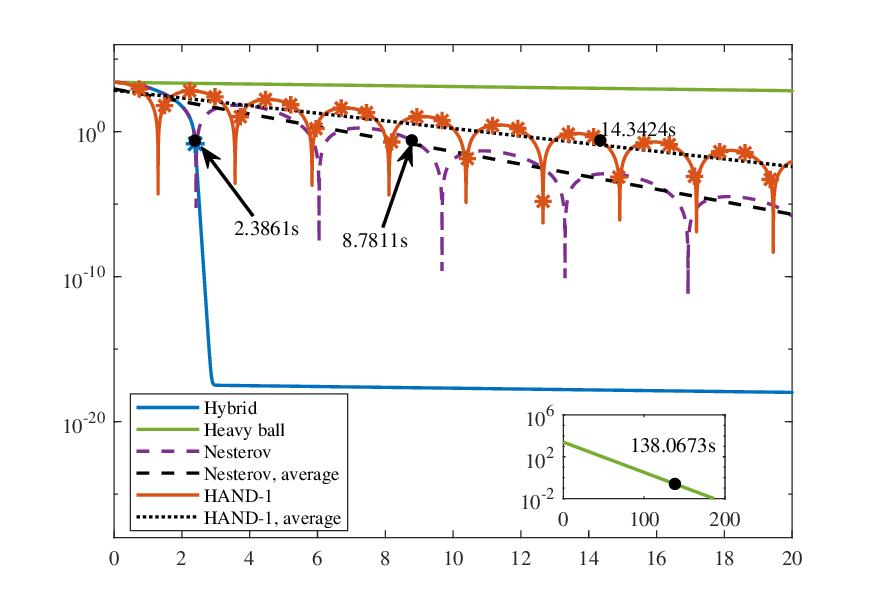}}
			\put(0,6.2){\rotatebox{90}{$L(\xi) - L^*$}}
			\put(10.2,0.4){$t[s]$}
		\end{picture}
	
		\vspace{-0.3cm}
	
		\caption{A comparison of the evolution of $L$ over time for Nesterov's method in \eqref{eqn:MJODE_ZetaNum}, heavy ball, HAND-1 from \cite{poveda2019inducing}, and our proposed uniting algorithm, for a function $L(\xi) := \xi^2$, with a single minimizer at $\xi^* = 0$. 
		As opposed to Fig. \ref{fig:MotivationalPlot}, which uses $\zeta = 2$ for $\HS_1$, this example uses $\zeta = 1$, which results in slower convergence of solutions to $\HS$ and $\HS_1$.}
		\label{fig:NSCLargerZeta}
	\end{figure}
	
	\item {\em A (hybrid) convergence rate preserving the rates of Nesterov's method and heavy ball.} 
	Numerical simulations\footnote{Code at \texttt{gitHub.com/HybridSystemsLab/UnitingTradeoff} \label{foot:Tradeoff}} in Fig. \ref{fig:NSCLargerZeta} and Section \ref{sec:Examples} show improved performance of our uniting algorithm over HAND-1.
	
	\item {\em Extension of the results on Nesterov's method in \cite{muehlebach2019dynamical}.} 
	We achieve such an extension by moving $\zeta$ into the numerator of the coefficient of $\nabla L$, effectively decoupling $\zeta$ and $M$, namely, \IfConf{
	
	\vspace{-0.2cm}
	
	}{}
	\begin{equation} \label{eqn:MJODE_ZetaNum}
		\ddot{\xi} + 2\bar{d}(t)\dot{\xi} + \frac{\zeta^2}{M}\nabla L(\xi + \bar{\beta}(t) \dot{\xi}) = 0.
	\end{equation}
\end{enumerate}}{
\begin{enumerate}[label={\arabic*)},leftmargin=*]
	\item {\em A uniting algorithm for fast convergence and UGAS of the minimizer:} 
	In Section \ref{sec:NonstronglyConvexNest} we propose a uniting algorithm that solves optimization problems of the form $\min_{\xi \in \reals^n} L(\xi)$
	with accelerated gradient methods.
	Designed using hybrid system tools,
	the algorithm unites Nesterov's method in \eqref{eqn:MJODE_ZetaNum_TR} globally and the heavy ball method in \eqref{eqn:HBF} with large $\lambda$ locally to guarantee fast convergence with UGAS of the minimizer $\xi^*$ of a convex objective function $L$; see Sections \ref{sec:NonstronglyConvexNest} and \ref{sec:ProofMainResult}. The establishment of UGAS solves the difficult problem of achieving such a property for Nesterov-like algorithms \cite{poveda2019inducing}, \cite{poveda2020heavy}.
	The algorithm we propose exploits measurements of $\nabla L$ and requires no knowledge of $L^* := L(\xi^*)$ or $\xi^*$. In practice, such measurements of $\nabla L$ are typically approximated from measurements of $L$. The algorithm, however, does not require measurements of the Hessian of $L$. 
	
	\item {\em Well-posedness and existence of solutions:} Nesterov's method can suffer from error accumulation, due to its velocity term \cite{flammarion2015averaging}. To overcome this issue, in Section \ref{sec:NonstronglyConvexNest} we prove well-posedness and existence of solutions for the proposed hybrid closed-loop algorithm. Hybrid systems that are {\em well-posed} are defined to be those hybrid systems, vaguely speaking, for which graphical limits of graphically convergent sequences of solutions, with no perturbations and with vanishing perturbations, respectively, are still solutions \cite[Chapter~6]{65}. It is important for our algorithm to be well-posed as we want to ensure robustness to small noise in measurements of the gradient of $L$.
	
	\item {\em Robustness to small perturbations:} Due to the well-posedness of the proposed hybrid uniting algorithm, we show that the established UGAS property is robust to small perturbations in measurements of the gradient of $L$ \cite[Theorem~7.21]{65}. We illustrate this robustness in Section \ref{sec:Examples} via numerical simulations that include small noise in measurements of the gradient. 
	
	\item {\em A (hybrid) convergence rate preserving the rates of Nesterov's method and heavy ball:} In Sections \ref{sec:NonstronglyConvexNest} and \ref{sec:ProofMainResult} we show that our uniting algorithm attains a rate of $\frac{1}{(t+2)^2}$ for the global algorithm and $\exp \left(-(1-m)\psi t\right)$, where $m \in (0,1)$ and $\alpha > 0$ are such that $\psi := \frac{m\alpha\gamma}{\lambda} > 0$ and\footnote{Although the constant $\nu$ does not appear in the rate for the local algorithm, such a constant is used in the forthcoming Proposition \ref{prop:HBFConvergenceRate} to derive this rate.} $\nu := \psi (\psi - \lambda) < 0$, for the local algorithm. The latter rate holds under the mild assumption on $L$ of quadratic growth away from the minimizer. As mentioned in Section \ref{sec:Background}, Nesterov-like algorithms do not necessarily assure UGAS of the minimizer. The HAND-1 algorithm for convex $L$, proposed in \cite{poveda2019inducing}, provides UGAS via a hybrid restarting mechanism that yields a convergence rate $\frac{1}{t^2}$. However, this convergence rate holds only until the first reset. The algorithm we propose not only renders the minimizer UGAS, but also has a (hybrid) convergence rate that preserves the rates of the individual optimization algorithms for all (hybrid) time. Moreover, the global rate of our algorithm is commensurate with that of HAND-1. In \IfConf{Fig. \ref{fig:NSCLargerZeta}}{Figure \ref{fig:NSCLargerZetaTR}} and Section \ref{sec:Examples}, our uniting algorithm is shown via numerical simulations\footnote{Code at \texttt{gitHub.com/HybridSystemsLab/UnitingTradeoff} \label{foot:Tradeoff_TR}} to have improved performance over the HAND-1 algorithm in \cite{poveda2019inducing}. 
	
\begin{figure}[thpb]
	\centering
	\setlength{\unitlength}{1.0pc} 
	%
	%
	\begin{picture}(24,17.5)(0,0)
		\footnotesize
%
%
		\put(0,0.5){\includegraphics[trim={0.5cm 0.2cm 1cm 0.3cm},clip,width=24\unitlength]{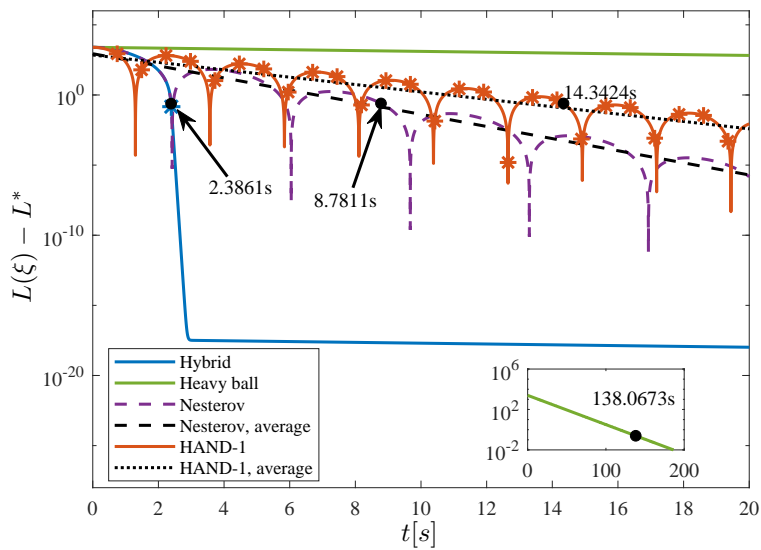}}
		\put(0,8){\rotatebox{90}{$L(\xi) - L^*$}}
		\put(12.3,0.5){$t[s]$}
	\end{picture}
	\caption{A comparison of the evolution of $L$ over time for Nesterov's method in \eqref{eqn:MJODE_ZetaNum_TR}, heavy ball, HAND-1 from \cite{poveda2019inducing}, and our proposed uniting algorithm, for a function $L(\xi) := \xi^2$, with a single minimizer at $\xi^* = 0$. Nesterov's method, shown in purple, settles to within $1\%$ of $\xi^*$ in about 8.8 seconds. The heavy ball algorithm, shown in green, settles to within $1\%$ of $\xi^*$ in about 138.1 seconds. HAND-1, shown in orange, settles to within $1\%$ of $\xi^*$ in about 14.3 seconds. The hybrid closed-loop system $\HS$, shown in blue, settles to within $1\%$ of $\minSet$ in about 2.4 seconds. As opposed to Figure \ref{fig:MotivationalPlotTR}, which uses $\zeta = 2$ for $\HS_1$, this example uses $\zeta = 1$, which results in slower convergence of solutions to $\HS$ and $\HS_1$ than in Figure \ref{fig:MotivationalPlotTR}.}	
	\label{fig:NSCLargerZetaTR}
\end{figure}
	
	\item {\em Extension of the results on Nesterov's method in \cite{muehlebach2019dynamical}:} In the process, in Section \ref{sec:ProofMainResult}, we extend the properties and convergence results for Nesterov's method in \cite{muehlebach2019dynamical}. In particular, while the convergence rate results in \cite{muehlebach2019dynamical} assume that $L(\xi_1^*) = 0$ at $\xi^* = 0$, and $\zeta = 1$ for \eqref{eqn:MJODENCVX}, 
	here we prove UGAS of the minimizer, with a convergence rate of $\frac{1}{(t+2)^2}$, for cost functions with a minimum value that is not necessarily zero, which holds for a generic parameter $\zeta > 0$. We achieve the relaxation on $\zeta$ by moving it into the numerator of the coefficient of the gradient, effectively decoupling $\zeta$ and $M$, namely,
	\begin{equation} \label{eqn:MJODE_ZetaNum_TR}
		\ddot{\xi} + 2\bar{d}(t)\dot{\xi} + \frac{\zeta^2}{M}\nabla L(\xi + \bar{\beta}(t) \dot{\xi}) = 0.
	\end{equation}
	Such a modification leads to faster convergence as $\zeta$ increases, and slower convergence as $\zeta \rightarrow 0$.
\end{enumerate}}
\IfConf{
	
	\vspace{-0.2cm}
	
}{}
\IfConf{While preliminary work in \cite{hustig2021uniting} proposed an algorithm uniting Nesterov's method globally and heavy ball locally for $\mathcal{C}^2$, strongly convex $L$, with different results and examples, the uniting algorithm proposed in this paper relaxes the conditions in \cite{hustig2021uniting} to $\mathcal{C}^1$, convex $L$ with a unique minimizer.}{Preliminary work in \cite{hustig2021uniting} proposed an algorithm uniting Nesterov's method globally and heavy ball locally for $\mathcal{C}^2$, strongly convex objective functions $L$, and included different results and examples that reflect such conditions, with proofs omitted due to space considerations. The uniting algorithm proposed in this paper relaxes the conditions in \cite{hustig2021uniting} to $\mathcal{C}^1$, convex $L$ with a unique minimizer. Such a relaxation is reflected in the results, examples, and proofs presented here.} \IfConf{A technical report version of this paper, with more details \cite{dhustigs2022unitingNSC}, is available online.}{}

\subsection{Notation}
The sets of real, positive real, and natural numbers are denoted by $\reals$, $\reals_{>0}$, and $\naturals$, respectively. The closed unit ball, of appropriate dimension, in the Euclidean norm is denoted as $\ball$. The set $\mathcal{C}^n$ represents the family of $n$-th continuously differentiable functions. For \IfConf{a vector}{vectors} $v \in \reals^n$\IfConf{,}{ and $w \in \reals^n$,} $\left|v\right| = \sqrt{v^{\top} v}$ denotes the Euclidean vector norm of $v$\IfConf{.}{, and $\langle v,w \rangle = v^{\top} w$ the inner product of $v$ and $w$.} For any $x \in \reals^n$ and $y \in \reals^m$, $\left(x,y\right) := [x^{\top},y^{\top}]^{\top}$. The closure of a set $S$ is denoted $\overline{S}$\IfConf{.}{ and the set of interior points of $S$ is denoted $\mathrm{int}(S)$.} Given a set $S \subset \reals^n \times \reals^m$, the projection of $S$ onto $\reals^n$ is defined as $\Pi(S) := \defset{x \in \reals^n\!\!}{\!\!\exists y \ \text{such that } (x,y) \in S}$. \IfConf{}{The distance of a point $x \in \reals^n$ to a set $S \in \reals^n$ is defined by $\left|x\right|_S = \inf_{y \in S} \left| y - x \right|$. }Given a set-valued mapping $M : \reals^m \rightrightarrows \reals^n$, the domain of $M$ is the set $\mathrm{dom} M = \defset{x \in \reals^m\!\!}{\!\!M(x) \neq \emptyset\!\!}$\IfConf{.}{, and the range of $M$ is the set $\mathrm{rge} \ M = \defset{y \in \reals^n\!\!}{\!\!\exists x \in \reals^m \mbox{ such that } y\in M(x)\!\!}$. A function $\alpha : \reals_{\geq 0} \rightarrow \reals_{\geq 0}$ is a class-$\mathcal{K}_{\infty}$ function, also written $\alpha \in \mathcal{K}_{\infty}$, if $\alpha$ is zero at zero, continuous, strictly increasing, and unbounded. A function $\beta : \reals_{\geq 0} \times \reals_{\geq 0} \rightarrow \reals_{\geq 0}$ is a class-$\mathcal{KL}$ function, also written $\beta \in \mathcal{KL}$, if it is nondecreasing in its first argument, nonincreasing in its second argument, $\lim_{r \rightarrow 0^+} \beta(r,s) = 0$ for each $s \in \reals_{\geq 0}$, and $\lim_{s \rightarrow \infty} \beta(r,s) = 0$ for each $r \in \reals_{\geq 0}$.}
\section{Uniting Optimization Algorithm}
	\label{sec:NonstronglyConvexNest}
	
\subsection{Problem Statement} \label{sec:ProblemStatement}
As illustrated in \IfConf{Fig. \ref{fig:MotivationalPlot}}{Figure \ref{fig:MotivationalPlotTR}}, the performance of Nesterov's accelerated gradient descent commonly suffers from oscillations near the minimizer. This is also the case for the heavy ball method when $\lambda > 0$ is small. However, when $\lambda$ is large, the heavy ball method converges slowly, albeit without oscillations. In Section \ref{sec:Intro} we discussed how Nesterov's algorithm guarantees a rate of $\frac{1}{(t+2)^2}$ for convex $L$. We also discussed how the heavy ball algorithm guarantees a rate of $\frac{1}{t}$ for convex $L$, although it was demonstrated in \cite{sebbouh2020convergence} that the heavy ball algorithm converges exponentially for convex $L$ when such an objective function also has the property of quadratic growth away from its minimizer. We desire to attain the rate $\frac{1}{(t+2)^2}$ globally and an exponential rate locally,
while avoiding oscillations via the heavy ball algorithm with large $\lambda$. We state the problem to solve as follows: 
\IfConf{\begin{prob}}{\begin{problem}}\label{problem:ProbStatement}
	Given a scalar, real-valued, continuously differentiable, and convex objective function $L$ with a unique minimizer, 
	design an optimization algorithm that, without knowing the function $L$ or the location of its minimizer, has the minimizer UGAS, with a convergence rate of $\frac{1}{(t+2)^2}$ globally and an exponential convergence rate locally,
	and with robustness to arbitrarily small noise in measurements of $\nabla L$.
\IfConf{\end{prob}}{\end{problem}}

\subsection{Modeling}
	
	In this section, we present an algorithm that solves \IfConf{Problem \ref{problem:ProbStatement}}{Problem $(\star)$}. We interpret the ODEs in \eqref{eqn:HBF} and \IfConf{\eqref{eqn:MJODE_ZetaNum}}{\eqref{eqn:MJODE_ZetaNum_TR}} as control systems consisting of a plant and a control algorithm \cite{191} \cite{220}. Defining $\xp_1$ as $\xi$ and $\xp_2$ as $\dot{\xi}$, the plant associated to these ODEs is given by 
	the double integrator
	\begin{equation}\label{eqn:HBFplant-dynamicsTR}
		\matt{\dot{\xp}_1\\ \dot{\xp}_2\\} = \matt{\xp_2\\ \up} =: \fp(\xp,u) \ \ \  (\xp,\up) \in \reals^{2n} \times \reals^n =: \Cp
	\end{equation}
	With this model, the optimization algorithms that we consider assign $u$ to a function of the state that involves the cost function, and such a function of the state may be time dependent. The control algorithm leading to \eqref{eqn:HBF} assigns $u$ to $-\lambda \xp_2 - \gamma \nabla L(\xp_1)$
	 where $\gamma > 0$ and $\lambda > 0$, and the control algorithm leading to \IfConf{\eqref{eqn:MJODE_ZetaNum}}{\eqref{eqn:MJODE_ZetaNum_TR}} assigns $u$ to $-2\bar{d}(t)\xp_2 - \frac{\zeta^2}{M}\nabla L(\xp_1 + \bar{\beta}(t)\xp_2)$
	where $\zeta > 0$, $M >0$ is the Lipschitz constant for $\nabla L$, and where $\bar{d}(t)$ and $\bar{\beta}(t)$ are defined, for all $t \geq 0$, as
%
%
	\begin{equation} \label{eqn:dBarBetaBar}
		\bar{d}(t) := \frac{3}{2(t+2)}, \quad \bar{\beta}(t) := \frac{t - 1}{t + 2}.
	\end{equation}
	The functions $\bar{d}$ and $\bar{\beta}$ are defined expressly as in \eqref{eqn:dBarBetaBar} for ease of analysis in the forthcoming Propositions \ref{prop:ConvergenceNSCVXNesterov}-\ref{prop:UGASNest}. Such a time-varying definition satisfies the linear combination of the damping terms mentioned below \eqref{eqn:MJODENCVX}. While constant terms can be used for \eqref{eqn:MJODENCVX}, when $L$ is strongly convex, constant damping terms are not adequate for convex $L$; see \cite{muehlebach2019dynamical}.
	The proposed logic-based algorithm ``unites'' the two optimization algorithms modeled by $\hcq$, where the logic variable $\xlogic \in Q := \{0,1\}$ indicates which algorithm is currently being used. The local and global algorithms, respectively, are defined as
	\begin{subequations} \label{eqn:StaticStateFeedbackLawsNSC}
		\begin{align}
			\hczero(\hp_0(\xp)) & = -\lambda \xp_2 - \gamma \nabla L(\xp_1) \label{eqn:StaticStateFeedbackLawLocal}\\
			\hcone(\hp_1(\xp,t),t) & = -2\bar{d}(t)\xp_2 - \frac{\zeta^2}{M}\nabla L(\xp_1 + \bar{\beta}(t)\xp_2)\label{eqn:StaticStateFeedbackLawGlobalNSCVX}
		\end{align}
	\end{subequations}
	where the algorithm defined by $\hcone$ plays the role of the global algorithm in uniting control (see, e.g., {\cite[Chapter~4]{220}), while the algorithm defined by $\hczero$ plays the role of the local algorithm.
	The outputs $\hp_0$ corresponding to the output for the heavy ball algorithm and $\hp_1$ corresponding to the output for Nesterov's algorithm are defined as
	\begin{equation} \label{eqn:H0H1NSCNesterovHBF}
		\hp_0(\xp)\! :=\! \matt{\xp_2\\ \nabla L(\xp_1)}\!, \hp_1(\xp,t)\! :=\! \matt{\xp_2\\ \nabla L(\xp_1 + \bar{\beta}(t) \xp_2)}.
	\end{equation}
	Namely, the algorithm exploits measurements of $\nabla L$, which in practice are typically approximated using measurements of $L$. The parameters $\lambda > 0$ and $\gamma > 0$ should be designed to achieve convergence without oscillations nearby the minimizer.
	
	We use the hybrid systems framework to design our algorithm. A hybrid system $\HS$ has data $(C,F,D,G)$ and is defined as \cite[Definition~2.2]{65} \IfConf{
			
			\vspace{-0.55cm}
			
		}{}
		
		\begin{equation}\label{eqn:GeneralH}
			\HS = 
			\begin{cases}
				\dot{x} & \!\!\!\!\! \in \f(x) \quad x \in C \\
				x^+ & \!\!\!\!\! \in G(x) \quad x \in D
			\end{cases}
		\end{equation}
		where $x \in \reals^n$ is the system state, $\f : \reals^n \rightrightarrows \reals^n$ is the flow map, $C \subset \reals^n$ is the flow set, $\g : \reals^n \rightrightarrows \reals^n$ is the jump map, and $D \subset \reals^n$ is the jump set. Since the ODE in \IfConf{\eqref{eqn:MJODE_ZetaNum}}{\eqref{eqn:MJODE_ZetaNum_TR}} is time varying, and since solutions to hybrid systems are parameterized by\footnote{The variable $t$ is the amount of time that has passed and $j$ is the number of jumps that have occurred.} $(t,j) \in \reals_{\geq 0} \times \naturals$, we employ the state $\tau$ to capture ordinary time as a state variable, in this way, leading to a time-invariant hybrid system. To encapsulate the plant, static state-feedback laws, and the time-varying nature of the ODE in \IfConf{\eqref{eqn:MJODE_ZetaNum}}{\eqref{eqn:MJODE_ZetaNum_TR}}, we define a hybrid closed-loop system $\HS$ with state $x := \left(\xp, \xlogic, \tau\right) \in \reals^{2n} \times \xlogicSpace \times \reals_{\geq 0}$ \IfConf{as}{as\footnote{Although 
		$h_0$ in \eqref{eqn:H0H1NSCNesterovHBF} does not depend on $\tau$, to simplify notation we keep $\tau$ as an argument in the general function $h_{\xlogic}$.}}
%
%
	\begin{subequations} \label{eqn:HS-TimeVarying}
		\begin{equation}
			\left.\begin{aligned} 
				\dot{\xp} & = \matt{\xp_2 \\ \hcq(\hp_{\xlogic}(\xp,\tau),\tau)} \label{eqn:FlowMap}\\
				\dot{\xlogic} & = 0 \\
				\dot{\tau} & = \xlogic
			\end{aligned}\right\} =: F(x) \qquad x \in C := C_0 \cup C_1
		\end{equation}
		\begin{equation}
			\left.\begin{aligned}
				\xp^+ & = \matt{\xp_1 \\ \xp_2}\\
				\xlogic^+ & = 1-\xlogic\\
				\tau^+ & = 0
			\end{aligned}\right\} =: G(x) \qquad x \in D := D_0 \cup D_1 \label{eqn:JumpMap}
		\end{equation}
	\end{subequations}
	The sets $C_0$, $C_1$, $D_0$, and $D_1$ are defined as
	\begin{subequations} \label{eqn:CAndDGradientsNestNSC}
		\begin{align}
			&C_0 := \mathcal{U}_0 \times \{0\} \times \{0\}, \ C_1 := \overline{\reals^{2n}\setminus \T_{1,0}} \times \{1\} \times \reals_{\geq 0}\\
			&D_0 := \T_{0,1} \times \{0\} \times \{0\}, \ D_1 := \T_{1,0}\times \{1\} \times \reals_{\geq 0}.
		\end{align}
	\end{subequations} 
%
%
The sets $\mathcal{U}_0$, $\T_{1,0}$, and $\T_{0,1}$ are precisely defined in Section \ref{sec:AssDesign}, using Lyapunov functions defined therein, but the idea behind their construction is as follows. The switch between $\hczero$ and $\hcone$ is governed by a {\em supervisory algorithm} implementing switching logic\IfConf{.}{; see Figure \ref{fig:FeedbackDiagramTR}.} The supervisor selects between these two optimization algorithms, based on the output of the plant in \eqref{eqn:H0H1NSCNesterovHBF} and the optimization algorithm currently applied. When $\xp \in \mathcal{U}_0$, $\xlogic = 0$, and $\tau = 0$ (i.e., $x \in C_0$), due to the design of $\mathcal{U}_0$ in Section \ref{sec:U0}, then the state $\xp$ is near the minimizer, which is denoted $\xp_1^*$, and the supervisor allows flows of \eqref{eqn:HS-TimeVarying} using $\hczero$ and $\dot{\tau} = q = 0$ to avoid oscillations. Conversely, when $\xp \in \overline{\reals^{2n}\setminus \T_{1,0}}$ and $\xlogic = 1$ (i.e., $x \in C_1$), due to the design of $\T_{1,0}$ in Section \ref{sec:DesignT10}, then the state $\xp$ is far from the minimizer and the supervisor allows flows of \eqref{eqn:HS-TimeVarying} using $\hcone$ and $\dot{\tau} = q = 1$ to converge quickly to the neighborhood of the minimizer. When $\xp \in \T_{1,0}$ and $\xlogic = 1$ (i.e., $x \in D_1$), then this indicates that the state $\xp$ is near the minimizer, and the supervisor assigns $u$ to $\hczero$, resets $\xlogic$ to $0$, and resets $\tau$ to $0$. Conversely, when $\xp \in \T_{0,1}$, $\xlogic = 0$, and $\tau = 0$ (i.e., $x \in D_0$), due to the design of $\T_{0,1}$ in Section \ref{sec:T01}, then this indicates that the state $\xp$ is far from the minimizer and the supervisor assigns $u$ to $\hcone$ and resets $\xlogic$ to $1$. \IfConf{}{
	The complete algorithm, defined in \eqref{eqn:HS-TimeVarying}-\eqref{eqn:CAndDGradientsNestNSC}, is summarized in Algorithm \ref{algo:HS-Uniting}.
	\begin{algorithm}[thpb]
		\caption{Uniting algorithm}
		\label{algo:HS-Uniting}
		\begin{algorithmic}[1]
			\IfConf{\State}{\STATE} Set $\xlogic(0,0)$ to $0$, $\tau(0,0)$ to $0$, and set $\xp(0,0)$ as an initial condition with an arbitrary value.
			\IfConf{\While{true}}{\WHILE{true}}
			\IfConf{\If{$\xp \in \T_{0,1}$, $\xlogic = 0$, and $\tau = 0$}}{\IF{$\xp \in \T_{0,1}$, $\xlogic = 0$, and $\tau = 0$}}
			\IfConf{\State}{\STATE} Reset $\xlogic$ to $1$.
			\IfConf{\ElsIf{$\xp \in \T_{1,0}$ and $\xlogic = 1$}}{\ELSIF{$\xp \in \T_{1,0}$ and $\xlogic = 1$}}
			\IfConf{\State}{\STATE} Reset $\xlogic$ to $0$ and $\tau$ to $0$.
			\IfConf{\ElsIf{$\xp \in \mathcal{U}_0$, $\xlogic = 0$, and $\tau = 0$}}{\ELSIF{$\xp \in \mathcal{U}_0$, $\xlogic = 0$, and $\tau = 0$}}
			\IfConf{\State}{\STATE} Assign $u$ to $\hczero(\hp_0(\xp))$ and update $\xp$, $\xlogic$, and $\tau$ according to \eqref{eqn:FlowMap}.
			\IfConf{\ElsIf{$\xp \in \overline{\reals^{2n}\setminus \T_{1,0}}$ and $\xlogic = 1$}}{\ELSIF{$\xp \in \overline{\reals^{2n}\setminus \T_{1,0}}$ and $\xlogic = 1$}}
			\IfConf{\State}{\STATE} Assign $u$ to $\hcone(\hp_1(\xp,\tau),\tau)$ and update $\xp$, $\xlogic$, and $\tau$ according to \eqref{eqn:FlowMap}.
			\IfConf{\EndIf}{\ENDIF}
			\IfConf{\EndWhile}{\ENDWHILE}
		\end{algorithmic}
	\end{algorithm} }

	The reason that the state $\tau$ in \eqref{eqn:HS-TimeVarying} changes at the rate $\xlogic$ during flows and is reset to $0$ at jumps is that when the state $x$ is in $C_1$, then $\dot{\tau} = \xlogic = 1$, which implies that $\tau$ behaves as ordinary time, so it is used to represent time in the time-varying algorithm $\hcone$. On the other hand, when the state $x$ is in $C_0$, then $\dot{\tau} = \xlogic = 0$ causes the state $\tau$ to stay at zero, which is an appropriate value for $\tau$ as it is not required by the time-invariant algorithm $\hczero$. Such an evolution ensures that the set to asymptotically stabilize is compact.
	
	\IfConf{}{\IfConf{Fig. \ref{fig:FeedbackDiagram}}{Figure \ref{fig:FeedbackDiagramTR}} shows the feedback diagram of this hybrid closed-loop system $\HS$.} We denote the closed-loop system resulting from $\hczero$ as $\HS_0$, which is given by \IfConf{
		
		\vspace{-0.5cm}	
		
	}{}
	\begin{equation}\label{eqn:H0}
		\dot{\xp} = \matt{\xp_2 \\ \hczero(\hp_0(\xp))} \qquad \xp \in \reals^{2n}
	\end{equation}
	and we denote the closed-loop system resulting from $\hcone$ as $\HS_1$, which is given by \IfConf{

\vspace{-0.9cm}	

}{}
\begin{equation} \label{eqn:H1}
	\dot{\xp} = \matt{\xp_2 \\ \hcone(\hp_1(\xp,\tau),\tau)}, \ \dot{\tau} = 1 \qquad \left(\xp,\tau\right) \in \reals^{2n} \times \reals_{\geq 0}.
\end{equation}

\IfConf{
%
%
%
%
%
}{
\begin{figure}[thpb]
	\centering
	\setlength{\unitlength}{1.0pc}
	
	\begin{picture}(30.8,12)(0,0) 
		\footnotesize
		%
		%
		\put(0,2){\includegraphics[scale=0.2,trim={1.5cm 0cm 0cm 0.8cm},clip,width=11.5\unitlength]{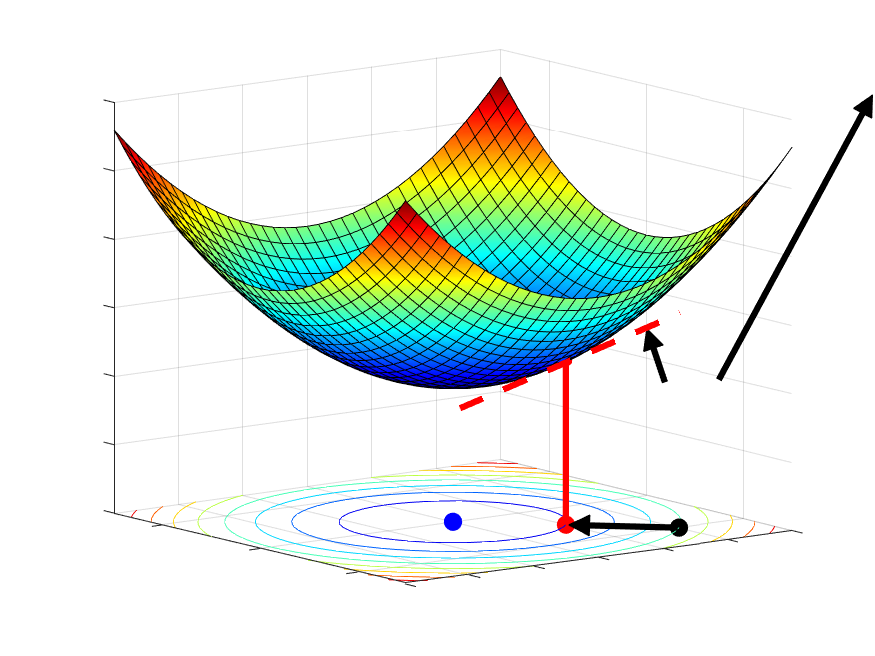}}			
		\put(11.5,0.5){\includegraphics[scale=0.4,trim={0cm 0cm 0cm 0cm},clip,width=19\unitlength]{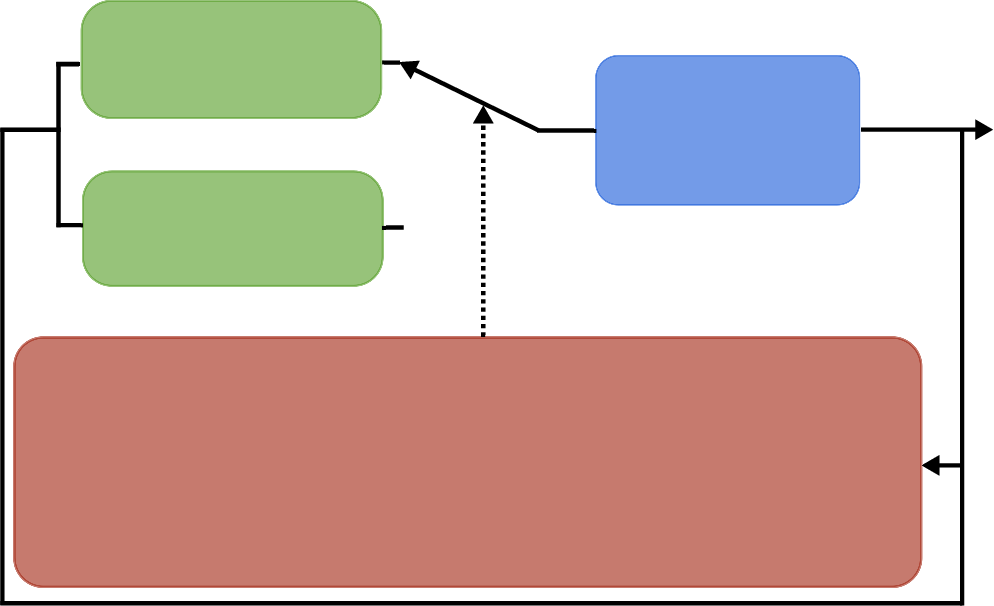}}
		\put(18.5,4.5){{\bf Supervisor}}
		\put(12.1,3.2){$\ \ \dot{\xlogic} = 0 \qquad \quad \ \dot{\tau} = \xlogic \ \ \ (\xp, \xlogic, \tau) \in C := C_0 \cup C_1$}
		\put(12.1,1.8){$\xlogic^+ = 1- \xlogic \ \ \  \tau^+ = 0 \ \ \ (\xp, \xlogic, \tau) \in D := D_0 \cup D_1$}
		\put(24.4,7.4){\bf plant}
		\put(24.1,9.9){$\dot{\xp}_1 = \xp_2$}
		\put(24.1,8.9){$\dot{\xp}_2 = u$}
		\put(13.5,11.2){$\hcone(\hp_1(\xp,\tau),\tau)$}
		\put(13.5,10.2){$\dot{\tau} = 1, \tau^+ = 0$}
		\put(13.3,9.1){\bf global ($\xlogic = 1$)}
		\put(14.3,7.5){$\hczero(\hp_0(\xp))$}
		\put(13.6,5.9){\bf local ($\xlogic = 0$)}
		\put(21,7){$\xlogic$}
		\put(22.1,10.1){$u$}
		\put(28.4,10.1){$\hp_{\xlogic}$}
		\put(11.6,10.1){$\hp_{\xlogic}$}
		\put(7.6,5.3){$\nabla L(\xp_1)$}
		\put(5.5,3.8){$\xp_1^*$}
		\put(6.7,3.3){$\xp_1$}
		\put(8.8,4){$\xp_{1_{\circ}}$}
	\end{picture}
	\caption{Feedback diagram of the hybrid closed-loop system $\HS$ (on the right), in \eqref{eqn:HS-TimeVarying}, uniting global and local optimization algorithms. An example optimization problem to solve is shown on the left and, for this example optimization problem, measurements of the gradient are used for the input of $\hcq$.}
	\label{fig:FeedbackDiagramTR}
\end{figure}}
	
\subsection{Design of the Hybrid Algorithm} \label{sec:AssDesign}
In order for the supervisor to determine when the state component $\xp_1$ is close to the minimizer of $L$, denoted $\minSet$, without knowledge of $\minSet$ or $L^* := L(\xp_1)$, we impose the following assumptions on $L$. 
\IfConf{\begin{assum}}{\begin{assumption}} \label{ass:LisNSCVX}
	The function $L$ is $\mathcal{C}^1$, convex\footnote{A function $L \!:\! \reals^n \rightarrow \reals$ is convex if, for all $u_1, w_1 \in \reals^n$, $L(u_1) \geq L(w_1) + \left\langle \nabla L(w_1), u_1 - w_1 \right\rangle $ \cite{boyd2004convex}. \label{foot:Convexity}}, and has a single minimizer $\xp_1^*$. 
\IfConf{\end{assum}}{\end{assumption}}
\IfConf{\begin{assum}[Quadratic growth of $L$]}{\begin{assumption}[Quadratic growth of $L$]} \label{ass:QuadraticGrowth}
		The function $L$ has quadratic growth away from its minimizer $\minSet$; i.e., there exists $\alpha > 0$ such that \IfConf{$L(\xp_1) - L^* \geq \alpha \left| \xp_1 - \xp_1^* \right|^2$ for all $\xp_1 \in \reals^{n}$.}{
		\begin{equation} \label{eqn:QuadraticGrowth}
			L(\xp_1) - L^* \geq \alpha \left| \xp_1 - \xp_1^* \right|^2 \quad \forall \xp_1 \in \reals^{n}.
		\end{equation}}
\IfConf{\end{assum}}{\end{assumption}}
\IfConf{}{\begin{remark}
	Assumption \ref{ass:LisNSCVX}, which is a common assumption used in the analysis of optimization algorithms \cite{boyd2004convex} \cite{nesterov2004introductory}, ensures that the objective function is continuously differentiable, which is necessary for well-posedness of $\HS$, as was explained in Section \ref{sec:Contributions}. Additionally, the convex property and the restriction that $L$ has a single minimizer $\xp_1^*$ in Assumption \ref{ass:LisNSCVX} rules out the possibility of the objective function having a continuum of minimizers or multiple isolated minimizers.
	Assumption \ref{ass:QuadraticGrowth}, which is used for the construction of $\mathcal{U}_0$, $\T_{1,0}$, and $\T_{0,1}$, is employed as a means of determining when the state $\xp$ is near the minimizer of $L$, via measurements of \IfConf{$\nabla L$}{the gradient}. Such an assumption is also commonly used in the analysis of convex optimization algorithms; see, e.g., \cite{drusvyatskiy2018error}, \cite{karimi2016linear}.
\end{remark}}

To make the switch back to $\hcone$, we impose the following assumption on $L$. 
\IfConf{\begin{assum}[Lipschitz Continuity of $\nabla L$]}{\begin{assumption}[Lipschitz Continuity of $\nabla L$]} \label{ass:Lipschitz}
	The function $\nabla L$ is Lipschitz continuous with constant $M > 0$, namely,  \IfConf{$\left|\nabla L(w_1) - \nabla L(u_1)\right| \leq M \left|w_1 - u_1\right|$ for all $w_1, u_1 \in \reals^n$.}{\\ $\left|\nabla L(w_1) - \nabla L(u_1)\right| \leq M \left|w_1 - u_1\right|$ for all $w_1, u_1 \in \reals^n$.}
\IfConf{\end{assum}}{\end{assumption}}
	\IfConf{}{\begin{remark}
		Assumption \ref{ass:Lipschitz} is used in the forthcoming construction of $\T_{0,1}$. Additionally, Assumption \ref{ass:Lipschitz} is commonly used in nonlinear analysis to ensure that the differential equations of the individual optimization algorithms, for example, those in \IfConf{\eqref{eqn:MJODE_ZetaNum}}{\eqref{eqn:MJODE_ZetaNum_TR}} and \eqref{eqn:HBF}, do not have solutions that escape in finite time, which is used to guarantee existence and completeness of maximal solutions to $\HS_{\xlogic}$ \cite[Theorem~3.2]{khalil2002nonlinear}.
	\end{remark}}

Under Assumptions \ref{ass:LisNSCVX} and \ref{ass:QuadraticGrowth}, the following lemma, used in some of the results to follow, relates the size of the gradient at a point to the distance from the point to $\minSet$. \IfConf{Its proof is in \cite{dhustigs2022unitingNSC}.}{}
\IfConf{\begin{lem}[Suboptimality]}{\begin{lemma}(Suboptimality):} \label{lemma:ConvexSuboptimality}
	Let $L$ satisfy Assumptions \ref{ass:LisNSCVX} and \ref{ass:QuadraticGrowth}, and let $\alpha > 0$ come from Assumption \ref{ass:QuadraticGrowth}. For some $\varepsilon > 0$, if $\xp_1 \in \reals^{n}$ is such that
	$\left| \nabla L(\xp_1) \right| \leq \varepsilon \alpha$, then $\left| \xp_1 - \xp_1^* \right| \leq \varepsilon$.
\IfConf{\end{lem}}{\end{lemma}} 
%
%
	
\IfConf{}{\begin{proof}
	Combining Assumption \ref{ass:LisNSCVX} and\IfConf{}{ \eqref{eqn:QuadraticGrowth} from} Assumption \ref{ass:QuadraticGrowth} with $u_1 = \xp_1^*$ and $w_1 = \xp_1$ yields
	\begin{align} \label{eqn:NearOptimalityC}
		\IfConf{\alpha \left| \xp_1 - \xp_1^* \right|^2 & \leq \left| L(\xp_1) - L^* \right| \leq \left| \left\langle \nabla L(\xp_1), \xp_1^* - \xp_1 \right\rangle \right|\nonumber\\ & \leq \left| \nabla L(\xp_1) \right| \left| \xp_1 - \xp_1^* \right|}{\alpha \left| \xp_1 - \xp_1^* \right|^2 \leq \left| L(\xp_1) - L^* \right| \leq \left| \left\langle \nabla L(\xp_1), \xp_1^* - \xp_1 \right\rangle \right| \leq \left| \nabla L(\xp_1) \right| \left| \xp_1 - \xp_1^* \right|}
	\end{align}
	where the first inequality holds since $L(\xp_1) \geq L^*$.
	Then, \IfConf{$\left| \xp_1 - \xp_1^* \right| \leq \frac{1}{\alpha} \left| \nabla L(\xp_1) \right|$, from which}{
	\begin{equation}\label{eqn:Suboptimality}
		\left| \xp_1 - \xp_1^* \right| \leq \frac{1}{\alpha} \left| \nabla L(\xp_1) \right|.
	\end{equation}
	From \eqref{eqn:Suboptimality},} we can deduce that $\left| \nabla L(\xp_1) \right| \leq \varepsilon \alpha$ implies $\left| \xp_1 - \xp_1^* \right| \leq \frac{1}{\alpha} \left(\varepsilon \alpha \right) = \varepsilon$.\IfConf{\hfill{} \qed}{}
\end{proof}} 
%
%
The suboptimality condition from Lemma \ref{lemma:ConvexSuboptimality} is typically used as a stopping condition for optimization, as it indicates that the argument of $L$ is close enough to the minimizer $\minSet$ \cite{boyd2004convex}. We exploit Lemma \ref{lemma:ConvexSuboptimality} to determine when the state component $\xp_1$ of the hybrid closed-loop system $\HS$ is close enough to the minimizer $\minSet$ so as to switch to the local optimization algorithm, $\hczero$, in this way activating $\HS_0$\IfConf{}{; see \IfConf{Fig. \ref{fig:FeedbackDiagram}}{Figure \ref{fig:FeedbackDiagramTR}}}. \IfConf{
	
	\vspace{-0.1cm}
	
}{}
\subsubsection{Design of the Set $\mathcal{U}_0$} \label{sec:U0}
\IfConf{T}{Recall from lines 7-8 of Algorithm \ref{algo:HS-Uniting} that t}he objective is to design $\mathcal{U}_0$ such that when $\xp \in \mathcal{U}_0$, $\xlogic = 0$, and $\tau = 0$, the state component $\xp_1$ is near $\minSet$ and the uniting algorithm allows flows of \eqref{eqn:HS-TimeVarying} with $\hczero$ and $\xlogic = 0$. For such a design, we use Assumptions \ref{ass:LisNSCVX} and \ref{ass:QuadraticGrowth} and the Lyapunov function\IfConf{\vspace{-0.2cm}}{}
	\begin{equation} \label{eqn:LyapunovHBF}
		V_0(\xp) := \gamma \left(L(\xp_1) - L^* \right) + \frac{1}{2} \left|\xp_2\right|^2 
	\end{equation}
defined for each $\xp \in \reals^{2n}$, where $\gamma > 0$. The choice of $V_0$ in \eqref{eqn:LyapunovHBF} \IfConf{is used in the proof of the forthcoming Proposition~\ref{prop:GAS-HBF} to establish UGAS of the minimizer for $\HS_0$ in \eqref{eqn:H0}.}{leads to $\dot{V}_0$ showing decrease of $V_1$ for each $\xp \in \reals^{2n}$, in the proof of the forthcoming Proposition~\ref{prop:GAS-HBF}. Such a property of $V_0$ is needed to establish UGAS of the minimizer for $\HS_0$ in \eqref{eqn:H0}.}
Given $\varepsilon_0 > 0$, $c_0 > 0$, and $\gamma > 0$ from $\hczero$ in \eqref{eqn:StaticStateFeedbackLawLocal}, let $\alpha > 0$ come from Assumption \ref{ass:QuadraticGrowth} such that
	\begin{equation} \label{eqn:UTilde0SetEquations}
		\tilde{c}_0 := \varepsilon_0 \alpha > 0, \quad
		d_0 := c_0 - \gamma \left(\frac{\tilde{c}_0^2}{\alpha}\right) > 0.
	\end{equation}
\IfConf{Then, $V_0$ in \eqref{eqn:LyapunovHBF} can be upper bounded as follows:}{Then, $V_0$ in \eqref{eqn:LyapunovHBF} can be upper bounded, using Assumption \ref{ass:LisNSCVX} as done to arrive to \eqref{eqn:NearOptimalityC}, as follows: for each $\xp \in \reals^{2n}$ 
\begin{align} \label{eqn:c0SublevelSet}
	\IfConf{V_0(\xp) = \gamma \left(L(\xp_1) - L^*\right) + \frac{1}{2}\left|\xp_2\right|^2 \leq & \gamma \left|\nabla L(\xp_1)\right| \left| \xp_1 - \xp_1^* \right|\nonumber\\ & + \frac{1}{2}\left|\xp_2\right|^2.}{V_0(\xp) = \gamma \left(L(\xp_1) - L^*\right) + \frac{1}{2}\left|\xp_2\right|^2 \leq \gamma \left|\nabla L(\xp_1)\right| \left| \xp_1 - \xp_1^* \right| + \frac{1}{2}\left|\xp_2\right|^2.}
\end{align}
Then,}
due to $L$ being $\mathcal{C}^1$, convex, and having a single minimizer $\xp_1^*$ by Assumption \ref{ass:LisNSCVX}, and due to $L$ having quadratic growth away from $\minSet$ by Assumption \ref{ass:QuadraticGrowth}, when $\left|\nabla L(\xp_1)\right| \leq \tilde{c}_0$, the suboptimality condition in Lemma \ref{lemma:ConvexSuboptimality} implies $\left| \xp_1 - \xp_1^* \right| \leq \frac{\tilde{c}_0}{\alpha}$, from where we get\IfConf{\vspace{-0.2cm}}{}
\begin{equation} \label{eqn:V0Bound}
	V_0(\xp) \leq \gamma \left(\frac{\tilde{c}_0^2}{\alpha}\right) + \frac{1}{2}\left|\xp_2\right|^2
\end{equation}
Then, by defining the set $\mathcal{U}_0$ as
\begin{equation}\label{eqn:U0}
	\mathcal{U}_0 := \defset{\xp \in \reals^{2n}\!\!}{\!\!\left| \nabla L(\xp_1) \right| \leq \tilde{c}_0,\frac{1}{2} \left|\xp_2\right|^2 \leq d_0\!\!},
\end{equation}
every $\xp \in \mathcal{U}_0$ belongs to the $c_0$-sublevel set of $V_0$. In fact, using the conditions in \eqref{eqn:UTilde0SetEquations} and \eqref{eqn:V0Bound}, we have that for each $\xp \in \mathcal{U}_0$, \IfConf{$V_0(\xp) \leq \gamma \left(\frac{\tilde{c}_0^2}{\alpha}\right) + \frac{1}{2}\left|\xp_2\right|^2 \leq c_0$.}{
\begin{equation} \label{eqn:V0c0SublevelSet}
	V_0(\xp) \leq \gamma \left(\frac{\tilde{c}_0^2}{\alpha}\right) + \frac{1}{2}\left|\xp_2\right|^2 \leq c_0.
\end{equation}}
Since $\hczero$ in \eqref{eqn:StaticStateFeedbackLawLocal} is such that the set $\{\minSet\} \times \{0\}$ is globally asymptotically stable for the closed-loop system resulting from controlling \eqref{eqn:HBFplant-dynamicsTR} by $\hczero$, as we show in the forthcoming Proposition \ref{prop:GAS-HBF}, the set $\mathcal{U}_0$ is contained in the basin of attraction induced by $\hczero$. 

\subsubsection{Design of the Set $\T_{1,0}$}
\label{sec:DesignT10}
\IfConf{T}{Recall from lines 5-6 of Algorithm \ref{algo:HS-Uniting} that t}he objective is to design $\T_{1,0}$ such that when $\xp \in \T_{1,0}$ and $\xlogic = 1$, the state component $\xp_1$ is near $\minSet$ and the supervisor resets $\xlogic$ to $0$, resets $\tau$ to $0$, and assigns $u$ to $\hczero(\hp_0(\xp))$. For such a design, we use Assumptions \ref{ass:LisNSCVX} and \ref{ass:QuadraticGrowth} and the Lyapunov function 
\IfConf{\vspace{-0.2cm}}{}
\begin{equation} \label{eqn:LyapunovNesterovNSCVX}
	V_1(\xp,\tau)\! :=\! \frac{1}{2}\left|\bar{a}(\tau)\left(\xp_1 - \xp_1^*\right) \!+\! \xp_2\right|^2 + \frac{\zeta^2}{M} (L(\xp_1) - L^*)\!
\end{equation}
defined for each $\xp \in \reals^{2n}$ and each $\tau \geq 0$, where $\zeta > 0$, $M > 0$ is the Lipschitz constant of $\nabla L$, and the function $\bar{a}$ is defined as
\begin{equation} \label{eqn:BarA}
	\bar{a}(\tau) := \frac{2}{\tau+2}.
\end{equation} 
The choice of $V_1$ in \eqref{eqn:LyapunovNesterovNSCVX} comes from \cite{muehlebach2019dynamical}, and \IfConf{is used to establish UGAS of the minimizer for $\HS_1$ in \eqref{eqn:H1} and the convergence rate $\frac{1}{\left(t+2\right)^2}$.}{such a choice leads to $\dot{V}_1$ showing decrease of $V_1$, for each $\xp \in \reals^{2n}$ and each $\tau \in \reals_{\geq 0}$, in the proof of the forthcoming Proposition \ref{prop:ConvergenceNSCVXNesterov}. Such a property of $V_1$ is needed to establish UGAS of the minimizer for $\HS_1$ in \eqref{eqn:H1}, and the convergence rate $\frac{1}{\left(t+2\right)^2}$.} In this same proof, the specific choice of $\bar{a}$ in \eqref{eqn:BarA}, which comes from \cite{muehlebach2019dynamical}, \IfConf{is also used to show decrease of $V_1$.}{leads to the elimination of the cross term $\left\langle \xp_1 - \xp^*_1,\xp_2 \right\rangle$ -- which has indeterminate sign -- in the upper bound on $\dot{V}_1$. In other words, without the definitions of $V_1$ in \eqref{eqn:LyapunovNesterovNSCVX} and $\bar{a}$ in \eqref{eqn:BarA}, decrease of $V_1$ and, consequently, UGAS of the minimizer for $\HS_1$and the convergence rate $\frac{1}{\left(t+2\right)^2}$ cannot be established.} More details on how \eqref{eqn:LyapunovNesterovNSCVX} and \eqref{eqn:BarA} are used in our analysis can be found in \IfConf{Proposition \ref{prop:ConvergenceNSCVXNesterov}.}{Appendix \ref{sec:ProofProp54}.} Given $c_{1,0} \in (0, c_0)$ and $\varepsilon_{1,0} \in (0,\varepsilon_0)$, where $c_0 > 0$ and $\varepsilon_0 > 0$ come from Section \ref{sec:U0}, let $\tilde{c}_0$ and $d_0$ be given in \eqref{eqn:UTilde0SetEquations}, and let $\alpha > 0$ come from Assumption \ref{ass:QuadraticGrowth} such that
\begin{subequations} \label{eqn:TTilde10SetEquations}
	\begin{align} 
		\tilde{c}_{1,0} & := \varepsilon_{1,0} \alpha \in (0, \tilde{c}_0) \\
		d_{1,0} & := c_{1,0} - \left(\frac{\tilde{c}_{1,0}}{\alpha}\right)^2 - \frac{\zeta^2}{M}\left(\frac{\tilde{c}_{1,0}^2}{\alpha}\right) \in (0, d_0)
	\end{align}
\end{subequations}
where $\zeta > 0$ comes from \IfConf{\eqref{eqn:MJODE_ZetaNum}}{\eqref{eqn:MJODE_ZetaNum_TR}}. Note that $\bar{a}$, defined via \eqref{eqn:BarA}, which is in $V_1$, equals $1$ when $\tau = 0$ and monotonically decreases toward zero (but being always positive) as $\tau$ tends to $\infty$. Namely, $\bar{a}$ is upper bounded by $1$.
Then, with $V_1$ given in \eqref{eqn:LyapunovNesterovNSCVX}\IfConf{,}{ and using Assumption \ref{ass:LisNSCVX} with $u_1 = \xp_1^*$ and $w_1 = \xp_1$, \IfConf{$V_1(\xp,\tau) \leq \left| \xp_1 - \xp_1^* \right|^2 + \left|\xp_2\right|^2 + \frac{\zeta^2}{M} \left| \nabla L(\xp_1) \right|\left| \xp_1 - \xp_1^* \right|$.}{
\begin{equation} \label{eqn:c10SublevelSetNestNSC}
	V_1(\xp,\tau) \leq \left| \xp_1 - \xp_1^* \right|^2 + \left|\xp_2\right|^2 + \frac{\zeta^2}{M} \left| \nabla L(\xp_1) \right|\left| \xp_1 - \xp_1^* \right|.
\end{equation}}
Then,} due to $L$ being $\mathcal{C}^1$, convex, and having a single minimizer $\xp_1^*$ by Assumption \ref{ass:LisNSCVX}, and due to $L$ having quadratic growth away from $\minSet$ by Assumption \ref{ass:QuadraticGrowth}, when $\left|\nabla L(\xp_1)\right| \leq \tilde{c}_{1,0}$, the suboptimality condition in Lemma \ref{lemma:ConvexSuboptimality} implies $\left| \xp_1 - \xp_1^* \right| \leq \frac{\tilde{c}_{1,0}}{\alpha}$, from where we get \IfConf{\vspace{-0.2cm}}{}
\begin{equation} \label{eqn:V1Bound}
	V_1(\xp,\tau) \leq \left(\frac{\tilde{c}_{1,0}}{\alpha}\right)^2 + \left|\xp_2\right|^2 + \frac{\zeta^2}{M} \left(\frac{\tilde{c}_{1,0}^2}{\alpha}\right). 
\end{equation}
Then, by defining $\T_{1,0}$ as
\begin{equation} \label{eqn:T10}
	\T_{1,0} := \defset{\xp \in \reals^{2n}\!\!}{\!\!\left| \nabla L(\xp_1) \right| \leq \tilde{c}_{1,0}, \left|\xp_2\right|^2 \leq d_{1,0}\!\!} 
\end{equation}
which, by construction, is contained in the interior of $\mathcal{U}_0$ defined in \eqref{eqn:U0}, every $\xp \in \T_{1,0}$ belongs to the $c_{1,0}$-sublevel set of $V_1$. In fact, using the conditions in \eqref{eqn:TTilde10SetEquations} and \eqref{eqn:V1Bound}, we have for each $\xp \in \T_{1,0}$, \IfConf{$V_1(\xp,\tau) \leq \left(\frac{\tilde{c}_{1,0}}{\alpha}\right)^2 + \left|\xp_2\right|^2 + \frac{\zeta^2}{M} \left(\frac{\tilde{c}_{1,0}^2}{\alpha}\right) \leq c_{1,0}$.}{
\begin{equation} \label{eqn:V1BoundNSCNest}
	V_1(\xp,\tau) \leq \left(\frac{\tilde{c}_{1,0}}{\alpha}\right)^2 + \left|\xp_2\right|^2 + \frac{\zeta^2}{M} \left(\frac{\tilde{c}_{1,0}^2}{\alpha}\right) \leq c_{1,0}.
\end{equation}}
The constants $\tilde{c}_0$, $\tilde{c}_{1,0}$, $d_0$, and $d_{1,0}$ in \eqref{eqn:UTilde0SetEquations} and \eqref{eqn:TTilde10SetEquations} comprise the hysteresis necessary to avoid chattering at the switching boundary. The idea behind these hysteresis boundaries is as follows. When $\xp \in \mathcal{U}_0$ and $\xlogic = 1$, we have that $\xp \in \overline{\reals^{2n}\setminus \T_{1,0}}$, and it is not yet time to switch to $\hczero$ but to continue to flow using $\hcone$. But once $\xp \in \T_{1,0}$ then $\xp$ is close enough to $\{\minSet\} \times \{0\}$, and the supervisor switches to $\hczero$. Note that $\T_{0,1} \cap \T_{1,0} = \emptyset$. \IfConf{Fig. \ref{fig:Hysteresis}}{Figure \ref{fig:HysteresisTR}} illustrates the hysteresis mechanism in the design of $\mathcal{U}_0$ and $\T_{1,0}$.

\IfConf{
\begin{figure}[thpb]
	\centering
	\setlength{\unitlength}{1.0pc}
	
	\vspace{-0.3cm}
	
	\begin{picture}(20,9.7)(0,0) 
		\footnotesize
%
%
		\put(0,0){\includegraphics[scale=0.2,trim={0cm 0cm 0cm 0cm},clip,width=20\unitlength]{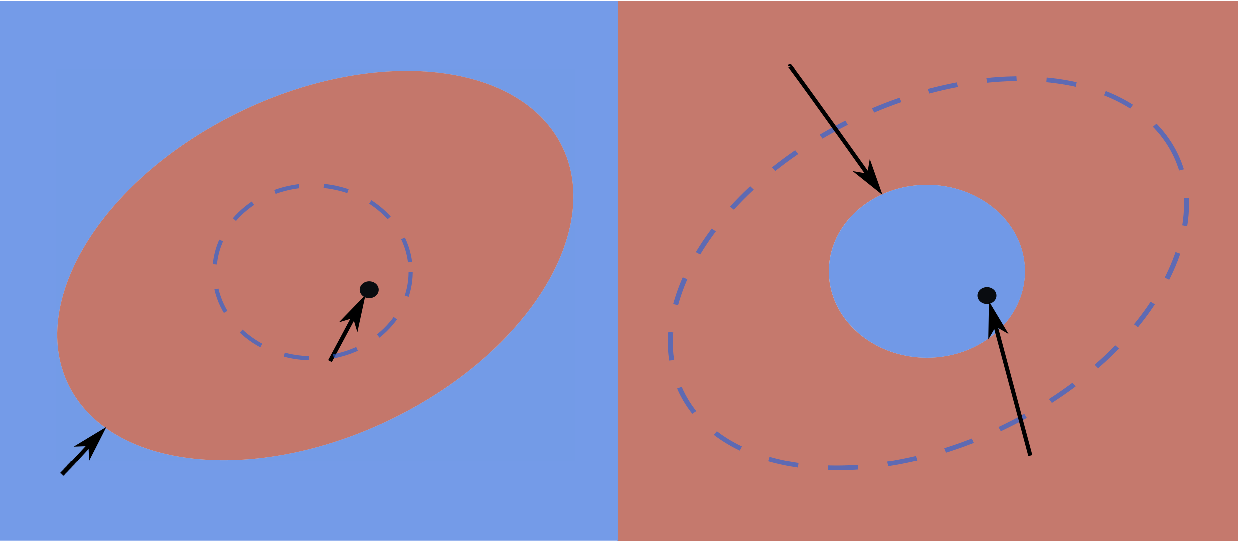}}		
		\put(4,9){$\xlogic = 0$}
		\put(14,9){$\xlogic = 1$}
		\put(0.1,0.4){$\gamma \left(\frac{\alpha}{M^2} \right) \left|\nabla L(\xp_1)\right|^2 + \frac{1}{2}\left|\xp_2\right|^2 = c_0$}
		\put(10,7.9){$\left|\nabla L(\xp_1)\right| = \tilde{c}_{1,0}, \frac{1}{2}\left|\xp_2\right|^2 = d_{1,0}$}
		\put(5,6.5){$\mathcal{U}_0$}
		\put(1,8){$\T_{0,1}$}
		\put(14.5,4.5){$\T_{1,0}$}
		\put(12,0.3){$\overline{\reals^{2n} \setminus \T_{1,0}}$}
		\put(3.3,2.3){$\{\minSet\} \times \{0\}$}
		\put(16.3,0.6){$\{\minSet\} \times \{0\}$}
	\end{picture}

	\vspace{-0.2cm}

	\caption{An illustration of hysteresis in the design of the sets $\mathcal{U}_0$, $\T_{1,0}$, and $\T_{0,1}$ on $\reals^{2n}$, via the constants $\tilde{c}_{1,0} \in (0, \tilde{c}_0)$, $d_{1,0} \in (0,d_0)$, and $c_0 > 0$. Left: due to the design of $\mathcal{U}_0$ in \eqref{eqn:U0}, every $\xp \in \mathcal{U}_0$ belongs to the $c_0$-sublevel set of the Lyapunov function $V_0$, where $V_0$ is defined via \eqref{eqn:LyapunovHBF}. Hence, the same value of $c_0 > 0$ is also used to define $\T_{0,1}$ as the closed complement of a sublevel set of $V_0$ with level equal to $c_0$. Right: the constants $\tilde{c}_{1,0} \in (0, \tilde{c}_0)$ and $d_{1,0} \in (0,d_0)$, defined via  \eqref{eqn:TTilde10SetEquations}, are chosen such that the set $\T_{1,0}$ in \eqref{eqn:T10} is contained in the interior of $\mathcal{U}_0$.}
	\label{fig:Hysteresis}
\end{figure}}{
\begin{figure}[thpb]
	\centering
	\setlength{\unitlength}{1.0pc}
	
	\begin{picture}(24,11.5)(0,0) 
		\footnotesize
		%
		%
		\put(0,0){\includegraphics[scale=0.2,trim={0cm 0cm 0cm 0cm},clip,width=24\unitlength]{Figures/Hysteresis.eps}}		
		\put(5,10.8){$\xlogic = 0$}
		\put(17,10.8){$\xlogic = 1$}
		\put(0.3,0.5){$\gamma \left(\frac{\alpha}{M^2} \right) \left|\nabla L(\xp_1)\right|^2 + \frac{1}{2}\left|\xp_2\right|^2 = c_0$}
		\put(13,9.5){$\left|\nabla L(\xp_1)\right| = \tilde{c}_{1,0}, \frac{1}{2}\left|\xp_2\right|^2 = d_{1,0}$}
		\put(5,7){$\U_0$}
		\put(1,9){$\T_{0,1}$}
		\put(17,5){$\T_{1,0}$}
		\put(13,0.5){$\overline{\reals^{2n} \setminus \T_{1,0}}$}
		\put(4.1,2.7){$\{\minSet\} \times \{0\}$}
		\put(19,1){$\{\minSet\} \times \{0\}$}
	\end{picture}
	\caption{An illustration of hysteresis in the design of the sets $\U_0$, $\T_{1,0}$, and $\T_{0,1}$ on $\reals^{2n}$, via the constants $\tilde{c}_{1,0} \in (0, \tilde{c}_0)$, $d_{1,0} \in (0,d_0)$, and $c_0 > 0$. Left: due to the design of $\U_0$ in \eqref{eqn:U0}, every $\xp \in \U_0$ belongs to the $c_0$-sublevel set of the Lyapunov function $V_0$, where $V_0$ is defined via \eqref{eqn:LyapunovHBF}. Hence, the same value of $c_0 > 0$ is also used to define $\T_{0,1}$ as the closed complement of a sublevel set of $V_0$ with level equal to $c_0$. Right: the constants $\tilde{c}_{1,0} \in (0, \tilde{c}_0)$ and $d_{1,0} \in (0,d_0)$, defined via \eqref{eqn:TTilde10SetEquations}, are chosen such that the set $\T_{1,0}$ in \eqref{eqn:T10} is contained in the interior of $\U_0$.}
	\label{fig:HysteresisTR}
\end{figure}}

\subsubsection{Design of the Set $\T_{0,1}$} \label{sec:T01}
\IfConf{T}{Recall from lines 3-4 of Algorithm \ref{algo:HS-Uniting} that t}he objective is to design $\T_{0,1}$ such that when $\xp \in \T_{0,1}$, $\xlogic = 0$, and $\tau = 0$, the state component $\xp_1$ is far from $\minSet$ and the supervisor resets $\xlogic$ to $1$ and assigns $u$ to $\hcone(\hp_1(\xp,\tau),\tau)$ so that $\hcone$ steers $\xp_1$ back to nearby $\minSet$. 
Given $c_0 > 0$, let $\alpha > 0$ come from Assumption \ref{ass:QuadraticGrowth}, and let $M > 0$ come from Assumption \ref{ass:Lipschitz}. Then, using Assumption \ref{ass:Lipschitz} with $u_1 = \xp_1^*$ and $w_1 = \xp_1$ yields $\left|\nabla L(\xp_1)\right| \leq M \left| \xp_1 - \xp_1^* \right|$
for all $\xp_1 \in \reals^n$. 
Since $L$ has quadratic growth away from $\minSet$ by Assumption \ref{ass:QuadraticGrowth}, then dividing both sides of $\left|\nabla L(\xp_1)\right| \leq M \left| \xp_1 - \xp_1^* \right|$
by $M$ and substituting into \IfConf{Assumption \ref{ass:QuadraticGrowth}}{\eqref{eqn:QuadraticGrowth}} leads to $L(\xp_1) - L^* \geq \frac{\alpha}{M^2} \left|\nabla L(\xp_1)\right|^2$,
where $\alpha > 0$ comes from Assumption \ref{ass:QuadraticGrowth}.
Then, $V_0$ in \eqref{eqn:LyapunovHBF} is lower bounded as follows: for each $\xp \in \reals^{2n}$, \IfConf{$V_0(\xp) = \gamma \left(L(\xp_1) - L^*\right) + \frac{1}{2}\left|\xp_2\right|^2 \geq \gamma \left(\frac{\alpha}{M^2} \right) \left|\nabla L(\xp_1)\right|^2 + \frac{1}{2}\left|\xp_2\right|^2$. Using such a lower bound}{
\begin{equation} \label{eqn:LowerBoundV0}
	V_0(\xp) = \gamma \left(L(\xp_1) - L^*\right) + \frac{1}{2}\left|\xp_2\right|^2 \geq \gamma \left(\frac{\alpha}{M^2} \right) \left|\nabla L(\xp_1)\right|^2 + \frac{1}{2}\left|\xp_2\right|^2.
\end{equation}
Using the right-hand side of \eqref{eqn:LowerBoundV0}} and the same $c_0 > 0$ as in Section \ref{sec:U0}, we define the set 
\begin{equation} \label{eqn:T01}
	\!\!\!\!\!\T_{0,1} := \defset{\!\xp \in \reals^{2n}\!\!\!}{\!\!\!\gamma \left(\frac{\alpha}{M^2} \right)\! \left|\nabla L(\xp_1)\right|^2 + \frac{1}{2}\left|\xp_2\right|^2 \geq c_0\!\!\!}.
\end{equation}
The set in \eqref{eqn:T01} defines the (closed) complement of a sublevel set of the Lyapunov function $V_0$ in \eqref{eqn:LyapunovHBF} with level equal to $c_0$. The constant $c_0$ is also a part of the hysteresis mechanism, as shown in \IfConf{Fig. \ref{fig:Hysteresis}}{Figure \ref{fig:HysteresisTR}}. When $\xp \in \mathcal{U}_0$, $\xlogic = 0$, and $\tau = 0$, then the supervisor does not need to switch to $\hcone$, as the state component $\xp$ is close enough to the minimizer to keep using $\hczero$. But if $\xp \in \T_{0,1}$ while $\xlogic = 0$ and $\tau = 0$, then $\xp$ is far enough from the minimizer, and the supervisor then switches to $\hcone$.

While the constants $\tilde{c}_0$, $\tilde{c}_{1,0}$, $d_0$, $d_{1,0}$, and the set $\T_{0,1}$ in \eqref{eqn:T01} depend on the constants $M > 0$ and $\alpha > 0$ which characterize the objective function $L$, as long as $M$ and $\alpha$ are positive, the uniform asymptotic stability property established in the forthcoming Theorem \ref{thm:GASNestNSC} still holds. As long as $M > 0$ and $\alpha > 0$ belong to a known set, the parameters $\tilde{c}_0$, $\tilde{c}_{1,0}$, $d_0$, and $d_{1,0}$ can still be tuned, treating such tuning as a worst-case tuning problem.

\subsection{Design of the Parameter $\lambda$}
\label{sec:DesignOfLambda}
The heavy ball parameter $\lambda > 0$ should be made large enough to avoid oscillations near the minimizer, as stated in Sections \ref{sec:Background}, \ref{sec:Contributions}, and \ref{sec:ProblemStatement}. To gain some intuition on how to tune $\lambda$, consider the quadratic objective function $L(\xp_1) = \frac{1}{2} a_1 \xp_1^2$, $a_1 > 0$, which was analyzed in detail in \cite{attouch2000heavy}. For such a case, solutions to the heavy ball algorithm are overdamped (i.e., converge slowly with no oscillations) when $\lambda > 2\sqrt{a_1}$, critically damped (i.e., the fastest convergence possible with no oscillations) when $\lambda = 2\sqrt{a_1}$, and underdamped (fast convergence with oscillations) when $\lambda < 2\sqrt{a_1}$. Therefore, setting $\lambda \geq 2\sqrt{a_1}$ gives the desired behavior of solutions to $\HS_0$, for such an objective function. 
More generally, setting $\lambda$ sufficiently large to avoid oscillations suffices, in practice. Numerically, $\lambda$ can be tuned as follows. Choose an arbitrarily large value of $\lambda$. If there is still oscillations or overshoot locally, despite the switch from $\hcone$ to $\hczero$ being made near the minimizer, then gradually increase $\lambda$ until the oscillations and overshoot disappear. See \IfConf{Example \ref{ex:NSC}}{Examples \ref{ex:Robustness}, \ref{ex:NSC}, and \ref{ex:tradeOff}} where $\lambda$ was tuned in such a way.

\subsection{Well-posedness of the hybrid closed-loop system $\HS$}

When $L$ satisfies Assumptions \ref{ass:LisNSCVX}, \ref{ass:QuadraticGrowth}, and \ref{ass:Lipschitz}, the hybrid closed-loop system $\HS$ in \eqref{eqn:HS-TimeVarying} satisfies the hybrid basic conditions \IfConf{in \cite[Assumption~6.5]{65}.}{from \cite{65} and \cite{220}, defined as follows.
	\IfConf{\begin{defn}[Hybrid basic conditions]}{\begin{definition}[Hybrid basic conditions]} \label{def:HBCs} A hybrid system $\HS$ is said to satisfy the hybrid basic conditions if its data $(C,F,D,G)$ is such that
			\begin{enumerate}[label={(A\arabic*)},leftmargin=*]
				\item $C$ and $D$ are closed subsets of $\reals^n$;
				\item $F : \reals^n \rightrightarrows \reals^n$ is outer semicontinuous and locally bounded relative to $C$, $C \subset \mathrm{dom} \ F$, and $F(x)$ is convex for every $x \in C$; 
				\item \label{item:A3} $G : \reals^n \rightrightarrows \reals^n$ is outer semicontinuous and locally bounded relative to $D$, and $D \subset \mathrm{dom} \ G$.
			\end{enumerate}
			\IfConf{\end{defn}}{\end{definition}}
		
	The hybrid basic conditions are defined in such a manner as to ensure that a hybrid system satisfying such conditions does not have any discontinuities. In particular, for a discontinuous system, if the state starts close to one of the points of a discontinuity and small measurement noise is present, the solution remains nearby such a point, even when the noise is arbitrarily small. The limit of such a solution as the noise goes to zero is a solution to a differential inclusion (or a difference inclusion, or a hybrid inclusion) which is a {\em Krasovskii regularization} of the discontinuous system. Such a solution, when the right-hand side of the Krasovskii regularization is bounded, is also a {\em Hermes solution}, and represents an equilibrium point of the discontinuous system, from which the state cannot converge to the set of interest. Conversely, when a hybrid system satisfies the hybrid basic conditions of closed sets $C$ and $D$ and outer semicontinuous maps $F$ and $G$, this ensures that any existing stability properties of such a hybrid system are robust to small perturbations. See \cite{65} and \cite{220} for more details.} 
The satisfaction of such conditions is demonstrated in the following lemma. \IfConf{Its proof is in \cite{dhustigs2022unitingNSC}.}{}
\IfConf{}{A hybrid closed-loop system $\HS$ that satisfies the hybrid basic conditions is said to be well-posed in the sense that the limit of a graphically convergent sequence of solutions to $\HS$ having a mild boundedness property is also a solution to $\HS$ \cite{65}.}

\IfConf{\begin{lem}[Well-posedness of $\HS$]}{\begin{lemma}(Well-posedness of $\HS$):} 
	\label{lemma:HBC}
	Let the function $L$ satisfy Assumptions \ref{ass:LisNSCVX}, \ref{ass:QuadraticGrowth}, and \ref{ass:Lipschitz}. Let the sets $\mathcal{U}_0$, $\T_{1,0}$, and $\T_{0,1}$ be defined via \eqref{eqn:T10}, and \eqref{eqn:T01}, respectively. Let the functions $\bar{d}$ and $\bar{\beta}$ be defined as in \eqref{eqn:dBarBetaBar}. Let $\hczero$ and $\hcone$ be defined via \eqref{eqn:StaticStateFeedbackLawsNSC}. Then, the hybrid closed-loop system $\HS$ in \eqref{eqn:HS-TimeVarying} satisfies the hybrid basic conditions.
\IfConf{\end{lem}}{\end{lemma}}

\IfConf{}{\begin{proof}
	See Section \ref{sec:ProofHBCLemma}.
\end{proof}} 

In Theorem \ref{thm:GASNestNSC} we show that $\HS$ has a compact pre-asymptotically stable set. In light of this property, Lemma \ref{lemma:HBC} is key as it leads to pre-asymptotic stability that is robust to small perturbations \cite[Theorem~7.21]{65}.
In the case of gradient-based algorithms, for instance, such perturbations can take the form of small noise in measurements of the gradient.

\IfConf{Under Assumptions \ref{ass:LisNSCVX}, \ref{ass:QuadraticGrowth}, and \ref{ass:Lipschitz}, every maximal solution to $\HS$ is complete; see \cite[Section~2.6]{dhustigs2022unitingNSC} for such a result.}{}

\IfConf{}{
\subsection{Existence of solutions to $\HS$}

Under Assumptions \ref{ass:LisNSCVX}, \ref{ass:QuadraticGrowth}, and \ref{ass:Lipschitz}, every maximal solution to $\HS$ is complete\footnote{A solution $x$ to $\HS$ is called maximal if it cannot be extended further. A solution is called complete if its domain is unbounded.} and bounded, as stated in the following lemma. Such a property is useful since it guarantees that nontrivial solutions to $\HS$ exist from each initial point in $C \cup D$, and that such solutions do not escape $C \cup D$. When every maximal solution is complete, then uniform global pre-asymptotic stability\footnote{Uniform global pre-asymptotic stability indicates the possibility of a maximal solution that is not complete, even though it may be bounded.} of the set $\mathcal{A}$ becomes UGAS. The following lemma also states that $\Pi(C_0) \cup \Pi(D_0) = \reals^{2n}$ and $\Pi(C_1) \cup \Pi(D_1) = \reals^{2n}$. Such a property ensures that nontrivial solutions to $\HS$, which exist from each initial point in $C \cup D$, also exist from any initial point in $\reals^{2n} \times \xlogicSpace \times \reals_{\geq 0}$. 

\IfConf{\begin{prop}[Existence of solutions to $\HS$]}{\begin{proposition}(Existence of solutions to $\HS$):}
	\label{prop:Existence}
	Let the function $L$ satisfy Assumptions \ref{ass:LisNSCVX}, \ref{ass:QuadraticGrowth}, and \ref{ass:Lipschitz}. Let the sets $\mathcal{U}_0$, $\T_{1,0}$, and $\T_{0,1}$ be defined via \eqref{eqn:T10}, and \eqref{eqn:T01}, respectively. Let the functions $\bar{d}$ and $\bar{\beta}$ be defined as in \eqref{eqn:dBarBetaBar}. Let $\hczero$ and $\hcone$ be defined via \eqref{eqn:StaticStateFeedbackLawsNSC}. Then, $\Pi(C_0) \cup \Pi(D_0) = \reals^{2n}$, $\Pi(C_1) \cup \Pi(D_1) = \reals^{2n}$, and each maximal solution $(t,j) \mapsto x(t,j) = (\xp(t,j),\xlogic(t,j),\tau(t,j))$ to $\HS$ in \eqref{eqn:HS-TimeVarying} is bounded and complete. 
\IfConf{\end{prop}}{\end{proposition}}

\IfConf{}{\begin{proof}
	See Section \ref{sec:ProofExistence}.
\end{proof}}
}

\subsection{Main Result}
\label{sec:MainResult}
\IfConf{The result in this section depends on the notion of UGAS, which is defined in \cite{220} and \cite{65} as follows.
	
\begin{defn}[UGAS]
	Given a hybrid closed-loop system $\HS$ as in \eqref{eqn:GeneralH}, a nonempty set $\mathcal{A} \subset \reals^n$ is said to be \underline{uniformly globally stable} for $\HS$ if there exists a class-$\mathcal{K}_{\infty}$ function $\alpha$ such that any solution $x$ to $\HS$ satisfies $\left|x(t,j)\right|_{\mathcal{A}} \leq \alpha\left(\left|x(0,0)\right|_{\mathcal{A}}\right)$ for all $(t,j) \in \dom x$; \underline{uniformly globally pre-attractive (UGpA)} for $\HS$ if for each $\eps > 0$ and $\delta > 0$ there exists $T > 0$ such that, for any solution $x$ to $\HS$ with $\left|x(0,0)\right|_{\mathcal{A}} \leq \delta$, $(t,j) \in \dom x$ and $t + j \leq T$ imply $\left|x(t,j)\right|_{\mathcal{A}} \leq \eps$; and \underline{uniformly globally pre-asymptotically stable (UGpAS)} for $\HS$ if it is both uniformly globally stable and uniformly globally pre-attractive.
\end{defn}

When every maximal solution is complete, then the prefix ``pre'' is dropped to obtain UGA and UGAS. In this section, we present a result that establishes UGAS of the set}{The result in this section depends on the notion of stability, uniform global stability, pre-attractivity, uniform global pre-attractivity, and uniform global pre-asymptotic stability (UGpAS) which are listed in the following definition, from \cite{220} and \cite{65}.
\IfConf{\begin{defn}[Stability and attractivity notions]}{\begin{definition}[Stability and attractivity notions]} \label{def:UGpAS}
		Given a hybrid \IfConf{closed-loop}{\\closed-loop} system $\HS$ as in \eqref{eqn:GeneralH}, a nonempty set $\mathcal{A} \subset \reals^n$ is said to be 
		\begin{itemize}[leftmargin=*]
			\item \underline{Stable} for $\HS$ if for each $\eps > 0$ there exists $\delta > 0$ such that each solution $x$ to $\HS$ with $\left|x(0,0)\right|_{\mathcal{A}} \leq \delta$ satisfies $\left|x(t,j)\right|_{\mathcal{A}} \leq \eps$ for all\footnote{The domain of $x$, namely, $\mathrm{dom} x \subset \reals_{\geq 0} \times \naturals$, is a hybrid time domain, which is a set such that for each $(T,J) \in \mathrm{dom} x$, $\mathrm{dom} x \cap \left( \left[0,T \right] \times \{0,1,\ldots, J\} \right) = \cup_{j=0}^{J} ([t_j,t_{j+1}],j)$ for a finite sequence of times $0 = t_0 \leq t_1 \leq t_2 \leq \ldots \leq t_{J+1}$.} $(t,j) \in \dom x$;
			\item \underline{Uniformly globally stable} for $\HS$ if there exists a class-$\mathcal{K}_{\infty}$ function $\alpha$ such that any solution $x$ to $\HS$ satisfies $\left|x(t,j)\right|_{\mathcal{A}} \leq \alpha\left(\left|x(0,0)\right|_{\mathcal{A}}\right)$ for all $(t,j) \in \dom x$;
			\item \underline{Pre-attractive} for $\HS$ if there exists $\mu > 0$ such that every solution $x$ to $\HS$ with $\left|x(0,0)\right|_{\mathcal{A}} \leq \mu$ is such that $(t,j) \mapsto \left|x(t,j)\right|_{\mathcal{A}}$ is bounded and if $x$ is complete then $\lim\limits_{(t,j) \in \dom x, \ t+j \rightarrow \infty} \left|x(t,j)\right|_{\mathcal{A}} = 0$;
			\item \underline{Uniformly globally pre-attractive (UGpA)} for $\HS$ if for each $\eps > 0$ and $\delta > 0$ there exists $T > 0$ such that, for any solution $x$ to $\HS$ with $\left|x(0,0)\right|_{\mathcal{A}} \leq \delta$, $(t,j) \in \dom x$ and $t + j \leq T$ imply $\left|x(t,j)\right|_{\mathcal{A}} \leq \eps$;
			\item \underline{Uniformly globally pre-asymptotically stable (UGpAS)} for $\HS$ if it is both uniformly globally stable and uniformly globally pre-attractive.
		\end{itemize}	
		\IfConf{\end{defn}}{\end{definition}}

In the notions involving convergence in Definition \ref{def:UGpAS}, when every maximal solution is complete, then the prefix ``pre'' is dropped to obtain attractivity, UGA, and UGAS. 
The prefix ``pre'' is in the notions involving convergence in Definition \ref{def:UGpAS} to allow for maximal solutions that are not complete. When every maximal solution is complete, such a property guarantees that nontrivial solutions exist from each initial point in $C \cup D$ to the hybrid system resulting from using our proposed uniting algorithm.

As was mentioned in Section \ref{sec:Background}, establishing UGAS for Nesterov's algorithm is a difficult problem to solve, due to its time-varying nature, as some solutions converge in a non-uniform way. We show in this section that our proposed uniting algorithm overcomes such a difficulty.

In this section, we present a result that establishes UGAS of the set}
\begin{align} \label{eqn:SetOfMinimizersHS-NSCVX}
	\IfConf{\mathcal{A} := & \defset{\xp \in \reals^{2n}}{\nabla L(\xp_1) = \xp_2 = 0} \times \{0\} \times \{0\}\nonumber\\
	= & \{\minSet\} \times \{0\} \times \{0\} \times \{0\}}{\mathcal{A} := \defset{\xp \in \reals^{2n}\!\!\!}{\!\!\!\nabla L(\xp_1) = \xp_2 = 0\!\!} \times \{0\} \times \{0\}
	= \{\minSet\} \times \{0\} \times \{0\} \times \{0\}}
\end{align}
and a hybrid convergence rate that, globally, is equal to $\frac{1}{(t+2)^2}$ while locally, is exponential, for the hybrid closed loop algorithm $\HS$ in \eqref{eqn:HS-TimeVarying} and \eqref{eqn:CAndDGradientsNestNSC}. Recall that the state $x := \left(\xp, \xlogic, \tau\right) \in \reals^{2n} \times \xlogicSpace \times \reals_{\geq 0}$. In light of this, the first component of $\mathcal{A}$, namely, $\{\minSet\}$, is the minimizer of $L$. The second component of $\mathcal{A}$, namely, $\{0\}$, reflects the fact that we need the velocity state $\xp_2$ to equal zero in $\mathcal{A}$ so that solutions are not pushed out of such a set. The third component in $\mathcal{A}$, namely, $\{0\}$, is due to the logic state ending with the value $\xlogic = 0$, namely using $\kappa_0$ as the state $\xp$ reaches the set of minimizers of $L$. The last component in $\mathcal{A}$ is due to $\tau$ being set to, and then staying at, zero when the supervisor switches to $\hczero$. 

\IfConf{\begin{thm}[UGAS of $\mathcal{A}$ for $\HS$]}{\begin{theorem}(UGAS of $\mathcal{A}$ for $\HS$):} \label{thm:GASNestNSC}
	Let the function $L$ satisfy Assumptions \ref{ass:LisNSCVX}, \ref{ass:QuadraticGrowth}, and \ref{ass:Lipschitz}. Let $\zeta > 0$, $\lambda > 0$, $\gamma > 0$, $c_{1,0} \in (0,c_0)$, and $\varepsilon_{1,0} \in (0,\varepsilon_0)$ be given. Let $\alpha > 0$ be generated by Assumption \ref{ass:QuadraticGrowth}, and let $M >0$ be generated by Assumption \ref{ass:Lipschitz}. Let $\tilde{c}_{1,0} \in (0,\tilde{c}_0)$ and $d_{1,0} \in (0,d_0)$ be defined via \eqref{eqn:UTilde0SetEquations} and \eqref{eqn:TTilde10SetEquations}. Let the sets $\mathcal{U}_0$, $\T_{1,0}$, and $\T_{0,1}$ be defined via \eqref{eqn:T10}, and \eqref{eqn:T01}, respectively. Let the functions $\bar{d}$ and $\bar{\beta}$ be defined as in \eqref{eqn:dBarBetaBar}, and let $\hczero$ and $\hcone$ be defined via \eqref{eqn:StaticStateFeedbackLawsNSC}. Then, the set $\mathcal{A}$, defined via \eqref{eqn:SetOfMinimizersHS-NSCVX}, is UGAS for $\HS$ given in \eqref{eqn:HS-TimeVarying}-\eqref{eqn:CAndDGradientsNestNSC}. Furthermore, each maximal solution $(t,j) \mapsto x(t,j) = (\xp(t,j), \xlogic(t,j), \tau(t,j))$ of the hybrid closed-loop algorithm $\HS$ starting from $C_1$ with $\tau(0,0)=0$ satisfies the following:
	\begin{enumerate}[label={\arabic*)},leftmargin=*]
		\item \label{item:1} The domain $\dom x$ of the solution $x$ is of the form\footnote{We define the interval $I^j := \defset{t\!\!}{\!\!(t,j) \in \dom x\!\!}$.} $\cup_{j=0}^{1} (I^j \times \{j\})$, with $I^0$ of the form $[t_0,t_1]$ and with $I^1$ of the form $[t_1,\infty)$ for some $t_1 \geq 0$ defining the time of the first jump. In other words, the system experiences at most one jump;
		\item \label{item:2} For each $t \in I^0$ such that\footnote{Note that at each $t \in I^0$, $q(t,0) = 1$, and at each $t \in I^1$, $q(t,1) = 0$.} $t \geq 0$
		\begin{align} \label{eqn:UnitingConvergenceRateNSCHS1}
			\IfConf{& L(\xp_1(t,0)) - L^* \nonumber\\
			& \leq \frac{4cM}{\zeta^2(t+2)^2} \left(\left|\xp_1(0,0) - \xp_1^* \right|^2 + \left| \xp_2(0,0) \right|^2 \right)}{L(\xp_1(t,0)) - L^* \leq \frac{4cM}{\zeta^2(t+2)^2} \left(\left|\xp_1(0,0) - \xp_1^* \right|^2 + \left| \xp_2(0,0) \right|^2 \right)}
		\end{align}
		where $L^* = L(\xp^*)$ and $c := \left(1 + \zeta^2\right)\exp \left(\sqrt{\frac{13}{4} + \frac{\zeta^4}{M}}\right)$. Namely, $L(\xp_1(t,0)) - L^*$ is \IfConf{$\mathcal{O}\!\left(\!\frac{4cM}{\zeta^2(t+2)^2}\!\right)$;}{\\ $\mathcal{O}\!\left(\!\frac{4cM}{\zeta^2(t+2)^2}\!\right)$;}
		\item \label{item:3} For each $t \in I^1$, $L(\xp_1(t,1)) - L^*$ \IfConf{is\\}{is} $\mathcal{O}\left(\exp \left(-(1-m)\psi t\right)\right)$,
		where $m \in (0,1)$ is such that
		$\psi := \frac{m\alpha\gamma}{\lambda} > 0$ and $\nu := \psi (\psi - \lambda) < 0$.
	\end{enumerate}
\IfConf{\end{thm}}{\end{theorem}}

As will be shown in the forthcoming proof of Theorem \ref{thm:GASNestNSC} in Section \ref{sec:ProofMainResult}, solutions starting from $C_1$ jump no more than once.
The UGAS of the hybrid closed-loop algorithm $\HS$ in Theorem \ref{thm:GASNestNSC} is proved as follows. First, in the forthcoming Proposition \ref{prop:GAS-HBF}, we establish UGAS of the set $\{\minSet\} \times \{0\}$ for the closed-loop algorithm $\HS_0$ in \eqref{eqn:H0} via Lyapunov theory and the application of an invariance principle. Then, in the forthcoming Proposition \ref{prop:UGASNest}, we prove UGAS of the set $\{\minSet\} \times \{0\} \times \reals_{\geq 0}$ for the closed-loop algorithm $\HS_1$ in \eqref{eqn:H1} via Lyapunov theory and a comparison principle. Then, UGAS of $\mathcal{A}$ for $\HS$ and item \ref{item:1} in Theorem \ref{thm:GASNestNSC} follow from a proof-by-contradiction employing the UGAS of $\HS_0$, the UGAS of $\HS_1$, and the construction of the sets $\mathcal{U}_0$, $\T_{1,0}$, and $\T_{0,1}$. The hybrid convergence rate of the closed-loop algorithm $\HS$ in items \ref{item:2} and \ref{item:3} of Theorem \ref{thm:GASNestNSC} is proved in the forthcoming Propositions \ref{prop:HBFConvergenceRate}, \ref{prop:ConvergenceNSCVXNesterov}, and \ref{prop:UpperBoundofV}.
\IfConf{\section{Numerical Example}}{\section{Numerical Examples}}
\label{sec:Examples}

\IfConf{}{In this section, we present multiple numerical examples to illustrate the hybrid closed-loop algorithm in \eqref{eqn:HS-TimeVarying} and \eqref{eqn:CAndDGradientsNestNSC}. Example \ref{ex:Robustness} first illustrates the operation of the nominal hybrid closed-loop system $\HS$, and then demonstrates the robustness of $\HS$ to different amounts of noise in measurements of $\nabla L$. Example \ref{ex:NSC} compares solutions to the hybrid closed-loop algorithm in \eqref{eqn:HS-TimeVarying} and \eqref{eqn:CAndDGradientsNestNSC} with solutions to $\HS_0$, $\HS_1$, and HAND-1 from \cite{poveda2019inducing}, with parameters chosen such that HAND-1 and $\HS$ are compared on equal footing. Example \ref{ex:NSC} then compares multiple solutions of $\HS$, starting from different initial values of $\xp_1$, to multiple solutions of HAND-1 from such initial values of $\xp_1$, to show that $\HS$ has a consistent percentage of improvement over HAND-1 for different solutions. Example \ref{ex:HHA} compares solutions to the hybrid closed-loop algorithm in \eqref{eqn:HS-TimeVarying} and \eqref{eqn:CAndDGradientsNestNSC} with solutions to $\HS_0$, $\HS_1$, HAND-1 from \cite{poveda2019inducing}, and the hybrid Hamiltonian algorithm (HHA) from \cite{teel2019first}, with parameters chosen such that HAND-1, HHA, and $\HS$ are compared on equal footing. Example \ref{ex:tradeOff} Illustrates the trade-off between speed of convergence and the resulting values of parameters for the uniting algorithm $\HS$, for different tunings of $\zeta > 0$. As in Example \ref{ex:NSC}, the parameter values for Example \ref{ex:tradeOff} are chosen such that HAND-1 and $\HS$ are compared on equal footing.}

\IfConf{}{
\IfConf{\begin{exmp}}{\begin{example}} \label{ex:Robustness}
%
	
	In this example, we simulate a solution to the nominal hybrid closed-loop system $\HS$ to illustrate how the uniting algorithm works. Then, we compare that same solution to solutions with different amounts of noise in measurements of $\nabla L$. For both the nominal system and the perturbed system, the choice of objective function, parameter values, and initial conditions are as follows. We use the objective function $L(\xp_1) := \xp_1^2$, the gradient of which is Lipschitz continuous with $M = 2$, and which has a single minimizer at $\minSet = 0$. This choice of objective function is made so that we can easily tune $\lambda$, as described in Section \ref{sec:DesignOfLambda}.
	We arbitrarily chose the heavy ball parameter value $\gamma = \frac{2}{3}$ and we tuned $\lambda$ to $200$ by choosing a value arbitrarily larger than $2\sqrt{a_1}$, where $a_1$ comes from Section \ref{sec:DesignOfLambda}, and gradually increasing it until there is no overshoot in the hybrid algorithm. For Nesterov's algorithm, we chose $\zeta = 2$. In Figure \ref{fig:NSCLargerZetaTR}, we stated that choosing $\zeta = 2$ leads to faster convergence, for Nesterov's method in \eqref{eqn:MJODE_ZetaNum_TR} and $\HS$, than choosing $\zeta = 1$. In general, convergence for such algorithms is faster as as $\zeta$ increases, and slower as $\zeta$ tends to zero.
	The parameter values for the uniting algorithm are $c_0 = 7000$, $c_{1,0} \approx 6819.68$, $\varepsilon_0 = 10$, $\varepsilon_{1,0} = 5$, and $\alpha = 1$, which yield the values $\tilde{c}_0 = 10$, $\tilde{c}_{1,0} = 5$, $d_0 = 6933$, and $d_{1,0} = 6744$, which are calculated via \eqref{eqn:UTilde0SetEquations} and \eqref{eqn:TTilde10SetEquations}. These values are chosen for proper tuning of the algorithm, in order to get nice performance, and the value of $c_{1,0}$ is chosen to exploit the properties of Nesterov's method for a longer time, so that the nominal solution gets closer to the minimizer faster. 
	Initial conditions for $\HS$ are $\xp_1(0,0) = 50$, $\xp_2(0,0) = 0$, $\xlogic(0,0) = 1$, and $\tau(0,0) = 0$. The plot \IfConf{on the top}{on the top} in \IfConf{Fig. \ref{fig:NominalNoise}}{Figure \ref{fig:NominalNoise_TR}} shows the solution to the nominal hybrid closed-loop algorithm\footnote{\label{foot:RobustnessURL}Code at \IfConf{\\}{} \texttt{gitHub.com/HybridSystemsLab/UnitingRobustness}} $\HS$, namely, the value of $\xp_1$ over time, with the time it takes for the solution to settle to within $1\%$ of $\minSet$ marked with a black dot and labeled in seconds. The jump at which the switch from $\HS_1$ to $\HS_0$ occurs is labeled with an asterisk. The solution converges quickly, without oscillations near the minimizer. \IfConf{
	
	\begin{figure}
		
		\vspace{-0.15cm}
		
		\centering
		\setlength{\unitlength}{1.0pc}
		
		\centering
		\subfloat{
			
			\begin{picture}(20,5.5)(0,0)
				\footnotesize
%
%
				\put(0,0.2){\includegraphics[scale=0.4,trim={0cm 0cm 0cm 0cm},clip,width=20\unitlength]{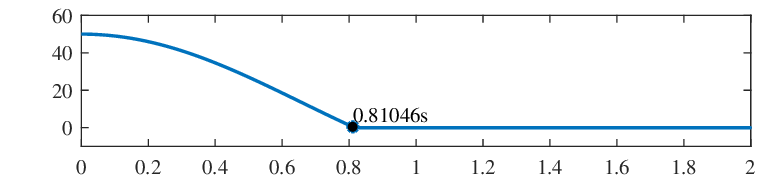}}
				\put(0,2.75){$\xp_1$}
				\put(10.6,0){$t$}
		\end{picture}}\\
		\subfloat{
			\footnotesize
			\begin{tabular}{|c|c|c|}
				\hline
				$\sigma$ & {\scriptsize $\lim\limits_{t+j \rightarrow \infty} \sup \left|\xp_1(t,j) - \xp_1^*\right|$}  & {\scriptsize $\lim\limits_{t+j \rightarrow \infty} \sup \left|L(\xp_1(t,j)) - L^*\right|$} \\
				\hline
				\hline
				$0.01$ & $8.857 \times 10^{-6}$ & $7.844 \times 10^{-11}$\\
				\hline
				$0.1$ & $8.011 \times 10^{-4}$ & $6.418 \times 10^{-7}$\\
				\hline
				$0.5$ & $9.039 \times 10^{-4}$ & $8.171 \times 10^{-7}$\\
				\hline
				$1$ & $6.982 \times 10^{-3}$ & $4.875 \times 10^{-5}$ \\
				\hline
				$5$ & $9.459 \times 10^{-3}$ & $8.947 \times 10^{-5}$ \\
				\hline
				$10$ & $1.450 \times 10^{-2}$ & $2.103 \times 10^{-4}$ \\
				\hline
				$15$ & $4.938 \times 10^{-2}$ & $2.438 \times 10^{-3}$ \\
				\hline
				$20$ & $5.992 \times 10^{-2}$ & $3.591 \times 10^{-3}$ \\
				\hline
				$25$ & $6.663 \times 10^{-2}$ & $4.439 \times 10^{-3}$ \\
				\hline
		\end{tabular}}
		\caption{Top: The evolution over time of $\xp_1$, for the nominal hybrid closed-loop algorithm $\HS$, for a function $L(\xp_1) := \xp_1^2$ with a single minimizer at $\minSet = 0$. The time at which the solution settles to within $1\%$ of $\minSet$ is marked with a dot and labeled in seconds. The jump is labeled with an asterisk. Bottom: Simulation results for perturbed solutions using zero mean Gaussian noise, with each simulation using a different value of the standard deviation $\sigma$. Results listed are for a large value of $t+j$.} 
		\label{fig:NominalNoise}
	\end{figure}
	}{\begin{figure}
		\centering
		\setlength{\unitlength}{1.0pc}
			\centering
			\subfloat{
				
				\begin{picture}(25,6)(0,0)
					\footnotesize
%
%
					\put(0,0.2){\includegraphics[scale=0.4,trim={0cm 0cm 0cm 0cm},clip,width=25\unitlength]{Figures/PlotsNominalTR.eps}}
					\put(0,3.5){$\xp_1$}
					\put(13.2,0){$t$}
			\end{picture}}\\
			\subfloat{
				\begin{tabular}{|c|c|c|}
					\hline
					$\sigma$ & $\lim\limits_{t+j \rightarrow \infty} \sup \left|\xp_1(t,j) - \xp_1^*\right|$ & $\lim\limits_{t+j \rightarrow \infty} \sup \left|L(\xp_1(t,j)) - L^*\right|$ \\
					\hline
					\hline
					$0.01$ & $8.857 \times 10^{-6}$ & $7.844 \times 10^{-11}$\\
					\hline
					$0.1$ & $8.011 \times 10^{-4}$ & $6.418 \times 10^{-7}$\\
					\hline
					$0.5$ & $9.039 \times 10^{-4}$ & $8.171 \times 10^{-7}$\\
					\hline
					$1$ & $6.982 \times 10^{-3}$ & $4.875 \times 10^{-5}$ \\
					\hline
					$5$ & $9.459 \times 10^{-3}$ & $8.947 \times 10^{-5}$ \\
					\hline
					$10$ & $1.450 \times 10^{-2}$ & $2.103 \times 10^{-4}$ \\
					\hline
					$15$ & $4.938 \times 10^{-2}$ & $2.438 \times 10^{-3}$ \\
					\hline
					$20$ & $5.992 \times 10^{-2}$ & $3.591 \times 10^{-3}$ \\
					\hline
					$25$ & $6.663 \times 10^{-2}$ & $4.439 \times 10^{-3}$ \\
					\hline
			\end{tabular}}
		\caption{Top: The evolution over time of $\xp_1$, for the nominal hybrid closed-loop algorithm $\HS$, for a function $L(\xp_1) := \xp_1^2$ with a single minimizer at $\minSet = 0$. The time at which the solution settles to within $1\%$ of $\minSet$ is marked with a dot and labeled in seconds. The jump is labeled with an asterisk. Bottom: Simulation results for perturbed solutions using zero mean Gaussian noise, with each simulation using a different value of the standard deviation $\sigma$. Results listed are for a large value of $t+j$.} 
		\label{fig:NominalNoise_TR}
	\end{figure}}
	
\IfConf{
\begin{figure}[thpb]
	\centering
	\setlength{\unitlength}{1.0pc}
	%
	%
	\begin{picture}(20,13.5)(0,0)
		\footnotesize
%
%
		\put(0.4,0.5){\includegraphics[scale=0.3,trim={0.5cm 0.4cm 0.8cm 0.4cm},clip,width=9\unitlength]{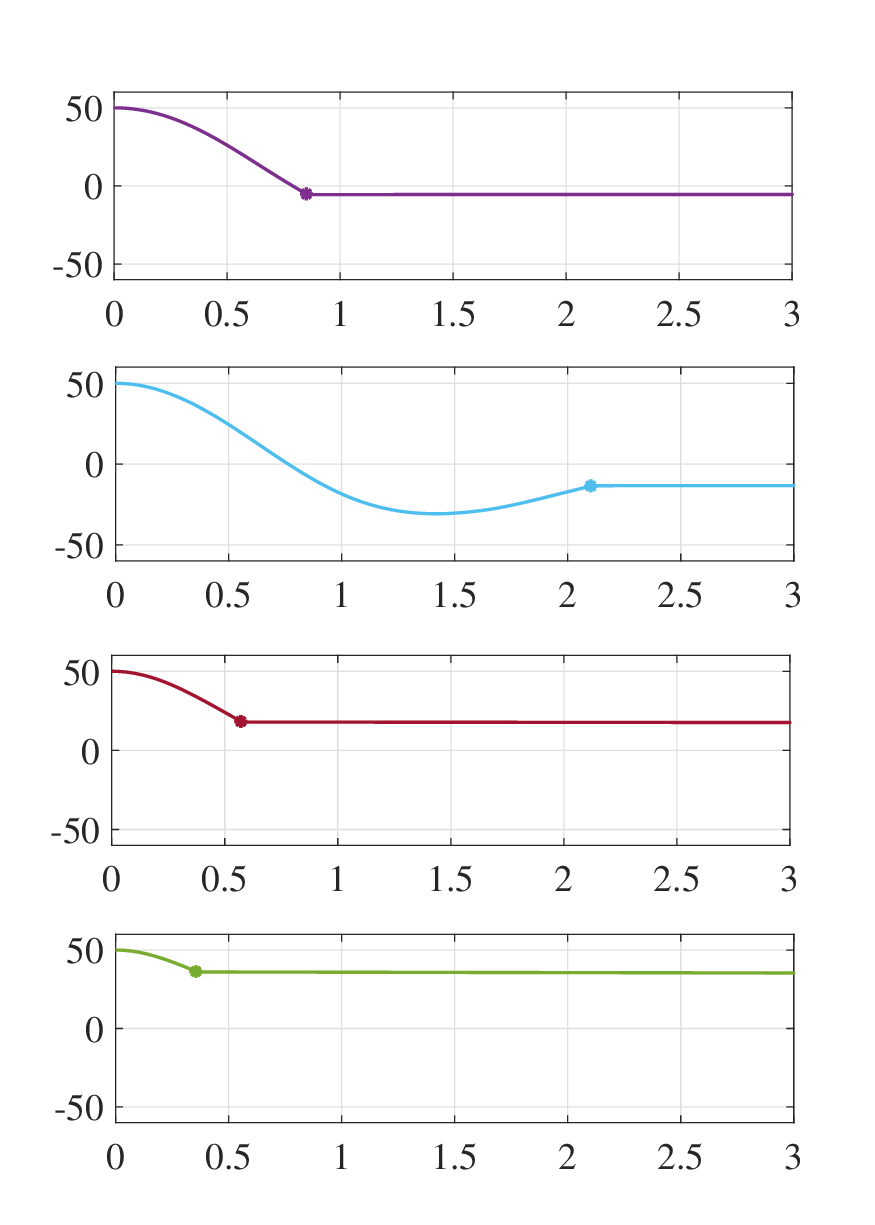}}
		\put(10,0.5){\includegraphics[scale=0.3,trim={0.5cm 0.4cm 0.7cm 0.4cm},clip,width=9\unitlength]{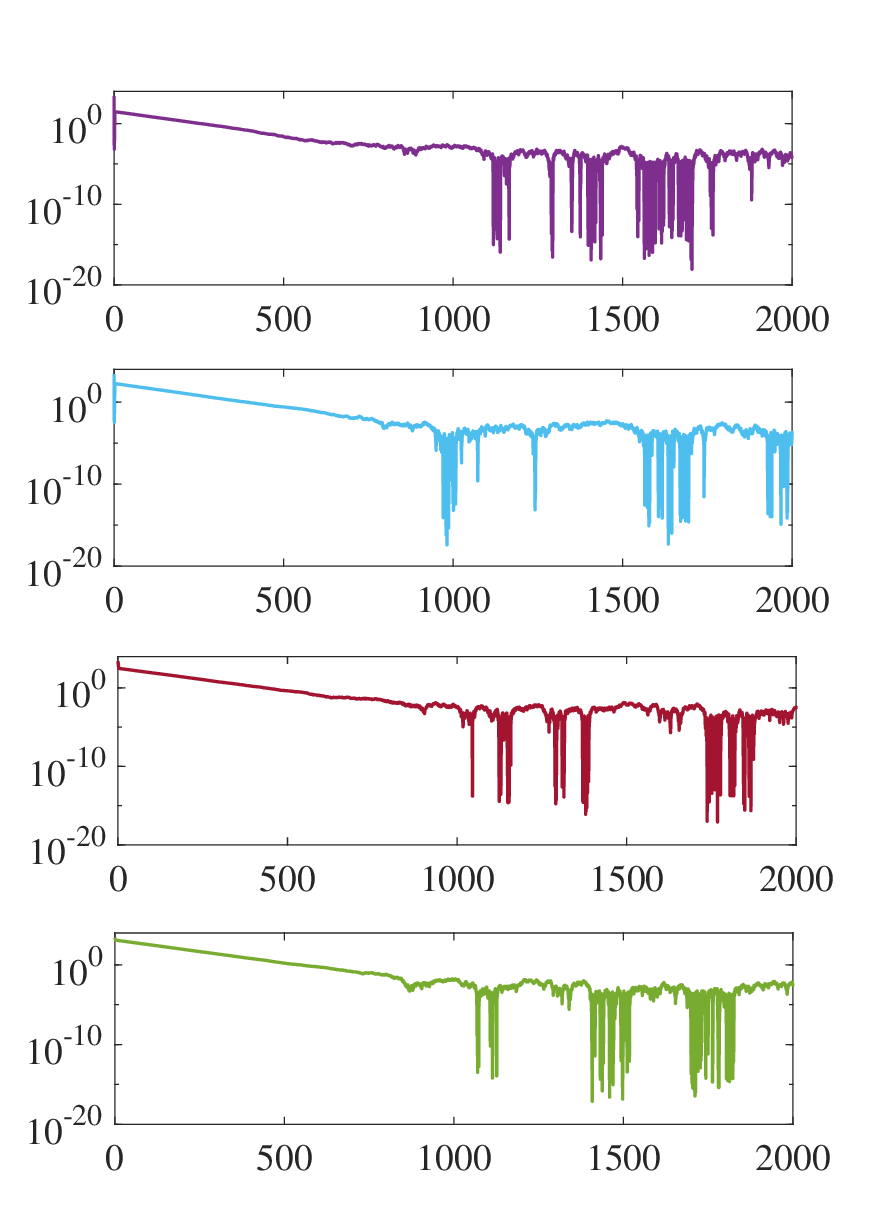}}
		\put(0.0,2.3){$\xp_1$}
		\put(0.0,5.45){$\xp_1$}
		\put(0.0,8.6){$\xp_1$}
		\put(0.0,11.8){$\xp_1$}
		\put(4.7,0.2){$t [s]$}
		\put(14.3,0.15){$t [s]$}
		\put(9.2,5.2){\rotatebox{90}{$L(\xp_1) - L^*$}}
		\put(4,3.6){{\scriptsize $\sigma = 25$}}
		\put(14,3.6){{\scriptsize $\sigma = 25$}}
		\put(4,6.7){{\scriptsize $\sigma = 20$}}
		\put(14,6.65){{\scriptsize $\sigma = 20$}}
		\put(4,10){{\scriptsize $\sigma = 10$}}
		\put(14,9.9){{\scriptsize $\sigma = 10$}}
		\put(4,13.1){{\scriptsize $\sigma = 5$}}
		\put(14,13.05){{\scriptsize $\sigma = 5$}}
	\end{picture}

	\vspace{-0.25cm}

	\caption{Simulation results for hybrid closed-loop algorithm $\HS$, for a function $L(\xp_1) := \xp_1^2$ with a single minimizer at $\minSet = 0$, with zero-mean Gaussian noise added to measurements of the gradient. Each subplot is labeled with the standard deviation used. Left subplots: the value of $\xp_1$ over time for each perturbed solution, with the jump in each solution labeled by an asterisk. Right subplots: the corresponding value of $L$ over time for each perturbed solution.}
	\label{fig:NoisyLog} 
\end{figure}}{
\begin{figure}[thpb]
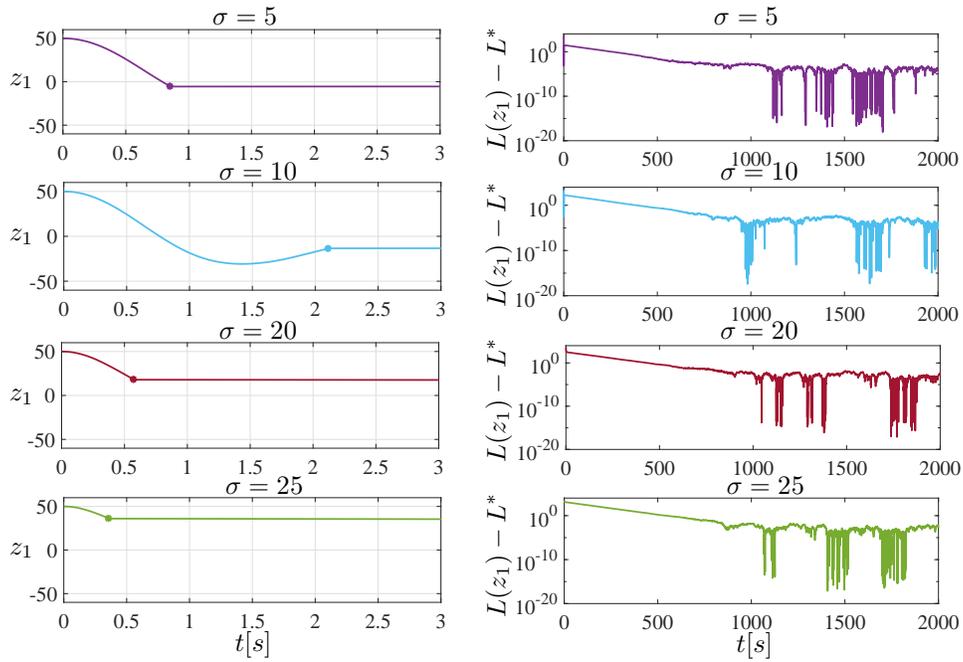

	\centering
	\setlength{\unitlength}{1.0pc}
	%
	%
	\begin{picture}(30.2,20.7)(0,0)
		\footnotesize
		%
		%
		\put(0.2,0.5){\includegraphics[scale=0.3,trim={0.5cm 0.4cm 0.8cm 0.4cm},clip,width=14\unitlength]{Figures/PlotsNoisy.eps}}
		\put(16,0.5){\includegraphics[scale=0.3,trim={0.5cm 0.4cm 0.7cm 0.4cm},clip,width=14\unitlength]{Figures/PlotsNoisyLog.eps}}
		\put(0.0,3.4){$\xp_1$}
		\put(0.0,8.25){$\xp_1$}
		\put(0.0,13.3){$\xp_1$}
		\put(0.0,18.2){$\xp_1$}
		\put(7.1,0.35){$t [s]$}
		\put(22.9,0.3){$t [s]$}
		\put(15,1.4){\rotatebox{90}{$L(\xp_1) - L^*$}}
		\put(15,6.3){\rotatebox{90}{$L(\xp_1) - L^*$}}
		\put(15,11.2){\rotatebox{90}{$L(\xp_1) - L^*$}}
		\put(15,16.1){\rotatebox{90}{$L(\xp_1) - L^*$}}
		\put(6.8,5.3){$\sigma = 25$}
		\put(22.6,5.3){$\sigma = 25$}
		\put(6.6,10.2){$\sigma = 20$}
		\put(22.4,10.2){$\sigma = 20$}
		\put(6.6,15.2){$\sigma = 10$}
		\put(22.4,15.2){$\sigma = 10$}
		\put(6.4,20.1){$\sigma = 5$}
		\put(22.2,20.1){$\sigma = 5$}
	\end{picture}
	\caption{Simulation results for hybrid closed-loop algorithm $\HS$, for a function $L(\xp_1) := \xp_1^2$ with a single minimizer at $\minSet = 0$, with zero-mean Gaussian noise added to measurements of the gradient. Each subplot is labeled with the standard deviation used. Left subplots: the value of $\xp_1$ over time for each perturbed solution, with the jump in each solution labeled by an asterisk. Right subplots: the corresponding value of $L$ over time for each perturbed solution.}
	\label{fig:NoisyLogTR}
\end{figure}}
	To show that the UGAS of $\mathcal{A}$, established in Theorem \ref{thm:GASNestNSC}, is robust to small perturbations, due to the hybrid closed-loop system $\HS$ satisfying the hybrid basic conditions by Lemma \ref{lemma:HBC}. we simulate the hybrid algorithm, using the objective function, parameter values, and initial conditions listed in the first paragraph of this example, with zero-mean Gaussian noise added to measurements of the gradient. 
	Separate simulations were run for each of the following standard deviations: $\sigma \in \{0.01,0.1,0.5,1,5,10,15,20,25\}$. \IfConf{Fig. \ref{fig:NoisyLog}}{Figure \ref{fig:NoisyLogTR}} shows some of these perturbed solutions, with each subplot labeled with the corresponding standard deviation used\footnote{Code found at same link as in Footnote \ref{foot:RobustnessURL}.}. The subplots on the left side of \IfConf{Fig. \ref{fig:NoisyLog}}{Figure \ref{fig:NoisyLogTR}} show the value of $\xp_1$ over time for different standard deviations, and the subplots on the right side of \IfConf{Fig. \ref{fig:NoisyLog}}{Figure \ref{fig:NoisyLogTR}} show the corresponding value of $L$ over time for such standard deviations. Note that, while all perturbed solutions shown in \IfConf{Fig. \ref{fig:NoisyLog}}{Figure \ref{fig:NoisyLogTR}} get close to the minimizer quickly, such perturbed solutions do not get as close to the minimizer as the solution to the nominal algorithm does; see the plot on the \IfConf{top}{top} in \IfConf{Fig. \ref{fig:NominalNoise}}{Figure \ref{fig:NominalNoise_TR}}. Also note that as the standard deviation gets larger, the corresponding perturbed solution stays slightly farther away from the minimizer. The results for all standard deviations are listed in the table \IfConf{on the bottom in}{in} \IfConf{Fig. \ref{fig:NominalNoise}}{Figure \ref{fig:NominalNoise_TR}}, showing the neighborhood of $\minSet$ that each solution settles to, for a large value of $t+j$, along with the corresponding value of $L$.
\IfConf{\end{exmp}}{\end{example}}} 

\IfConf{\begin{exmp}}{\begin{example}} \label{ex:NSC}
	In this example, to show the effectiveness of the uniting algorithm, we compare the hybrid closed-loop algorithm $\HS$, defined via \eqref{eqn:HS-TimeVarying} and \eqref{eqn:CAndDGradientsNestNSC}, with the individual closed-loop optimization algorithms $\HS_0$ and $\HS_1$ and with the HAND-1 algorithm from \cite{poveda2019inducing} which, in \cite{poveda2019inducing}, is designed and analyzed for convex functions $L$ satisfying Assumptions \ref{ass:LisNSCVX} and \ref{ass:Lipschitz}. \IfConf{The bound for HAND-1 is $L(\xp_1(t,0)) - L^* \leq \frac{B}{t^2}$ for all $\left(t,j\right) \in \dom (\xp, \tau)$ such that $j = 0$, $\xp_1(0,0) = \xp_2(0,0)$, $\tau(0,0) = T_{\min}$, $\xp_1(0,0) \in K_0 := \{\xp_1^*\} + r\ball$, where $B  := \frac{r^2}{2c_1} + T^2_{\min} \left(L(\xp_1(0,0)) - L^* \right) > 0$, $r \in \reals_{> 0}$, $c_1 > 0$. Such a rate is only guaranteed until the first jump.}{First, we compare the convergence rates of $\HS$ and HAND-1 analytically. Using an alternate state space representation, namely, $\xp_1 := \xi$ and $\xp_2 := \xi + \frac{\tau}{2}\dot{\xi}$, the HAND-1 algorithm has state $(\xp, \tau) \in \reals^{2n+1}$ and data $(C,F,D,G)$
	\begin{equation} \label{eqn:HANDFAndG}
		F(\xp, \tau) := \matt{\frac{2}{\tau}(\xp_2 - \xp_1)\\-2c_1\tau \nabla L(\xp_1) \\1} \  (\xp,\tau) \in C, \quad G(\xp, \tau) := \matt{\xp\\T_{\min}} \ (\xp,\tau) \in D
	\end{equation}
	where $c_1 > 0$ and the flow and jump sets are\\ $C := \defset{(\xp, \tau) \in \reals^{2n+1}\!\!\!}{\!\!\!\tau \in [T_{\min}, T_{\max}]\!\!}$ and\\ $D := \defset{(\xp, \tau) \in \reals^{2n+1}\!\!\!}{\!\!\!\tau \in [T_{\text{med}},T_{\max}]\!\!}$,
	with $0 < T_{\min} < T_{\text{med}} < T_{\max} < \infty$, and $T_{\text{med}} \geq \sqrt{\frac{B}{\delta_{\text{med}}}} + T_{\min} > 0$, $\delta_{\text{med}} > 0$. It is shown in \cite{poveda2019inducing} that each maximal solution $(t,j) \mapsto (\xp(t,j),\tau(t,j))$ to the HAND-1 algorithm satisfies
	\begin{equation}\label{eqn:PNRateHAND-1}
		L(\xp_1(t,0)) - L^* \leq \frac{B}{t^2}
	\end{equation}
	for all $\left(t,j\right) \in \dom (\xp, \tau)$ such that $j = 0$, $\xp_1(0,0) = \xp_2(0,0)$, $\tau(0,0) = T_{\min}$, $\xp_1(0,0) \in K_0 := \{\xp_1^*\} + r\ball$, where $B  := \frac{r^2}{2c_1} + T^2_{\min} \left(L(\xp_1(0,0)) - L^* \right) > 0$, $r \in \reals_{> 0}$, $c_1 > 0$. 
	
	For the hybrid closed-loop algorithm $\HS$, the coefficient of the bound on $\HS_1$ from \eqref{eqn:UnitingConvergenceRateNSCHS1}, namely,
		\begin{equation} \label{eqn:BoundH1Again}
			L(\xp_1(t,0)) - L^* \leq \frac{4cM}{\zeta^2(t+2)^2} \left(\left|\xp_1(0,0) - \xp_1^* \right|^2 + \left| \xp_2(0,0) \right|^2 \right) 
		\end{equation}
	for each $t \in I^0$, $t \geq 0$, at which $q(t,0) = 1$, and for each $\zeta > 0$, and $M > 0$, is $\frac{4cM}{\zeta^2}\left(\left|\xp_1(0,0) - \xp_1^* \right|^2 + \left| \xp_2(0,0) \right|^2 \right)$, where $c := \left(1 + \zeta^2\right)\exp \left(\sqrt{\frac{13}{4} + \frac{\zeta^4}{M}}\right)$. The coefficient of the bound in HAND-1 is $B := \frac{r^2}{2c_1} + T^2_{\min} \left(L(\xp_1(0,0)) - L^* \right)$. Since, as $t \rightarrow \infty$, $\frac{1}{(t+2)^2} \rightarrow \frac{1}{t^2}$, then, comparing the coefficients of the bounds, the bound in \eqref{eqn:BoundH1Again} is slightly better than \eqref{eqn:PNRateHAND-1} since $\frac{r^2}{2c_1}$ is very large for small $t$. Neglecting the $\frac{r^2}{2c_1}$ term, however, the bound on $\HS_1$ \eqref{eqn:BoundH1Again} matches \eqref{eqn:PNRateHAND-1}. The rate for HAND-1, nevertheless, is only guaranteed until the first jump. After this, there is no characterized bound for HAND-1. In contrast, $\HS$ has a characterized bound for the domain of every solution such that $t \geq 0$. Namely, it has rate $\frac{1}{(t+2)^2}$ until the state $\xp$ is within a small neighborhood of the minimizer -- where the rate then switches to $\exp \left(-(1-m)\psi t\right)$, where, given $\gamma > 0$ and $\lambda > 0$, $m \in (0,1)$ is such that $\psi = \frac{m\alpha\gamma}{\lambda} > 0$ and $\nu = \psi (\psi - \lambda) < 0$.} 
	
	Next, we compare $\HS_0$, $\HS_1$, $\HS$, and HAND-1 in simulation. To compare these algorithms, \IfConf{the choice of objective function, parameter values, and initial conditions are as follows. We use the objective function $L(\xp_1) := \xp_1^2$, the gradient of which is Lipschitz continuous with $M = 2$, and which has a single minimizer at $\minSet = 0$. This choice of objective function is made so that we can easily tune $\lambda$, as described in Section \ref{sec:DesignOfLambda}.
	We arbitrarily\footnote{Although the choice of $\gamma$ is arbitrary, we have found in general that choosing $\gamma \in (0,1)$ works well, in practice.} chose the heavy ball parameter value $\gamma = \frac{2}{3}$ and we tuned $\lambda$ to $200$ by choosing a value arbitrarily larger than $2\sqrt{a_1}$, where $a_1$ comes from Section \ref{sec:DesignOfLambda}, and gradually increasing it until there is no overshoot in the hybrid algorithm.}{we use the same objective function $L$, heavy ball parameter values $\lambda$ and $\gamma$, Lipschitz parameter $M$, Nesterov parameter $\zeta$, 
	and uniting algorithm parameter values $c_0$, $c_{1,0}$, $\varepsilon_0$, $\varepsilon_{1,0}$, $\alpha$, $\tilde{c}_0$, $\tilde{c}_{1,0}$, $d_0$, and $d_{1,0}$ as in Example \ref{ex:Robustness}.} \IfConf{In Figure \ref{fig:NSCLargerZeta}, we stated that choosing $\zeta = 2$ leads to faster convergence, for Nesterov's method in \eqref{eqn:MJODE_ZetaNum} and $\HS$, than choosing $\zeta = 1$. In general, convergence for such algorithms is faster as $\zeta$ increases, and slower as $\zeta$ tends to zero.}{}
	Given $\zeta = 2$ for Nesterov's algorithm and $\HS$, the HAND-1 parameters $c_1 = 0.5$ and $T_{\min} = \frac{1+\sqrt{7}}{2}$ are chosen such that the resulting gain coefficients for $\xp_1$ and $\xp_2$ are the same for both $\HS$ and HAND-1, so that these algorithms are compared on equal footing\footnote{Although there exist parameter values for which HAND-1 has faster, oscillation-free performance, due to the way $\HS$ and HAND-1 relate to each other, they are compared fairly for a particular set of parameters.}. The remaining HAND-1 parameters, $r$ and $\delta_{\text{med}}$, have different values depending on the initial conditions $\xp_1(0,0) = \xp_2(0,0)$, listed in \IfConf{\cite[Table~2]{dhustigs2022unitingNSC}}{Table \ref{table:PercentImprovementNSCTR}}, which leads to different values of $T_{\text{med}}$ and $T_{\max}$, for each solution. Such values are chosen such that $T_{\text{med}} \geq \sqrt{\frac{B}{\delta_{\text{med}}}} + T_{\min} > 0$. Additionally, we choose $T_{\max} = T_{\text{med}} + 1$. The parameter values for the uniting algorithm are $\varepsilon_0 = 10$, $\varepsilon_{1,0} = 5$, and $\alpha = 1$. The remaining parameter values $c_0$ and $c_{1,0}$ are different depending on the initial condition $\xp_1(0,0)$ and are listed in \IfConf{\cite[Table~2]{dhustigs2022unitingNSC}}{Table \ref{table:PercentImprovementNSCTR}}, which leads to different values of $d_0$, calculated via \eqref{eqn:UTilde0SetEquations}, and $d_{1,0}$ calculated via \eqref{eqn:TTilde10SetEquations}. These values are chosen for proper tuning of the algorithm, in order to get nice performance, and for exploiting the properties of Nesterov's method as long as we want. Initial conditions for all solutions to $\HS$ are $\xp_2(0,0) = 0$, $\xlogic(0,0) = 1$, and $\tau(0,0) = 0$, with values of $\xp_1(0,0)$ listed in \IfConf{\cite[Table~2]{dhustigs2022unitingNSC}}{Table \ref{table:PercentImprovementNSCTR}}. Initial conditions for all solutions to HAND-1 are $\tau(0,0) = T_{\min}$, with values of $\xp_1(0,0) = \xp_2(0,0)$ listed in \IfConf{\cite[Table~2]{dhustigs2022unitingNSC}}{Table \ref{table:PercentImprovementNSCTR}}.
	\IfConf{}{
	\begin{table}[thpb]
%
%
		\begin{center}
			\begin{tabular}{|c|c|c|}
				\hline
							& { Average time} & { Average $\%$}\\
				Algorithm & { to converge (s)} & { improvement} \\
				\hline
				\hline
				$\HS$ & $0.811$ & -- \\
				\hline
				$\HS_0$ & $690.759$ & $99.9$ \\
				\hline
				$\HS_1$ & $4.409$ & $81.6$ \\
				\hline
				HAND-1 & $8.649$  & $90.6$ \\
				\hline
			\end{tabular}
			\IfConf{
				\caption{Average times for which $\HS$, $\HS_0$, $\HS_1$, and HAND-1 settle to within $1\%$ of $\minSet$, and the average percent improvement of $\HS$ over each algorithm. Percent improvement is calculated via \eqref{eqn:PercentImprovement}. The objective function used for this table is $L(\xp_1) := \xp_1^2$.}		
			}{
				\caption{Average times for which $\HS$, $\HS_0$, $\HS_1$, and HAND-1 settle to within $1\%$ of $\minSet$, and the average percent improvement of $\HS$ over each algorithm. Percent improvement is calculated via \eqref{eqn:PercentImprovement}. The objective function used for this table is $L(\xp_1) := \xp_1^2$.}	
			}
			\label{table:TimePercentImpNSC}
		\end{center}
	\end{table}}

	\IfConf{The time that it takes for each algorithm to settle within $1\%$ of $\xp^*_1$, averaged over solutions starting from ten different values\footnote{Code at\IfConf{\\}{} \texttt{gitHub.com/HybridSystemsLab/UnitingDifferentICs}} of $\xp_1(0,0)$, are as follows: $0.811$ seconds for $\HS$, $690.759$ seconds for $\HS_0$, $4.409$ seconds for $\HS_1$, and $8.649$ seconds for HAND-1.}{\IfConf{Table \ref{table:TimePercentImpNSC}}{Table \ref{table:TimePercentImpNSC}} shows the time that each algorithm takes to settle within $1\%$ of $\minSet$, averaged over solutions starting from ten different values\footnote{Code at\IfConf{\\}{} \texttt{gitHub.com/HybridSystemsLab/UnitingDifferentICs}\label{foot:DifferentICsCode}} of $\xp_1(0,0)$ (listed in the first column of \IfConf{\cite[Table~2]{dhustigs2022unitingNSC}}{Table \ref{table:PercentImprovementNSCTR}}), and the average percent improvement of $\HS$ over $\HS_0$, $\HS_1$, and HAND-1, which is calculated using the following formula}
\IfConf{Using the formula}{}
	\begin{equation} \label{eqn:PercentImprovement}
		\left(\frac{\left(\text{Time of } \HS_0, \HS_1, \text{ or HAND-1}\right) - \text{Time of } \HS}{\text{Time of } \HS_0, \HS_1, \text{ or HAND-1}}\right) \times 100 \% .
	\end{equation}
\IfConf{t}{As can be seen in 
\IfConf{Table \ref{table:TimePercentImpNSC}}{Table \ref{table:TimePercentImpNSC}}, $\HS$ converges faster than the other algorithms, and t}he average percent improvement of $\HS$ over each of the other algorithms\IfConf{}{ in Table \ref{table:TimePercentImpNSC}} is $99.9\%$ over $\HS_0$, $81.6\%$ over $\HS_1$, and $90.6\%$ over HAND-1. \IfConf{
\begin{figure}[thpb]
	\centering
	\setlength{\unitlength}{1.0pc} 
	%
	%
	\begin{picture}(20,9.5)(0,0)
		\footnotesize
%
%
		\put(0,0.5){\includegraphics[scale=0.4,trim={0.7cm 0.1cm 1cm 0.3cm},clip,width=20\unitlength]{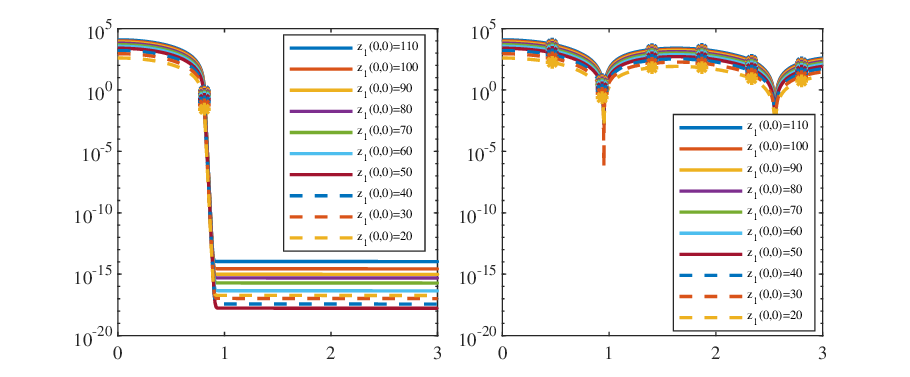}}
		\put(0,3.3){\rotatebox{90}{$L(\xp_1) - L^*$}}
		\put(5.3,0.2){$t[s]$}
		\put(15,0.2){$t[s]$}
	\end{picture}

	\vspace{-0.25cm}

	\caption{The evolution of $L$ over time, from different initial conditions, for $\HS$ (left) and HAND-1 (right). All solutions are for the objective function $L(\xp_1) := \xp_1^2$, and the parameters used for HAND-1 and $\HS$ are listed in \cite[Table~2]{dhustigs2022unitingNSC}, with different values of $c_0$ and $c_{1,0}$ for each solution of $\HS$, leading to different values of $d_0$ calculated via \eqref{eqn:UTilde0SetEquations} and $d_{1,0}$ calculated via \eqref{eqn:TTilde10SetEquations}, and different values of $r$ and $\delta_{\text{med}}$ for each solution of HAND-1, leading to different values of $T_{\text{med}}$ and $T_{\max}$. Jumps are marked with asterisks.}
	\label{fig:NSCTrajectoriesPlotsLog}
\end{figure}
}{
\begin{figure}[thpb]
	\centering
	\setlength{\unitlength}{1.0pc} 
	%
	%
	\begin{picture}(30.8,15)(0,0)
		\footnotesize
		%
		%
		\put(0,0.5){\includegraphics[scale=0.4,trim={0.8cm 0.1cm 1cm 0.3cm},clip,width=30.8\unitlength]{Figures/UnitingNSCTrajectoriesLog.eps}}
		\put(0,6){\rotatebox{90}{$L(\xp_1) - L^*$}}
		\put(8.3,0.4){$t[s]$}
		\put(23.3,0.4){$t[s]$}
	\end{picture}
	\caption{The evolution of $L$ over time, from different initial conditions, for $\HS$ (left) and HAND-1 (right). All solutions are for the objective function $L(\xp_1) := \xp_1^2$, and the parameters used for HAND-1 and $\HS$ are listed in Table \ref{table:PercentImprovementNSCTR}, with different values of $c_0$ and $c_{1,0}$ for each solution of $\HS$, leading to different values of $d_0$ calculated via \eqref{eqn:UTilde0SetEquations} and $d_{1,0}$ calculated via \eqref{eqn:TTilde10SetEquations}, and different values of $r$ and $\delta_{\text{med}}$ for each solution of HAND-1, leading to different values of $T_{\text{med}}$ and $T_{\max}$. Jumps are marked with asterisks.}
	\label{fig:NSCTrajectoriesPlotsLogTR}
\end{figure}}
	
\IfConf{}{
\begin{table}[thpb]
	\begin{center}
		\begin{tabular}{|c|c|c|c|c|c|c|c|}
			\hline
			& & & & & \multicolumn{2}{|c|}{{\footnotesize Time to converge (s)\par}} & {\footnotesize$\%$ Improve-}\\
			\cline{6-7} 
			$\xp_1(0,0)$ & $c_0$ & $c_{1,0}$ & $r$ & $\delta_{\text{med}}$ & $\HS$  & HAND-1  & {\footnotesize ment} \\
			\hline
			\hline
			$110$ & $34000$ & $32719.231$ & $111$ & $240700$ & 0.811 & 8.649 & 90.6 \\
			\hline
			$100$ & $28000$ & $27053.704$ & $101$ & $199000$ & 0.811 & 8.65 & 90.6 \\
			\hline
			$90$ & $23000$ & $21927.75$ & $91$ & $161300$ & 0.811 & 8.648 & 90.6 \\
			\hline
			$80$ & $18000$ & $17341.37$ & $81$ & $127550$ & 0.811 & 8.65 & 90.6 \\
			\hline
			$70$ & $14000$ & $13294.565$ & $71$ & $97700$ & 0.811 & 8.649 & 90.6 \\
			\hline
			$60$ & $10500$ & $9787.333$ & $61$ & $71875$ & 0.811 & 8.648 & 90.6 \\
			\hline
			$50$ & $7000$ & $6819.676$ & $51$ & $50000$ & 0.810 & 8.65 & 90.6 \\
			\hline
			$40$ & $5000$ & $4391.593$ & $41$ & $32075$ & 0.811 & 8.65 & 90.6 \\
			\hline
			$30$ & $3000$ & $2503.083$ & $31$ & $18110$ & 0.811 & 8.648 & 90.6 \\
			\hline
			$20$ & $2000$ & $1154.148$ & $21$ & $8112$ & 0.811 & 8.648 & 90.6 \\
			\hline
		\end{tabular}
		\caption{Times for which $\HS$ and HAND-1 settle to within $1\%$ of $\minSet$, and percent improvement of $\HS$ over HAND-1, for solutions from different initial conditions, shown in Figure \ref{fig:NSCTrajectoriesPlotsLogTR}. The objective function used for this table is $L(\xp_1) := \xp_1^2$.}	
		\label{table:PercentImprovementNSCTR}
	\end{center}
\end{table}
}
	
	\IfConf{Fig. \ref{fig:NSCTrajectoriesPlotsLog}}{Figure \ref{fig:NSCTrajectoriesPlotsLogTR}} compares different solutions for $\HS$ and HAND-1, from different values of $\xp_1(0,0)$,
	for the objective function
	$L(\xp_1) := \xp_1^2$. \IfConf{\cite[Table~2]{dhustigs2022unitingNSC}}{Table \ref{table:PercentImprovementNSCTR}} lists the times for which each solution settles to within $1\%$ of $\minSet$ for both $\HS$ and HAND-1, and shows the percent improvement of $\HS$ over HAND-1. 
	As can be seen in \IfConf{Fig. \ref{fig:NSCTrajectoriesPlotsLog}}{Figure \ref{fig:NSCTrajectoriesPlotsLogTR}} and in \IfConf{\cite[Table~2]{dhustigs2022unitingNSC}}{Table \ref{table:PercentImprovementNSCTR}}, the percent improvement of $\HS$ over HAND-1 for all solutions is $90.6\%$, which shows consistency in the performance of $\HS$ versus HAND-1.
	
	\IfConf{}{The bound for HAND-1, shown in \eqref{eqn:PNRateHAND-1} and which holds only until the first reset, is only guaranteed when $\xp_1(0,0) = \xp_2(0,0)$. This leads to a required nonzero velocity for HAND-1 in most scenarios, which leads to overshoot. In contrast, $\HS$ has no such constraint on $\xp_2(0,0)$, which can be set to zero in all scenarios. The lack of such a constraint on the initial condition $\xp_2(0,0)$ for the hybrid closed-loop algorithm $\HS$ is essential to its improved performance over HAND-1, as the overshoot in solutions to HAND-1 due to $\xp_1(0,0) = \xp_2(0,0)$ leads to a slower convergence time than for $\HS$, as seen in Table \ref{table:TimePercentImpNSC}. Moreover, as described previously in this example, no bound for HAND-1 is characterized after the first reset, whereas the (hybrid) convergence bound characterized for $\HS$ holds for the domain of every solution such that $t \geq 1$.}
\IfConf{\end{exmp}}{\end{example}}
\IfConf{
	
\vspace{-0.2cm}
	
For an illustration of the robustness of the UGAS property in Theorem \ref{thm:GASNestNSC} to perturbations, see \cite[Example~4.1]{dhustigs2022unitingNSC}. For a comparison between $\HS$, HAND-1, and the hybrid Hamiltonian algorithm (HHA) in \cite{teel2019first}, see \cite[Example~4.3]{dhustigs2022unitingNSC}. For an illustration of the trade-off between speed of convergence and the resulting values of parameters for the uniting algorithm $\HS$, for different values of $\zeta > 0$, see \cite[Example~4.3]{dhustigs2022unitingNSC}. 

\vspace{-0.15cm}

}{
\begin{example} \label{ex:HHA}
	In this example, we compare the hybrid closed-loop algorithm $\HS$, defined via \eqref{eqn:HS-TimeVarying} and \eqref{eqn:CAndDGradientsNestNSC}, with the individual closed-loop optimization algorithms $\HS_0$ and $\HS_1$, the HAND-1 algorithm from \cite{poveda2019inducing}, and the hybrid Hamiltonian algorithm (HHA) from \cite{teel2019first}. In \cite{teel2019first}, the HHA algorithm is designed and analyzed for objective functions $L$ which satisfy the following assumptions:
	\begin{assumption}\label{ass:HHA_Assumps} \hfill
		\begin{enumerate}[label={(L\arabic*)},leftmargin=*]
			\item The function $L$ is $\mathcal{C}^1$;
			\item \label{item:L2} The function $L$ has compact sublevel sets;
			\item \label{item:L3} $\nabla L$ is Lipschitz continuous with constant $M \in \left(0, \bar{M}\right]$.
		\end{enumerate}
	\end{assumption}
	
	First, we compare the convergence rates of $\HS$ and HHA analytically. Letting $H : \reals^{2n} \rightarrow \reals_{\geq 0}$ denote the separable Hamiltonian
	\begin{equation} \label{eqn:Hamiltonian}
		H(\xp) := L(\xp^*_1) - L^* + \frac{1}{2} \left|\xp_2\right|^2
	\end{equation}
	and letting $J$ and $R_0$ be defined as
	\begin{equation}
		J := \matt{0 & I_n\\-I_n & 0}, R_0 := \matt{I_n & 0\\0 & 0}
	\end{equation}
	the HHA algorithm has state $\left(\xp,\tau\right) \in \reals^{2n+1}$ and data $(C,F,D,G)$
	\begin{subequations} \label{eqn:HHA}
		\begin{align}
			F(\xp,\tau) & := \matt{J\nabla H(\xp) \\ 1} = \matt{ \xp_2\\-\nabla L(\xp)\\ 1} \quad (\xp,\tau) \in C\\
			G(\xp,\tau) & := \matt{ R_0\xp\\ 0} = \matt{\xp_1\\ \xp_2\\0} \quad (\xp,\tau) \in D
		\end{align}
	\end{subequations}
	where the flow and jump sets are
	\begin{subequations} \label{eqn:CAndDHHA}
		\begin{align}
			C & := C_0 \times \left[0,\bar{T}\right]\\
			D & := \left(C_0 \times \left\{\bar{T}\right\} \right) \cup \left(D_0 \times \left[0, \bar{T}\right] \right)
		\end{align}
	\end{subequations}
	where $C_0 := \defset{\xp \in \reals^{2n}\!\!}{\!\!\left\langle \nabla L(\xp_1),\xp_2 \right\rangle \leq 0\!\!}$,\\ $D_0 := \defset{\xp \in \reals^{2n}\!\!}{\!\!\left\langle \nabla L(\xp_1),\xp_2 \right\rangle = 0, \left|\xp_2\right|^2 \geq \left|\nabla L(\xp_1)\right|^2/\bar{M}\!\!}$, and $\bar{T} \in \left(0,\infty\right)$ is the timeout parameter of the timer $\tau$. The following assumption is imposed on HHA in \eqref{eqn:HHA}-\eqref{eqn:CAndDHHA}:
	\begin{assumption}\label{ass:Timeout}
		No solution to the flow dynamics starting from $\left(\xp_1(0,0),0,0\right)$ with $\nabla L(\xp_1) \neq 0$ causes the timer to timeout, i.e., reaches $\left(C_0 \setminus D_0\right) \times \left\{\bar{T}\right\}$.
	\end{assumption}
	It is shown in \cite{teel2019first} that Assumption \ref{ass:Timeout} holds for quadratic functions $L$ with $A = A^{\top} > 0$, when $\bar{T} := \frac{n\pi}{2\sqrt{\lambda_{\min}(A)}}$, where $\lambda_{\min}(A)$ denotes the minimum eigenvalue of $A$.

The forthcoming convergence rate results in \cite{teel2019first} hold when $L$ satisfies the following assumption:
\begin{assumption}\label{ass:PL}
	The function $L$ satisfies the Polyak-\L ojasiewicz inequality with $\varsigma > 0$, namely,
	\begin{equation}
		\left|\nabla L(\xp_1) \right|^2 \geq 2\varsigma \left(L(\xp_1) - L^*\right).
	\end{equation}
\end{assumption}
It is shown in \cite{teel2019first} that when $L$ satisfies Assumptions \ref{ass:HHA_Assumps} and \ref{ass:PL}, and HHA in \eqref{eqn:HHA}-\eqref{eqn:CAndDHHA} satisfies Assumption \ref{ass:Timeout}, each maximal solution $\left(t,j\right) \mapsto (\xp(t,j),\tau(t,j))$ to HHA in \eqref{eqn:HHA}-\eqref{eqn:CAndDHHA} satisfies
\begin{equation}\label{eqn:ConvergenceBoundHHA}
	L(\xp_1(t,j)) - L^* \leq \left(L(\xp_1(0,0)) - L^*\right) \min \left\{1,\exp\left(-\Omega \left(t - \theta\right)\right)\right\} 
\end{equation}
for each $\left(t,j\right) \in \dom (\xp,\tau)$ such that $\xp_2(0,0) = 0$ and $\tau(0,0) = 0$, where 
\begin{equation}\label{eqn:OmegaTheta}
	\Omega := \ln \left(1+\frac{\varsigma}{M}\right)\bar{T}^{-1}, \quad \theta := \bar{T}
\end{equation} 
and where $\varsigma > 0$ comes from Assumption \ref{ass:PL} and $M > 0$ comes from Assumption \ref{ass:HHA_Assumps}. Comparing the bound in \eqref{eqn:ConvergenceBoundHHA} with the bound on $\HS_1$ in \eqref{eqn:H1}, from \eqref{eqn:UnitingConvergenceRateNSCHS1}, the HHA algorithm reaches the neighborhood of the minimizer faster than the hybrid closed-loop algorithm $\HS$ in \eqref{eqn:HS-TimeVarying}-\eqref{eqn:CAndDGradientsNestNSC}. The assumptions imposed on $L$ for HHA, however, are different from the assumptions on $L$ for the hybrid closed-loop algorithm $\HS$. While some assumptions on $L$ are weaker for HHA -- for instance, Assumption \ref{ass:PL} -- other assumptions are stronger, such as items \ref{item:L2} and \ref{item:L3} of Assumption \ref{ass:HHA_Assumps}. Moreover, Assumption \ref{ass:Timeout} is imposed on HHA, which is not imposed on the hybrid closed-loop algorithm $\HS$. Table \ref{table:AssumpsProperties} summarizes the assumptions and the convergence rate results for HHA, $\HS$, and HAND-1.  

\begin{table}[thpb]
	\begin{center}
		\begin{tabular}{|c|c|c|}
			\hline
			Algorithm & Assumptions & Convergence Rate \\
			\hline
			\hline
			\multirow{3}{*}{$\HS$} & Assumptions & $\frac{1}{\left(t+2\right)^2}$ globally and \\
			&\ref{ass:LisNSCVX}, \ref{ass:QuadraticGrowth}, & $\exp \left(-(1-m)\psi t\right)$ locally, where\\
			& and \ref{ass:Lipschitz} & $m \in (0,1)$ s.t. $\psi := \frac{m\alpha\gamma}{\lambda} > 0$\\
			& & and $\nu := \psi (\psi - \lambda) < 0$. \\
			\hline
			\multirow{2}{*}{HAND-1} & Assumptions & $\frac{1}{t^2}$ until the first\\
			& \ref{ass:LisNSCVX} and \ref{ass:Lipschitz} & reset (at $j = 1$).\\
			\hline
			\multirow{4}{*}{HHA} & Assumptions & $\min \left\{1,\exp\left(-\Omega \left(t - \theta\right)\right)\right\}$, \\
			& \ref{ass:HHA_Assumps}, \ref{ass:Timeout}, & where $\Omega$ and $\theta$ are\\
			& and \ref{ass:PL} & defined via \eqref{eqn:OmegaTheta}.\\
			\hline
		\end{tabular}
		\caption{A comparison of assumptions and convergence rates, for HAND-1 from \cite{poveda2019inducing}, HHA from \cite{teel2019first}, and the hybrid closed-loop algorithm $\HS$.}		
		\label{table:AssumpsProperties}
	\end{center}
\end{table}

Next, we compare $\HS_0$, $\HS_1$, $\HS$, HAND-1, and HHA in simulation. To compare these algorithms, the choice of objective function, parameter values, and initial conditions are as follows. We use the objective function $L(\xp_1) := \xp_1^2$, the gradient of which is Lipschitz continuous with $M = 2$, and which has a single minimizer at $\minSet = 0$. This choice of objective function is made so that we can easily tune $\lambda$, as described in Section \ref{sec:DesignOfLambda}. We arbitrarily chose the heavy ball parameter value $\gamma = \frac{2}{3}$ and we tuned $\lambda$ to $40$ by choosing a value arbitrarily larger than $2\sqrt{a_1}$, where $a_1$ comes from Section \ref{sec:DesignOfLambda}, and gradually increasing it until there is no overshoot in the hybrid algorithm. Since the gain coefficient of $\xp_1$ for HHA is fixed at $-1$, then to ensure that the gain coefficients of $\xp_1$ for $\HS$ and HHA are the same, the only value possible for $\zeta$, for the objective function $L$ chosen, is $\sqrt{2}$, for Nesterov's algorithm. In Addition, since $\xp_2$ does not appear in the dynamics of $\dot{\xp}_2$ of HHA in \eqref{eqn:HHA}, we could not compare gain coefficients for $\xp_2$ for both $\HS$ and HHA. For the remaining HHA parameters, we chose $\bar{M} = 2$, since this satisfies $\bar{M} \geq M$. To ensure Assumption \ref{ass:Timeout} is satisfied, we chose $\bar{T} := \frac{\pi}{2\sqrt{1}}$.
Given $\zeta = \sqrt{2}$, the HAND-1 parameters $c_1 = 0.25$ and $T_{\min} = \frac{1+\sqrt{13}}{2}$ are chosen such that the resulting gain coefficients for $\xp_1$ and $\xp_2$ are the same for both $\HS$ and HAND-1, so that these algorithms are compared on equal footing. The remaining HAND-1 parameters, we chose $r = 51$ and $\delta_{\text{med}} = 8650$ such that $T_{\text{med}} \geq \sqrt{\frac{B}{\delta_{\text{med}}}} + T_{\min} > 0$, and we chose $T_{\max} = T_{\text{med}} + 1$ to ensure resets happen at the proper times. The uniting algorithm parameters are $c_0 = 7000$, $c_{1,0} \approx 1354.025$, $\varepsilon_0 = 10$, $\varepsilon_{1,0} = 5$, and $\alpha = 1$, which yield the values $\tilde{c}_0 = 10$, $\tilde{c}_{1,0} = 5$, $d_0 \approx 6933.3$, and $d_{1,0} \approx 1304.0$, which are calculated via \eqref{eqn:UTilde0SetEquations} and \eqref{eqn:TTilde10SetEquations}. These values are chosen for proper tuning of the algorithm, in order to get nice performance, and for exploiting the properties of Nesterov's method as long as we want. Initial conditions for $\HS$ are $\xp_1(0,0) = 50$, $\xp_2(0,0) = 0$, $\xlogic(0,0) = 1$, and $\tau(0,0) = 0$ Initial conditions for HAND-1 are $\xp_1(0,0) = \xp_2(0,0) = 50$ and $\tau(0,0) = T_{\min}$. Initial conditions for HHA are $\xp_1(0,0) = 50$, $\xp_2(0,0) = 0$, and $\tau(0,0) = 0$.

\begin{figure}[thpb]
	\centering
	\setlength{\unitlength}{1.0pc} 
	%
	%
	\begin{picture}(24,19)(0,0)
		\footnotesize
%
%
		\put(0,0.5){\includegraphics[trim={0.5cm 0.2cm 1cm 0.3cm},clip,width=24\unitlength]{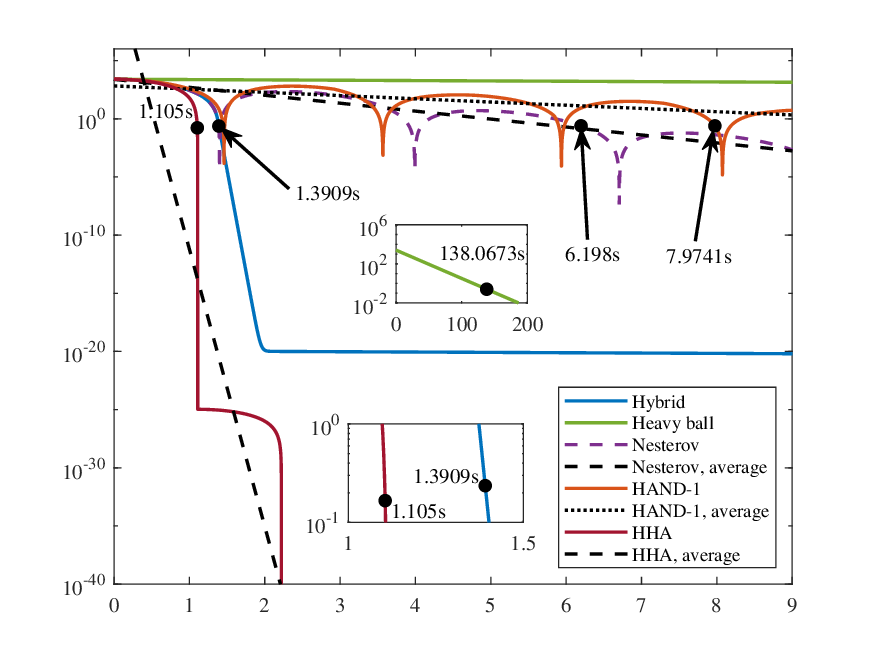}}
		\put(0,8){\rotatebox{90}{$L(\xp_1) - L^*$}}
		\put(12.3,0.5){$t[s]$}
	\end{picture}
	\caption{A comparison of the evolution of $L$ over time for Nesterov's method in \eqref{eqn:MJODE_ZetaNum_TR}, heavy ball, HAND-1 from \cite{poveda2019inducing}, HHA from \cite{teel2019first}, and our proposed uniting algorithm, for a function $L(\xp_1) := \xp_1^2$, with a single minimizer at $\xp_1^* = 0$. Nesterov's method, shown in purple, settles to within $1\%$ of $\xp_1^*$ in about 6.2 seconds. The heavy ball algorithm, shown in green, settles to within $1\%$ of $\xp_1^*$ in about 138.1 seconds. HAND-1, shown in orange, settles to within $1\%$ of $\xp_1^*$ in about 8.0 seconds. HHA, shown in red, settles to within $1\%$ of $\xp_1^*$ in about 1.1 seconds. The hybrid closed-loop system $\HS$, shown in blue, settles to within $1\%$ of $\minSet$ in about 1.4 seconds.}	
	\label{fig:HHAComparison}
\end{figure}

\begin{table}[thpb]	
	\begin{center}
		\begin{tabular}{|c|c|c|}
			\hline
			& { Time to} & { $\%$}\\
			Algorithm & { converge (s)} & { improvement} \\
			\hline
			\hline
			$\HS$ & $1.390$ & -- \\
			\hline
			$\HS_0$ & $138.1$ & $99.0$ \\
			\hline
			$\HS_1$ & $6.191$ & $77.5$ \\
			\hline
			HAND-1 & $7.974$  & $82.6$ \\
			\hline
			HHA & $1.105$  & -- \\
			\hline
		\end{tabular}
		\caption{Average times for which $\HS$, $\HS_0$, $\HS_1$, HAND-1, and HHA settle to within $1\%$ of $\minSet$, and the percent improvement of $\HS$ over each algorithm. Percent improvement is calculated via \eqref{eqn:PercentImprovement}. The objective function used for this table is $L(\xp_1) := \xp_1^2$.}		
		\label{table:TimePercentImpHHA}
	\end{center}
\end{table}

Figure \ref{fig:HHAComparison} and Table \ref{table:TimePercentImpHHA} show the time that each algorithm takes to settle within\footnote{Code at \texttt{github.com/HybridSystemsLab/UnitingComparison}} $1\%$ of $\minSet$. Table \ref{table:TimePercentImpHHA} shows the percent improvement $\HS$ over $\HS_0$, $\HS_1$, and HAND-1, which is calculated using \eqref{eqn:PercentImprovement}. Since HHA converges exponentially, as shown in the bound in \eqref{eqn:ConvergenceBoundHHA}, such a rate explains why HHA is faster than $\HS$ in Figure \ref{fig:HHAComparison} and Table \ref{table:TimePercentImpHHA}. The hybrid closed-loop algorithm $\HS$ converges faster than $\HS_0$, $\HS_1$, and HAND-1, however. 
\end{example}

\begin{example} \label{ex:tradeOff}
	This example explores the trade-off that results from using different values of $\zeta > 0$ for the uniting algorithm. Particularly, for $\zeta = 1$, we first compare the uniting algorithm in simulation with the individual optimization algorithms $\HS_0$, $\HS_1$, and the HAND-1 algorithm from \cite{poveda2019inducing}, using the same objective function as in Example \ref{ex:NSC}, and next we compare the resulting solutions with those in Table \ref{table:TimePercentImpNSC}. Recall that the objective function in Example \ref{ex:NSC} is $L(\xp_1) := \xp_1^2$, the gradient of which is Lipschitz continuous with $M = 2$, and which has a single minimizer at $\minSet = 0$. Since the gain coefficient of $\nabla L$ is proportional to $\zeta^2$, we choose different parameters for the HAND-1 algorithm for the simulation depicted in\footnote{Code found at same link as in Footnote \IfConf{\ref{foot:Tradeoff}}{\ref{foot:Tradeoff_TR}}} Figure \ref{fig:NSCLargerZetaTR}, so that the gain coefficients of $\xp_1$ and $\xp_2$ are the same for HAND-1 and $\HS$ in this simulation. Namely, given $\zeta = 1$, for HAND-1 we choose $T_{\min} = 3$ and $c_1 = 0.25$. For the other HAND-1 parameters, we choose $r = 51$ and $\delta_{\text{med}} = 4010$ such that $T_{\text{med}} \geq \sqrt{\frac{B}{\delta_{\text{med}}}} + T_{\min} > 0$, and we again choose $T_{\max} = T_{\text{med}} + 1$ to ensure resets happen at the proper times.
	We arbitrarily choose $\gamma = \frac{2}{3}$, and we tuned $\lambda$ to $40$ by choosing a value arbitrarily larger than $2\sqrt{a_1}$ and gradually increasing until there was no overshoot in the hybrid algorithm.
	The uniting algorithm parameters are $c_0 = 320$, $c_{1,0} \approx 271.584$, $\varepsilon_0 = 10$, $\varepsilon_{1,0} = 5$, and $\alpha = 1$, which yield the values $\tilde{c}_0 = 10$, $\tilde{c}_{1,0} = 5$, $d_0 \approx 253.333$, and $d_{1,0} \approx 234.084$, which are calculated via \eqref{eqn:UTilde0SetEquations} and \eqref{eqn:TTilde10SetEquations}. These values are chosen for proper tuning of the algorithm, in order to get nice performance, and for exploiting the properties of Nesterov's method as long as we want. Initial conditions for $\HS$ are $\xp_1(0,0) = 50$, $\xp_2(0,0) = 0$, $\xlogic(0,0) = 1$, and $\tau(0,0) = 0$, and for HAND-1 are $\xp_1(0,0) = \xp_2(0,0) = 50$ and $\tau(0,0) = T_{\min}$.
	
	First, we compare solutions to each algorithm within Figure \ref{fig:NSCLargerZetaTR} itself.
	Table \ref{table:TimePercentImpNSCZeta1} shows the time that each algorithm takes to settle within $1\%$ of $\minSet$, averaged over solutions starting from ten different values of $\xp_1(0,0)$ (listed in the first column of Table \ref{table:PercentImprovementNSCTR}), and the percent improvement of $\HS$ over $\HS_0$, $\HS_1$, and HAND-1, which is calculated using \eqref{eqn:PercentImprovement}.
	While the closed-loop algorithm $\HS$ still converges faster than all the other algorithms in Figure \ref{fig:NSCLargerZetaTR} and Table \ref{table:TimePercentImpNSCZeta1}, the improvement over $\HS_0$, $\HS_1$, and HAND-1 is smaller than it is in Table \ref{table:TimePercentImpNSC}. 
	
	\begin{table}[thpb]
		\begin{center}
			\begin{tabular}{|c|c|c|}
				\hline
				Algorithm & Average time to converge (s) & Average $\%$ improvement  \\
				\hline
				\hline
				$\HS$ & $2.387$ & -- \\
				\hline
				$\HS_0$ & $138.066$ & $98.3$ \\
				\hline
				$\HS_1$ & $8.782$ & $72.8$ \\
				\hline
				HAND-1 & $14.343$  & $83.4$ \\
				\hline
			\end{tabular}
			\caption{Times for which $\HS$, $\HS_0$, $\HS_1$, and HAND-1 settle to within $1\%$ of $\minSet$, and percent improvement of $\HS$ over each algorithm, as shown in Figure \ref{fig:NSCLargerZetaTR}. Percent improvement is calculated via \eqref{eqn:PercentImprovement}. The objective function used for this table is $L(\xp_1) := \xp_1^2$.}
			\label{table:TimePercentImpNSCZeta1}
		\end{center}
	\end{table}
	
	Next, we compare solutions using $\zeta = 1$, in Figure \ref{fig:NSCLargerZetaTR}, with solutions using $\zeta = 2$, in Table \ref{table:TimePercentImpNSC}. Since $\zeta > 0$ scales time in solutions to \IfConf{\eqref{eqn:MJODE_ZetaNum}}{\eqref{eqn:MJODE_ZetaNum_TR}}, Then smaller values of $\zeta$ result in slower settling to within $1\%$ of $\minSet$ for $\HS_1$ with less frequent oscillations, as seen in Figure \ref{fig:NSCLargerZetaTR} with $\zeta = 1$ (about $8.8$ seconds), while larger values of $\zeta$ result in settling to within $1\%$ of $\minSet$ for $\HS_1$ faster, with more frequent oscillations, as seen in Figure \ref{fig:MotivationalPlotTR} and Table \ref{table:TimePercentImpNSC} with $\zeta = 2$ (about $4.5$ seconds). For the uniting algorithm, this translates to faster settling to within $1\%$ of $\minSet$ with $\zeta = 2$ (about $0.8$ seconds), in Figure \ref{fig:MotivationalPlotTR} and Table \ref{table:TimePercentImpNSC}, compared with slower settling to within $1\%$ of $\minSet$ with $\zeta = 1$ (about $2.4$ seconds), in Figure \ref{fig:NSCLargerZetaTR}, but with no oscillations, in both cases, due to the switch to $\HS_0$. In both Figure \ref{fig:NSCLargerZetaTR}, and Table \ref{table:TimePercentImpNSC}, the uniting algorithm converges more quickly than the HAND-1 algorithm, when both algorithms are tuned to have the same gain coefficients for the $\xp_1$ and $\xp_2$ terms.  Although larger $\zeta$ results in faster convergence, the trade-off is that even though the $\xp_2$ (velocity) term generally reduces quickly as it approaches the neighborhood of the minimizer for any size of $\zeta$, the $\xp_2$ still ends up relatively larger near the minimizer than it is when $\zeta$ is smaller. The consequence is that, when $\zeta$ is larger, $d_{1,0}$ needs to be set much larger so that the uniting algorithm can still make the switch to $\HS_0$ at the proper time. This also means that $c_{1,0}$ needs to be set much larger, due to the definition of $d_{1,0}$ in \eqref{eqn:TTilde10SetEquations}. Additionally, $c_0$ and $d_0$, also need to be set larger to ensure the algorithm still has adequate hysteresis. Recall that, in Example \ref{ex:NSC}, for $\zeta = 2$, we have the parameter values $c_0 = 7000$, $c_{1,0} \approx 6819.676$, $d_0 = 6933$, and $d_{1,0} = 6744$, which are quite large, while for the simulation shown in Figure \ref{fig:NSCLargerZetaTR} these same parameters have much smaller values, as listed in the second paragraph of this example. 
\end{example}}
\section{Proof of Theorem \ref{thm:GASNestNSC}}
\label{sec:ProofMainResult}
\IfConf{This section provides a proof of Theorem \ref{thm:GASNestNSC} from Section \ref{sec:MainResult}. Section \ref{sec:PropertiesOfH0} establishes UGAS of $\{\minSet\} \times \{0\}$ and an exponential convergence rate for $\HS_0$. Section \ref{sec:PropertiesOfH1} establishes UGAS of $\{\minSet\} \times \{0\} \times \reals_{\geq 0}$ and a convergence rate $\frac{1}{(t+2)^2}$ for $\HS_1$. Section \ref{sec:UGAS} uses the properties in Sections \ref{sec:PropertiesOfH0} and \ref{sec:PropertiesOfH1} to establish UGAS of $\mathcal{A}$, defined via \eqref{eqn:SetOfMinimizersHS-NSCVX}, for $\HS$. Finally, Section \ref{sec:ConvRateH} proves the convergence rate of $\HS$ using the convergence rates of the individual closed-loop algorithms $\HS_0$ and $\HS_1$ established in Sections \ref{sec:PropertiesOfH0} and \ref{sec:PropertiesOfH1}, respectively.}{
	This section provides a proof of Theorem \ref{thm:GASNestNSC} from Section \ref{sec:MainResult}. The proof consists of the following steps.
	\begin{itemize}[leftmargin=*]
		\item Section \ref{sec:PropertiesOfH0} establishes UGAS of $\{\minSet\} \times \{0\}$, and an exponential convergence rate for the closed-loop algorithm $\HS_0$;
		\item Section \ref{sec:PropertiesOfH1} establishes UGAS of $\{\minSet\} \times \{0\} \times \reals_{\geq 0}$, and a convergence rate $\frac{1}{(t+2)^2}$ for the closed-loop algorithm $\HS_1$;
		\item Section \ref{sec:UGAS} uses the properties in Sections \ref{sec:PropertiesOfH0} and \ref{sec:PropertiesOfH1} and a proof-by-contradiction to establish UGAS of $\mathcal{A}$, defined via \eqref{eqn:SetOfMinimizersHS-NSCVX}, for $\HS$;
		\item Section \ref{sec:ConvRateH} proves the convergence rate of $\HS$ using the convergence rates of the individual closed-loop algorithms $\HS_0$ and $\HS_1$ established in Sections \ref{sec:PropertiesOfH0} and \ref{sec:PropertiesOfH1}, respectively.
	\end{itemize}}

	\subsection{Properties of $\HS_0$}
	\label{sec:PropertiesOfH0}
	The following result establishes that the closed-loop algorithm $\HS_0$ \IfConf{}{in \eqref{eqn:H0}} has the set $\{\minSet\} \times \{0\}$ UGAS. To prove it, we use an invariance principle. \IfConf{Its proof is in \cite{dhustigs2022unitingNSC}.}{}
	
	\IfConf{\begin{prop}[UGAS of $\{\minSet\} \times \{0\}$ for $\HS_0$]}{\begin{proposition}(UGAS of $\{\minSet\} \times \{0\}$ for $\HS_0$)} \label{prop:GAS-HBF}
		Let $L$ satisfy Assumptions \ref{ass:LisNSCVX}, \ref{ass:QuadraticGrowth}, and \ref{ass:Lipschitz}. For each $\lambda > 0$ and $\gamma > 0$, the set $\{\minSet\} \times \{0\}$ is UGAS for the closed-loop algorithm $\HS_0$ in \eqref{eqn:H0}. 
	\IfConf{\end{prop}}{\end{proposition}} \IfConf{

\vspace{-0.3cm}

 }{}\IfConf{}{\begin{proof}
		By Proposition \ref{prop:Existence}, each maximal solution to the closed-loop algorithm $\HS_0$, defined via \eqref{eqn:H0}, is bounded, complete, and unique. Recall that, in the proof of Proposition \ref{prop:Existence}, it was shown that $V_0$ in \eqref{eqn:LyapunovHBF} satisfies \IfConf{$\dot{V}_0(\xp) \leq -\lambda \left|\xp_2\right|^2 \leq 0$}{\eqref{eqn:VdotHBF}} for all $\xp \in \reals^{2n}$, since $\lambda$ is positive. Therefore, by an application of \IfConf{\cite[Theorem~A.3]{dhustigs2022unitingNSC}}{Theorem \ref{thm:hybrid Lyapunov theorem}}, since $\gamma > 0$ and $\lambda > 0$, the set $\{\minSet\} \times \{0\}$ is stable for the closed-loop algorithm $\HS_0$. 
		Since by Lemma \ref{lemma:HBC} $\HS_0$ satisfies the hybrid basic conditions, then, using the invariance principle in \IfConf{\cite[Theorem~A.6]{dhustigs2022unitingNSC}}{Theorem \ref{thm:HybridInvariancePrinciple}}, every maximal solution that is complete and bounded approaches the largest weakly invariant set for $\HS_0$ in \eqref{eqn:H0} that is contained in 
		\begin{equation}\label{eqn:LaSalleSetVHBF}
			\defset{\xp \in \reals^{2n}\!\!}{\!\!\dot{V}_0(\xp) = 0\!\!} \cap \defset{\xp \in \reals^{2n}\!\!}{\!\!V_0(\xp) = r\!\!}, \ r \geq 0.
		\end{equation}
		Such a set is nonempty only when $r = 0$ and, precisely, is equal to $\{\minSet\} \times \{0\}$. This property can be seen by noticing that 
		$\defset{\xp \in \reals^{2n}}{\dot{V}_0(\xp) = 0\!\!}
		= 
		\defset{\xp \in \reals^{2n}}{\xp_2 = 0\!\!}$, and that after setting $\xp_2$ to zero in
		\eqref{eqn:H0} we obtain $\matt{\dot{\xp}_1\\ 0\\} = \matt{0\\ -\gamma \nabla L(\xp_1)}$.
		For any solution to this system, its $\xp_1$ component satisfies $0 = \gamma \nabla L(\xp_1)$, which, since $\gamma > 0$ and since $\nabla L(\xp_1) = 0$ only when $\xp_1$ is the minimizer of $L$,		
		leads to $\xp_1 = \minSet$. Then, the only maximal solution that starts and stays in \eqref{eqn:LaSalleSetVHBF} is the solution from $\{\minSet\} \times \{0\}$, for which $r = 0$. Then, every bounded and complete solution to the closed-loop algorithm $\HS_0$ converges to $\{\minSet\} \times \{0\}$. The arguments above involving the Lyapunov theorem in \IfConf{\cite[Theorem~A.3]{dhustigs2022unitingNSC}}{Theorem \ref{thm:hybrid Lyapunov theorem}} and the invariance principle in \IfConf{\cite[Theorem~A.6]{dhustigs2022unitingNSC}}{Theorem \ref{thm:HybridInvariancePrinciple}} yield global pre-asymptotic stability of $\{\minSet\} \times \{0\}$ for $\HS_0$. Since by Proposition \ref{prop:Existence}, each maximal solution to $\HS_0$ is complete, then $\{\minSet\} \times \{0\}$ is globally asymptotically stable for the closed-loop algorithm $\HS_0$. Since $\HS_0$ satisfies the hybrid basic conditions by Lemma \ref{lemma:HBC}, then, by \IfConf{\cite[Theorem~A.4]{dhustigs2022unitingNSC}}{Theorem \ref{thm:GASImpliesUGAS}}, $\{\minSet\} \times \{0\}$ is UGAS for $\HS_0$.\IfConf{\hfill{} \qed}{}
	\end{proof}} 
%
%
Next, we establish the convergence rate of the closed-loop algorithm $\HS_0$. To do so, we use the following Lyapunov function, proposed in \cite[Lemma~4.2]{sebbouh2020convergence}\IfConf{:}{, for $\HS_0$:} \IfConf{\vspace{-0.2cm}}{}
	\begin{align} \label{eqn:AlternateVHBF}
		\IfConf{V(\xp) := & \gamma \left(L(\xp_1) - L^*\right) + \frac{1}{2}\left|\psi(\xp_1 - \xp_1^*) + \xp_2\right|^2 \nonumber\\
		& + \frac{\nu}{2}\left|\xp_1 - \xp_1^*\right|^2}{V(\xp) := \gamma \left(L(\xp_1) - L^*\right) + \frac{1}{2}\left|\psi(\xp_1 - \xp_1^*) + \xp_2\right|^2 + \frac{\nu}{2}\left|\xp_1 - \xp_1^*\right|^2}
	\end{align}
	where, given $\lambda > 0$, $\psi > 0$ is chosen such that $\nu := \psi \left(\psi - \lambda\right) < 0$. When $L$ satisfies Assumption \ref{ass:LisNSCVX}, the following lemma, which is a version of \cite[Lemma~4.2]{sebbouh2020convergence} tailored for 
	the unperturbed heavy ball algorithm in \eqref{eqn:H0}, gives an upper bound on the change of \IfConf{$V$}{the Lyapunov function} in \eqref{eqn:AlternateVHBF}. \IfConf{Its proof is in \cite{dhustigs2022unitingNSC}.}{}
	\IfConf{\begin{lem}}{\begin{lemma}} \label{lemma:IntermediateDotV}
		Let $L$ satisfy Assumption \ref{ass:LisNSCVX}, and let $\lambda > 0$ and $\gamma > 0$, which come from $\HS_0$ in \eqref{eqn:H0}, be given. For each $\psi > 0$ such that $\nu := \psi (\psi - \lambda) < 0$, the following bound is satisfied for each $\xp \in \reals^{2n}$: \IfConf{$\dot{V}(\xp) \leq -\psi \left(a(\xp_1) + 2\nu c(\xp_1)\right) + 2(\psi - \lambda)b(\xp)$,}{
		\begin{equation} \label{eqn:dotVHBFAlternateLemma}
			\dot{V}(\xp) \leq -\psi \left(a(\xp_1) + 2\nu c(\xp_1)\right) + 2(\psi - \lambda)b(\xp)
		\end{equation}}
		where $V$ is defined in \eqref{eqn:AlternateVHBF}, $a(\xp_1) := \gamma \left(L(\xp_1) - L^*\right)$, $b(\xp) := \frac{1}{2}\left|\psi(\xp_1 - \xp_1^*) + \xp_2\right|^2$, and $c(\xp_1) := \frac{1}{2}\left| \xp_1 - \xp_1^* \right|^2$. 
	\IfConf{\end{lem}}{\end{lemma}} 
		
	\IfConf{}{\begin{proof} Since $L$ is $\mathcal{C}^1$, convex, and has a single minimizer $\xp_1^*$, and
		since \IfConf{$\nabla V(\xp) = \left[\gamma \nabla L(\xp_1) + \psi \left( \psi \left( \xp_1 - \xp_1^* \right) + \xp_2 \right) + \nu \left( \xp_1 - \xp_1^* \right) \right.$\\$\left. \psi \left( \xp_1 - \xp_1^* \right) + \xp_2 \right]$}{$\nabla V(\xp) =$\\$\left[\gamma \nabla L(\xp_1) + \psi \left( \psi \left( \xp_1 - \xp_1^* \right) + \xp_2 \right) + \nu \left( \xp_1 - \xp_1^* \right) \ \; \psi \left( \xp_1 - \xp_1^* \right) + \xp_2 \right]$}, then we evaluate the derivative of $V$, defined via \eqref{eqn:AlternateVHBF}, using the map $\xp \mapsto F_P(\hczero(\hp_0(\xp)))$, where $F_P$ is defined via \eqref{eqn:HBFplant-dynamicsTR}, $\hczero$ is defined in \eqref{eqn:StaticStateFeedbackLawLocal}, and $\hp_0$ is defined via \eqref{eqn:H0H1NSCNesterovHBF}. For each $\xp \in \reals^{2n}$, we obtain \IfConf{
		\begin{align} \label{eqn:AlternateDotVHBF}
			\dot{V}(\xp) = & \left\langle \nabla V(z),F_P(\hczero(\hp_0(\xp))) \right\rangle \\
			& = -\gamma \psi \left\langle \nabla L(\xp_1),\xp_1 - \xp_1^* \right\rangle\nonumber\\
			& + \left(\nu + \psi(\psi - \lambda) \right) \left\langle \xp_2, \xp_1 - \xp_1^*\right\rangle + (\psi - \lambda)\left|\xp_2\right|^2.\nonumber
		\end{align}
		}{\begin{align} \label{eqn:AlternateDotVHBF_TR}
			\dot{V}(\xp) = & \left\langle \nabla V(z),F_P(\hczero(\hp_0(\xp))) \right\rangle = \left\langle \nabla V(z), \matt{\xp_2 \\ \hczero(\hp_0(\xp))} \right\rangle \\
			= & \gamma \left\langle \nabla L(\xp_1), \xp_2 \right\rangle + \psi \left\langle \xp_2, \psi \left( \xp_1 - \xp_1^* \right) + \xp_2 \right\rangle + \nu \left\langle \xp_2, \xp_1 - \xp_1^* \right\rangle \nonumber\\
			& - \lambda \left\langle \xp_2, \psi \left( \xp_1 - \xp_1^* \right) + \xp_2 \right\rangle - \gamma \left\langle \nabla L(\xp_1), \psi \left( \xp_1 - \xp_1^* \right) + \xp_2 \right\rangle \nonumber\\
			= & -\gamma \psi \left\langle \nabla L(\xp_1),\xp_1 - \xp_1^* \right\rangle + \left(\nu + \psi(\psi - \lambda) \right) \left\langle \xp_2, \xp_1 - \xp_1^*\right\rangle + (\psi - \lambda)\left|\xp_2\right|^2.\nonumber
		\end{align}}
		Note that $\left|\psi\left( \xp_1 - \xp_1^* \right) + \xp_2\right|^2 = \left|\xp_2\right|^2 + 2\psi \left\langle \xp_2, \xp_1 - \xp_1^* \right\rangle + \psi^2 \left| \xp_1 - \xp_1^* \right|^2$, from where we obtain \IfConf{$\left|\xp_2\right|^2 = $\\$\left|\psi\left( \xp_1 - \xp_1^* \right) + \xp_2\right|^2 - 2\psi \left\langle \xp_2, \xp_1 - \xp_1^* \right\rangle - \psi^2 \left| \xp_1 - \xp_1^* \right|^2$}{$\left|\xp_2\right|^2 =  \left|\psi\left( \xp_1 - \xp_1^* \right) + \xp_2\right|^2 - 2\psi \left\langle \xp_2, \xp_1 - \xp_1^* \right\rangle - \psi^2 \left| \xp_1 - \xp_1^* \right|^2$}.
		Substituting the expression for $\left|\xp_2\right|^2$ into \IfConf{\eqref{eqn:AlternateDotVHBF}}{\eqref{eqn:AlternateDotVHBF_TR}}, we arrive at, for all $\xp \in \reals^{2n}$, \IfConf{
		\begin{align} \label{eqn:Substituting}
			\dot{V}(\xp) = & -\gamma\psi \left\langle \nabla L(\xp_1),\xp_1 - \xp_1^* \right\rangle + 2(\psi - \lambda)b(\xp)\nonumber\\
			& - 2\psi \nu c(\xp_1)
		\end{align}
		}{\begin{align} \label{eqn:Substituting_TR}
			\dot{V}(\xp) = & -\gamma \psi \left\langle \nabla L(\xp_1),\xp_1 - \xp_1^* \right\rangle + (\psi - \lambda) \left|\psi\left( \xp_1 - \xp_1^* \right) + \xp_2\right|^2\nonumber\\
			& + \left(\nu - \psi(\psi - \lambda) \right) \left\langle \xp_2, \xp_1 - \xp_1^*\right\rangle - \psi^2 (\psi - \lambda) \left| \xp_1 - \xp_1^* \right|^2\nonumber\\
			= & -\gamma\psi \left\langle \nabla L(\xp_1),\xp_1 - \xp_1^* \right\rangle + 2(\psi - \lambda)b(\xp) - 2\psi \nu c(\xp_1)
		\end{align}}
		since $\nu = \psi(\psi - \lambda)$,
		where $b(\xp) = \frac{1}{2}\left|\psi\left( \xp_1 - \xp_1^* \right) + \xp_2\right|^2$ and $c(\xp_1) = \frac{1}{2}\left| \xp_1 - \xp_1^* \right|^2$. 
		Since $L$ is $\mathcal{C}^1$, convex, and has a unique minimizer by Assumption \ref{ass:LisNSCVX}, then using the definition of convexity in Footnote \ref{foot:Convexity} with $u_1 = \xp_1^*$ and $w_1 = \xp_1$, we get $- \left(L(\xp_1) - L^*\right) \geq - \left\langle \nabla L(\xp_1), \xp_1 - \xp_1^* \right\rangle$.
		Substituting it into \IfConf{\eqref{eqn:Substituting}}{\eqref{eqn:Substituting_TR}} yields, for all $\xp \in \reals^{2n}$, $\dot{V}(\xp) \leq -\psi a(\xp_1) + 2(\psi - \lambda)b(\xp) - 2\psi \nu c(\xp_1)$,
		where $a(\xp_1) = \gamma \left(L(\xp) - L^*\right)$, and \eqref{eqn:dotVHBFAlternateLemma} is satisfied.\IfConf{\hfill{} \qed}{}
	\end{proof}} 
%
%
We employ Lemma 5.2 to show that when $L$ satisfies Assumptions \ref{ass:LisNSCVX} and \ref{ass:QuadraticGrowth}, the convergence rate of the closed-loop algorithm $\HS_0$ in \eqref{eqn:H0} is exponential. This is supported by the following proposition, which is a version of \cite[Theorem~3.2]{sebbouh2020convergence} tailored for the unperturbed heavy ball algorithm $\HS_0$ in \eqref{eqn:H0}. \IfConf{Its proof is in \cite{dhustigs2022unitingNSC}.}{}
	\IfConf{\begin{prop}[Convergence rate for $\HS_0$]}{\begin{proposition}(Convergence rate for $\HS_0$)} \label{prop:HBFConvergenceRate} 
		Let $L$ satisfy Assumptions \ref{ass:LisNSCVX} and \ref{ass:QuadraticGrowth}, let $\alpha > 0$ come from \IfConf{Assumption \ref{ass:QuadraticGrowth}}{\eqref{eqn:QuadraticGrowth}}, and let $\lambda > 0$ and $\gamma > 0$ come from $\HS_0$ in \eqref{eqn:H0}. For each $m \in (0,1)$ such that $\psi := \frac{m\alpha\gamma}{\lambda} > 0$ and $\nu := \psi (\psi - \lambda) < 0$,
		each maximal solution $t \mapsto \xp(t)$ to the closed-loop algorithm $\HS_0$  satisfies
		\begin{equation} \label{eqn:ConvergenceRateHBFNSCVX}
			\IfConf{L(\xp_1(t)) - L^* = \mathcal{O}\left(\exp \left(-(1-m)\psi t\right)\right)}{L(\xp_1(t)) - L^* = \mathcal{O}\left(\exp \left(-(1-m)\psi t\right)\right) \quad \forall t \in \dom \xp \ (= \reals_{\geq 0}).}
		\end{equation}
	\IfConf{for all $t \in \dom \xp \ (= \reals_{\geq 0})$.}{}
	\IfConf{\end{prop}}{\end{proposition}} 
	
	\IfConf{}{\begin{proof}
			By Lemma \ref{lemma:IntermediateDotV}, the bound in \eqref{eqn:dotVHBFAlternateLemma} is satisfied for $V$ in \eqref{eqn:AlternateVHBF} for each $\xp \in \reals^{2n}$ since, by Assumption \ref{ass:LisNSCVX}, $L$ is $\mathcal{C}^1$, convex, and has a single minimizer $\xp_1^*$. Then, since $\psi = \frac{m\alpha\gamma}{\lambda} > 0$ is such that $\nu = \psi (\psi - \lambda) < 0$ and $c$ is nonnegative, this leads to
		\begin{equation} \label{eqn:NuIsNeg}
			V(\xp) = a(\xp_1) + b(\xp) + \nu c(\xp_1) \leq a(\xp_1) + b(\xp) \quad \forall \xp \in \reals^{2n}
		\end{equation}
		where $a$, $b$, and $c$ are defined below \eqref{eqn:dotVHBFAlternateLemma}.
		By Assumption \ref{ass:QuadraticGrowth}, $L$ has quadratic growth away from $\minSet$.  Therefore, we have, for all $\xp \in \reals^{2n}$, \IfConf{
			
			\vspace{-1cm}
			
		\begin{align} \label{eqn:InequalityFollowingFromQG}
			a(\xp_1) + 2\nu c(\xp_1) = a(\xp_1) - 2\left|\nu\right| c(\xp_1) \geq \left(1 - \frac{\left|\nu\right|}{\alpha \gamma}\right) a(\xp_1).
		\end{align}
	
		\vspace{-0.9cm}
	
		}{\begin{align} \label{eqn:InequalityFollowingFromQG_TR}
			a(\xp_1) + 2\nu c(\xp_1) = & a(\xp_1) - 2\left|\nu\right| c(\xp_1)  = \gamma \left(L(\xp_1) - L^*\right) - \left|\nu\right| \left| \xp_1 - \xp_1^* \right|^2 \\
			\geq &  \gamma \left(L(\xp_1) - L^*\right) - \frac{\left|\nu\right|\left(L(\xp_1) - L^*\right)}{\alpha}
			= \left(1 - \frac{\left|\nu\right|}{\alpha \gamma}\right) a(\xp_1).\nonumber
		\end{align}}
		Observe that, for each $m \in (0,1)$ such that $\psi = \frac{m\alpha\gamma}{\lambda} > 0$ and $\nu = \psi (\psi - \lambda) < 0$,
		we have \IfConf{$\left|\nu\right| = \psi\left(\lambda - \psi\right) \leq \lambda \psi = m \alpha \gamma$.}{
			\begin{equation} \label{eqn:NuUpperBound}
				\left|\nu\right| = \psi\left(\lambda - \psi\right) \leq \lambda \psi = m \alpha \gamma
			\end{equation}}
		It follows from \IfConf{this upper bound on $\left|\nu\right|$ and}{} \IfConf{\eqref{eqn:InequalityFollowingFromQG}}{\eqref{eqn:InequalityFollowingFromQG_TR} and \eqref{eqn:NuUpperBound}} that
		\begin{equation} \label{eqn:LowerBoundAandC}
			a(\xp_1) + 2\nu c(\xp_1) \geq (1-m)a(\xp_1) \IfConf{\vspace{0.2cm}}{}
		\end{equation} 
		for all $\xp \in \reals^{2n}$. Noticing that from \eqref{eqn:NuIsNeg} we have $a(\xp_1) + \frac{1}{(1-m)}b(\xp) \geq a(\xp_1) + b(\xp) \geq V(\xp)$, substituting \eqref{eqn:LowerBoundAandC} into \eqref{eqn:dotVHBFAlternateLemma} we have \IfConf{ \vspace{-0.15cm}
			\begin{align} \label{eqn:dotVContainingV}
				\dot{V}(\xp) & \leq - (1 - m) \psi a(\xp_1) + \psi b(\xp) + (\psi - 2\lambda) b(\xp) \nonumber\\
				& \leq - (1-m)\psi a(\xp_1) + \psi b(\xp)\nonumber\\
				& \leq -(1-m)\psi V(\xp)
			\end{align}
		}{\begin{align} \label{eqn:dotVContainingV_TR}
			\dot{V}(\xp) & \leq - (1-m) \psi a(\xp_1) + 2(\psi - \lambda)b(\xp) \nonumber\\
			& \leq - (1 - m) \psi a(\xp_1) + \psi b(\xp) + (\psi - 2\lambda) b(\xp) \nonumber\\
			& \leq - (1-m)\psi a(\xp_1) + \psi b(\xp)\nonumber\\
			& \leq -(1-m)\psi \left(a(\xp_1) + \frac{1}{(1-m)}b(\xp)\right)\nonumber\\
			& \leq -(1-m)\psi V(\xp)
		\end{align}}
		for all $\xp \in \reals^{2n}$. The \IfConf{second}{third} inequality comes from the fact that we choose $\psi = \frac{m\alpha\gamma}{\lambda} > 0$ such that $\psi - \lambda < 0$ and, consequently, $\psi - 2\lambda < 0$. Applying Gr\"{o}nwall's inequality to \IfConf{\eqref{eqn:dotVContainingV}}{\eqref{eqn:dotVContainingV_TR}} shows that every maximal solution $t \mapsto \xp(t)$ to \eqref{eqn:HBF} satisfies $V(\xp(t)) \leq V(\xp(0)) \exp \left(-(1-m)\psi t\right)$
		for all $t \ \in \dom \xp \ (= \reals_{\geq 0})$. Therefore, each maximal solution $t \mapsto \xp(t)$ to the closed-loop algorithm $\HS_0$ in \eqref{eqn:H0} satisfies \eqref{eqn:ConvergenceRateHBFNSCVX} for all $t \ \in \dom \xp \ (= \reals_{\geq 0})$. \IfConf{\hfill{} \qed}{}
	\end{proof}}
%
%
\subsection{Properties of $\HS_1$}
	\label{sec:PropertiesOfH1}
	When $L$ satisfies Assumptions \ref{ass:LisNSCVX} and \ref{ass:Lipschitz}, then we can derive an upper bound, for all $t \geq 0$, \IfConf{on $V_1$}{on the Lyapunov function $V_1$} in \eqref{eqn:LyapunovNesterovNSCVX} along solutions to $\HS_1$. To derive such a bound, we extend \cite[Proposition~3.2]{muehlebach2019dynamical}\IfConf{ to functions $L$ with generic $L^*$, $\xp^*_1$, and $\zeta > 0$}{, which assumes $L^* = 0$, $\xp_1^*=0$, and $\zeta = 1$, to the general case of $L^* \in \reals$, a single minimizer $\xp_1^* \in \reals^n$, and $\zeta > 0$}, in the following proposition. \IfConf{Its proof is in \cite{dhustigs2022unitingNSC}.}{}
	
	\IfConf{\begin{prop}}{\begin{proposition}} \label{prop:ConvergenceNSCVXNesterov} Let $L$ satisfy Assumptions \ref{ass:LisNSCVX} and \ref{ass:Lipschitz}. Then, each maximal solution $t \mapsto (\xp(t), \tau(t))$ to the closed-loop algorithm $\HS_1$ in \eqref{eqn:H1} with $\tau(0) = 0$ satisfies \IfConf{$V_1(\xp(t),t) \leq \frac{4}{(t+2)^2}V_1(\xp(0),0)$}{
		\begin{equation} \label{eqn:ConvergenceRateNCVX}
			V_1(\xp(t),t) \leq \frac{4}{(t+2)^2}V_1(\xp(0),0)
		\end{equation}}
		for all $t \geq 0$, where $V_1$ is defined via \eqref{eqn:LyapunovNesterovNSCVX}.
	\IfConf{\end{prop}}{\end{proposition}} 
\IfConf{}{\begin{proof}
	See Section \ref{sec:ProofProp54}.
	\end{proof}}
The following proposition establishes that the closed-loop algorithm $\HS_1$ has a convergence rate $\frac{1}{(t+2)^2}$ for all $t \geq 0$. To prove it, we use Proposition \ref{prop:ConvergenceNSCVXNesterov}. \IfConf{Its proof is in \cite{dhustigs2022unitingNSC}.}{} \IfConf{}{This proposition is a new result, which was not analyzed in \cite{muehlebach2019dynamical}.}
	\IfConf{\begin{prop}[Convergence rate for $\HS_1$]}{\begin{proposition}(Convergence rate for $\HS_1$)} \label{prop:UpperBoundofV}
		Let $L$ satisfy Assumptions \ref{ass:LisNSCVX} and \ref{ass:Lipschitz}. Let $\zeta > 0$ and $M > 0$ come from Assumption \ref{ass:Lipschitz}. Then, for each maximal solution $t \mapsto (\xp(t), \tau(t))$ to the closed-loop algorithm $\HS_1$ in \eqref{eqn:H1} with $\tau(0)=0$, the following holds:\IfConf{\vspace{-0.2cm}}{}
		\begin{align} \label{eqn:BoundOnL}
			&\frac{\zeta^2}{M}(L(\xp_1(t)) - L^*) \\
			& \leq V_1(\xp(t),t)
			\leq \frac{4c}{(t + 2)^2}\left(\left|\xp_1(0) - \xp_1^* \right|^2 + \left|\xp_2(0)\right|^2\right) \nonumber 
		\end{align}
		for all $t \geq 0$, where $c := \left(1 + \zeta^2\right)\exp \left(\sqrt{\frac{13}{4} + \frac{\zeta^4}{M}}\right)$. 
	\IfConf{\end{prop}}{\end{proposition}} \IfConf{\vspace{-0.4cm}}{}
	
	\IfConf{}{\begin{proof}
		The proof consists of the following steps. 
		\begin{enumerate}[label={\arabic*)},leftmargin=*]
			\item First, we use the definition of convexity in Footnote \ref{foot:Convexity} and the Lipschitz continuity of $\nabla L$ in Assumption \ref{ass:Lipschitz}, to show that $V_1$ satisfies
			\begin{equation}\label{eqn:upperBoundV1NSC}
				V_1(\xp,\tau) \leq \alpha_2 \left|\xp\right|_{\mathcal{A}_2}^2 := \left(1 + \zeta^2\right)\left|\xp\right|_{\mathcal{A}_2}^2
			\end{equation}
			where $1 + \zeta^2 > 0$;
			\item Then, we use the Lipschitz continuity of $\nabla L$ in Assumption \ref{ass:Lipschitz} and the comparison principle to show that the bound in step 1) along $t \mapsto \xp(t)$ satisfies \IfConf{$V_1(\xp(t),t) \leq \left(1 + \zeta^2\right) \exp \left(\!\!\sqrt{\frac{13}{4} + \frac{\zeta^4}{M}}t \right)\! \left|\xp(0)\right|_{\mathcal{A}_2}^2$}{$V_1(\xp(t),t) \leq \alpha_2 \exp \left(2\left(\sqrt{\frac{13}{4} + \frac{\zeta^4}{M}}\right)t \right) \left(\left|\xp_1(0) - \xp_1^* \right|^2 + \left|\xp_2(0)\right|^2\right)$} for all $t \geq 0$;
			\item Next, we show that at $t = 0$, $V_1(\xp(0),0)$ is upper bounded by\\ $c \left(\left|\xp_1(0) - \xp_1^* \right|^2 + \left|\xp_2(0)\right|^2\right)$, where \IfConf{$c = $\\$ \left(1 + \zeta^2\right)\exp \left(\sqrt{\frac{13}{4} + \frac{\zeta^4}{M}}\right)$}{$c = \left(1 + \zeta^2\right)\exp \left(\sqrt{\frac{13}{4} + \frac{\zeta^4}{M}}\right)$};
			\item Finally, we combine the bound in 3) with \eqref{eqn:ConvergenceRateNCVX} to get \eqref{eqn:BoundOnL} for all $t \geq 0$.
		\end{enumerate} 
		Proceeding with step 1), the Lyapunov function $V_1$, defined via \eqref{eqn:LyapunovNesterovNSCVX}, can be upper bounded by a class-$\mathcal{K}_{\infty}$ function, namely, defining the set 
		\begin{equation} \label{eqn:SetAForConvRate}
				\mathcal{A}_2 := \{\minSet\} \times \{0\}
		\end{equation}
		then, $V_1$ satisfies 
		\begin{equation} \label{eqn:UpperBoundOnV}
			V_1(\xp,\tau) \leq \alpha_2 \left|\xp\right|_{\mathcal{A}_2}^2
		\end{equation}
		for all $(\xp, \tau) \in \reals^{2n} \times \reals_{\geq 0}$,
		and with $\alpha_2$ derived as follows. Since $\bar{a}$, defined via \eqref{eqn:BarA}, equals $1$ at $\tau = 0$ and $\bar{a}$ is monotonically decreasing toward zero (but being always positive) as $\tau$ tends to $\infty$, then $\bar{a}$ is upper bounded by $1$, and, consequently, the first term of $V_1$ can be upper bounded, for all $(\xp, \tau) \in \reals^{2n} \times \reals_{\geq 0}$, as follows:
		\begin{equation} \label{eqn:UpperBound1}
			\frac{1}{2} \left|\frac{2}{(\tau+2)} \left(\xp_1 - \xp_1^*\right) + \xp_2\right|^2 \leq \left| \xp_1 - \xp_1^* \right|^2 + \left|\xp_2\right|^2.
		\end{equation}
		The second term of $V_1$ can be bounded as follows. Since by Assumption \ref{ass:LisNSCVX}, $L$ is $\mathcal{C}^1$, convex, and has a single minimizer $\xp_1^*$, then, since $\nabla L(\xp_1^*) = 0$, we can upper bound $L(\xp_1) - L^*$ in the following manner, using the definition of convexity in Footnote \ref{foot:Convexity} and the Lipschitz continuity of $\nabla L$ in Assumption \ref{ass:Lipschitz}, using $u_1 = \xp_1^*$ and $w_1 = \xp_1$: $\left| L(\xp_1) - L^* \right| \leq \left| \left\langle \nabla L(\xp_1), \xp_1^* - \xp_1 \right\rangle \right| \leq \left| \nabla L(\xp_1) \right| \left| \xp_1 - \xp_1^* \right| \leq M \left| \xp_1 - \xp_1^* \right|^2$,
		for all $\xp_1 \in \reals^n$. Therefore, since $L(\xp_1) \geq L^*$, we can upper bound the second term of $V_1$ as follows:
		\begin{equation} \label{eqn:UpperBound2}
			\frac{\zeta^2}{M}(L(\xp_1) - L^*) \leq \zeta^2 \left| \xp_1 - \xp_1^* \right|^2 \leq \zeta^2\left(\left| \xp_1 - \xp_1^* \right|^2 + \left|\xp_2\right|^2 \right)
		\end{equation}
		for all $\xp \in \reals^{2n}$. Using \eqref{eqn:UpperBound1} and \eqref{eqn:UpperBound2} $V_1(\xp,\tau)$ is upper bounded as in \eqref{eqn:upperBoundV1NSC} for each $\xp \in \reals^{2n}$ and each $\tau \in \reals_{\geq 0}$.
	
		Next, for step 2), in order to apply the comparison principle,
		we define the system
			\begin{equation}\label{eqn:MJ_CTDynamicsNCVX}
				\IfConf{\matt{\xp_2 \\ -2\bar{d}(t)\xp_2 - \frac{\zeta^2}{M} \nabla L(\xp_1 + \bar{\beta}(t) \xp_2)} =: f(\xp,t) \ \ \xp \in \reals^{2n}.}{
				\matt{\dot{\xp_1} \\ \dot{\xp_2}} = \matt{\xp_2 \\ -2\bar{d}(t)\xp_2 - \frac{\zeta^2}{M} \nabla L(\xp_1 + \bar{\beta}(t) \xp_2)} =: f(\xp,t) \quad \xp \in \reals^{2n}.}
		\end{equation}
		Since $\nabla L$ is Lipschitz continuous with constant $M > 0$ by Assumption \ref{ass:Lipschitz}, then using Assumption \ref{ass:Lipschitz} with $w_1 = \xp_1 + \bar{\beta}(t)\xp_2$ and $u_1 = \xp_1^*$ yields, for each $\xp_1, \xp_2 \in \reals^n$ and each $t \in \reals_{\geq 0}$, \IfConf{$\left|\nabla L(\xp_1 + \bar{\beta}(t)\xp_2)\right| \leq M \left|\xp_1 - \xp_1^* + \bar{\beta}(t)\xp_2 \right|$.}{
		\begin{equation}\label{eqn:LipschitzMJGradL}
			\left|\nabla L(\xp_1 + \bar{\beta}(t)\xp_2)\right| \leq M \left|\xp_1 - \xp_1^* + \bar{\beta}(t)\xp_2 \right|.
		\end{equation}}
		Then, since $\left|\bar{d}(t)\right| \leq \frac{3}{4}$ and $\left|\bar{\beta}(t)\right| \leq 1$ for all $t \geq 0$, we have \IfConf{$\left|f(\xp,t)\right|^2 \leq \left|\xp_2\right|^2 + \frac{9}{4} \left|\xp_2\right|^2 + \frac{\zeta^4}{M^2} \left|\nabla L(\xp_1 + \bar{\beta}(t)\xp_2)\right|^2 \leq \left(\frac{13}{4} + \frac{\zeta^4}{M} \right) \left|\xp\right|^2_{\mathcal{A}_2}$
		}{\begin{align} \label{eqn:Squaredf_TR}
			\left|f(\xp,t)\right|^2 & = \left|\xp_2\right|^2 + \left|-2\bar{d}(t)\xp_2 - \frac{\zeta^2}{M} \nabla L(\xp_1 + \bar{\beta}(t)\xp_2)\right|^2 \nonumber\\
			& \leq \left|\xp_2\right|^2 + \frac{9}{4} \left|\xp_2\right|^2 + \frac{\zeta^4}{M^2} \left|\nabla L(\xp_1 + \bar{\beta}(t)\xp_2)\right|^2 \nonumber\\
			& \leq \frac{13}{4} \left|\xp_2\right|^2 + \frac{\zeta^4}{M} \left(\left|\xp_1 - \xp_1^*\right|^2 + \left|\xp_2\right|^2\right) \nonumber\\
			& = \frac{\zeta^4}{M} \left|\xp_1 - \xp_1^*\right|^2 + \left(\frac{13}{4} + \frac{\zeta^4}{M} \right) \left|\xp_2\right|^2 \nonumber\\
			& \leq \left(\frac{13}{4} + \frac{\zeta^4}{M} \right) \left|\xp\right|^2_{\mathcal{A}_2}
		\end{align}}
		for all $\xp \in \reals^{2n}$ and all $t \in \reals_{\geq 0}$, where \IfConf{$\mathcal{A}_2 = \{\minSet\} \times \{0\}$}{$\mathcal{A}_2$ is defined via \eqref{eqn:SetAForConvRate}}. \IfConf{}{The second inequality in \IfConf{\eqref{eqn:Squaredf}}{\eqref{eqn:Squaredf_TR}} comes from applying \eqref{eqn:LipschitzMJGradL}.} 
		The comparison principle \cite[Lemma~3.4]{khalil2002nonlinear}, leads to the following bound of the norm of the solution to \eqref{eqn:MJ_CTDynamicsNCVX}: \IfConf{$\left|\xp(t)\right|_{\mathcal{A}_2} \leq \exp \left(\frac{1}{2}\sqrt{\frac{13}{4} + \frac{\zeta^4}{M}} t \right) \left|\xp(0) \right|_{\mathcal{A}_2}$,}{
		\begin{equation}
			\left|\xp(t)\right|_{\mathcal{A}_2} \leq \exp \left(\frac{1}{2}\sqrt{\frac{13}{4} + \frac{\zeta^4}{M}} t \right) \left|\xp(0) \right|_{\mathcal{A}_2}
		\end{equation}}
		for all $t \geq 0$.
		Then, \eqref{eqn:UpperBoundOnV} along $t \mapsto \left(\xp(t) \right)$ reduces to, for all $t \geq 0$\IfConf{, $V_1(\xp(t),t) \leq \left(1+\zeta^2\right) \left|\xp(t)\right|_{\mathcal{A}_2}^2 \leq \left(1+\zeta^2\right) \exp \left(\sqrt{\frac{13}{4} + \frac{\zeta^4}{M}} t \right) \left|\xp(0) \right|^2_{\mathcal{A}_2}$.}{
		\begin{align}
			V_1(\xp(t),t) & \leq \left(1+\zeta^2\right) \left|\xp(t)\right|_{\mathcal{A}_2}^2 \nonumber\\
			& \leq \left(1+\zeta^2\right) \exp \left(\sqrt{\frac{13}{4} + \frac{\zeta^4}{M}} t \right) \left( \left|\xp_1(0) - \xp_1^* \right|^2 + \left|\xp_2(0)\right|^2 \right).
		\end{align}}
		In step 3), we evaluate this bound at $t = 0$.
		Finally, for step 4), taking $c = \left(1 + \zeta^2\right)\exp \left(\!\sqrt{\frac{13}{4} + \frac{\zeta^4}{M}}\right)$, combining \eqref{eqn:ConvergenceRateNCVX} with 3) at $t = 0$ yields \eqref{eqn:BoundOnL}
		for all $t \geq 0$.\IfConf{\hfill{} \qed}{}
	\end{proof}} 

	The following proposition establishes that the closed-loop system $\HS_1$ in \eqref{eqn:H1} has the set
	\begin{equation} \label{eqn:setA_UGA}
		\mathcal{A}_1 := \{\minSet\} \times \{0\} \times \reals_{\geq 0}
	\end{equation}
	UGAS. 
	To prove it, we use Proposition \ref{prop:UpperBoundofV} and \cite[Theorem~3.18]{65}.
	\IfConf{Its proof is in \cite{dhustigs2022unitingNSC}.}{This proposition is a new result, which was not analyzed in \cite{muehlebach2019dynamical}.}

\IfConf{\begin{prop}[UGAS of $\mathcal{A}_1$ in \eqref{eqn:setA_UGA} for $\HS_1$]}{\begin{proposition}(UGAS of $\mathcal{A}_1$ in \eqref{eqn:setA_UGA} for $\HS_1$)} \label{prop:UGASNest}
	Let $L$ satisfy Assumptions \ref{ass:LisNSCVX} and \ref{ass:Lipschitz}. Let $\zeta > 0$ and let $M > 0$ come from Assumption \ref{ass:Lipschitz}. Then, the set $\mathcal{A}_1$ in \eqref{eqn:setA_UGA} is UGAS for $\HS_1$.
\IfConf{\end{prop}}{\end{proposition}} 

\IfConf{}{\begin{proof}
	By Proposition \ref{prop:Existence}, each maximal solution to $\HS_1$ in \eqref{eqn:H1} is complete and unique. Next, since $L$ is $\mathcal{C}^1$, convex, and has a unique minimizer by Assumption \ref{ass:LisNSCVX}, then $\mathcal{A}_1 \subset \reals^{2n} \times \reals_{\geq 0}$, defined via \eqref{eqn:setA_UGA}, is closed by construction, satisfying the first assumption of \cite[Theorem~3.18]{65}. Then, since by Assumption \ref{ass:LisNSCVX}, $L$ is $\mathcal{C}^1$, then $V_1$ in \eqref{eqn:LyapunovNesterovNSCVX} is continuously differentiable and therefore, since $\reals^{2n} \times \reals_{\geq 0} \subset \dom V_1$, $V_1$ is a Lyapunov function candidate for $\HS_1$ by \cite[Definition~3.16]{65}, satisfying the second assumption of \cite[Theorem~3.18]{65}. 
	
	Next, since the distance of the state $\tau$ to $\reals_{\geq 0}$ is always zero, then 
	we show that $V_1$ in \eqref{eqn:LyapunovNesterovNSCVX} is radially unbounded in $\xp$, relative to $\mathcal{A}_2$, defined via \eqref{eqn:SetAForConvRate}. 
	Since $L$ has quadratic growth away from $\xp^*_1$ with constant $\alpha > 0$, by Assumption \ref{ass:QuadraticGrowth}, and due to $\bar{a}$, defined via \eqref{eqn:BarA}, equaling $1$ at $\tau = 0$ and monotonically decreasing toward zero (but being always positive) as $\tau$ tends to $\infty$, then we lower bound $V_1$ as follows: \IfConf{ 
		
		\vspace{-0.85cm}
		
	\begin{align}\label{eqn:Alpha2Func}
		V_1(\xp, \tau) = & \frac{1}{2}\left|\bar{a}(\tau)\left(\xp_1 - \xp_1^*\right) + \xp_2\right|^2 + \frac{\zeta^2}{M} (L(\xp_1) - L^*) \nonumber\\
		\geq & \frac{1}{2}\left|\bar{a}(\tau)\left(\xp_1 - \xp_1^*\right) + \xp_2\right|^2 +
		\frac{\alpha \zeta^2}{M} \left|\xp_1 - \xp_1^*\right|^2 \nonumber\\
		\geq & \left(\frac{\bar{a}^2(\tau)}{2} + \frac{\alpha \zeta^2}{M}\right)\left|\xp_1 - \xp_1^*\right|^2 \nonumber\\
		& + \bar{a}(\tau)\left\langle \xp_1 - \xp^*_1,\xp_2 \right\rangle + \frac{1}{2} \left|\xp_2\right|^2 \nonumber\\
		\geq & \matt{\left(\xp_1 - \xp^*_1\right)^{\top} & \xp^{\top}_2} B \matt{\xp_1 - \xp^*_1\\ \xp_2} 
	\end{align}}{
	\begin{align}\label{eqn:Alpha2FuncTR}
		V_1(\xp, \tau) = & \frac{1}{2}\left|\bar{a}(\tau)\left(\xp_1 - \xp_1^*\right) + \xp_2\right|^2 + \frac{\zeta^2}{M} (L(\xp_1) - L^*) \\
		\geq & \frac{1}{2}\left|\bar{a}(\tau)\left(\xp_1 - \xp_1^*\right) + \xp_2\right|^2 +
		\frac{\alpha \zeta^2}{M} \left|\xp_1 - \xp_1^*\right|^2 \nonumber\\
		\geq & \frac{\bar{a}^2(\tau)}{2} \left|\xp_1 - \xp^*_1\right|^2 + \bar{a}(\tau)\left\langle \xp_1 - \xp^*_1,\xp_2 \right\rangle + \frac{1}{2} \left|\xp_2\right|^2 +
		\frac{\alpha \zeta^2}{M} \left|\xp_1 - \xp_1^*\right|^2\nonumber\\
		\geq & \left(\frac{\bar{a}^2(\tau)}{2} + \frac{\alpha \zeta^2}{M}\right)\left|\xp_1 - \xp_1^*\right|^2 + \frac{\bar{a}(\tau)}{2}\left\langle \xp_1 - \xp^*_1,\xp_2 \right\rangle \nonumber\\
		& + \frac{\bar{a}(\tau)}{2}\left\langle \xp_1 - \xp^*_1,\xp_2 \right\rangle + \frac{1}{2} \left|\xp_2\right|^2 \nonumber\\
		\geq & \matt{\left(\xp_1 - \xp^*_1\right)^{\top} & \xp^{\top}_2} P \matt{\xp_1 - \xp^*_1\\ \xp_2} \nonumber
	\end{align}}
	for each $\xp \in \reals^{2n}$ and each $\tau \in \reals_{\geq 0}$, where
	\begin{equation} \label{eqn:MatrixB}
		P := \matt{\left(\frac{\bar{a}^2(\tau)}{2} + \frac{\alpha \zeta^2}{M}\right) & \frac{\bar{a}(\tau)}{2}\\ \frac{\bar{a}(\tau)}{2} & \frac{1}{2}}.
	\end{equation}
	Next, we show that $P$ in \eqref{eqn:MatrixB} is positive definite, so that there exists $\alpha_1$ such that\footnote{It was already shown that there exists $\alpha_2$ such that the upper bound on $V_1$ in \eqref{eqn:upperBoundV1NSC} holds.}
	\begin{equation} \label{eqn:RadiallyUnboundedV1NSC}
		\alpha_1\left|\xp\right|^2_{\mathcal{A}_2} \leq \matt{\left(\xp_1 - \xp^*_1\right)^{\top} & \xp^{\top}_2} P \matt{\xp_1 - \xp^*_1\\ \xp_2}  \leq V_1(\xp,\tau)
	\end{equation}
	for each $\xp \in \reals^{2n}$ and each $\tau \in \reals_{\geq 0}$. To that end, we show that the leading principal minors of $P$ in \eqref{eqn:MatrixB} are strictly positive, as follows. Since $\bar{a}(\tau) \in (0,1]$ for each $\tau \in \reals_{\geq 0}$, $\alpha > 0$ from Assumption \ref{ass:QuadraticGrowth}, $\zeta > 0$, and $M > 0$ from Assumption \ref{ass:Lipschitz}, we have: 
	\begin{subequations}
		\begin{align}
			& \left(\frac{\bar{a}^2(\tau)}{2} + \frac{\alpha \zeta^2}{M}\right) > 0 \\
			&\det(P) = \frac{1}{2}\left(\frac{\bar{a}^2(\tau)}{2} + \frac{\alpha \zeta^2}{M}\right) - \left(\frac{\bar{a}(\tau)}{2}\right)^2 = \frac{\alpha \zeta^2}{2M} > 0.
		\end{align}
	\end{subequations} 
	Therefore, since the leading principal minors of $P$ are strictly positive, then $P$ is positive definite. Hence, there exists $\alpha_1$ such that \eqref{eqn:RadiallyUnboundedV1NSC} is true, and $V_1$ is radially unbounded in $\xp$, relative to $\mathcal{A}_2$. 
	
	By Proposition \ref{prop:ConvergenceNSCVXNesterov}, $V_1$ satisfies \IfConf{$\dot{V}_1(\xp,\tau) \leq -\bar{a}(\tau)V_1(\xp,\tau)$}{\eqref{eqn:MJ_dotVNCVX}} for each $\xp \in \reals^{2n}$ and $\tau \in \reals_{\geq 0}$. Since $L$ is $\mathcal{C}^1$, convex, and has a unique minimizer by Assumption \ref{ass:LisNSCVX}, then $L$ is positive definite with respect to $\minSet$ and, consequently, $V_1$ is positive definite with respect to $\mathcal{A}_1$ in \eqref{eqn:setA_UGA}. Then, since $\bar{a}(\tau) \in (0,1]$ for each $\tau \geq 0$, $\rho\left(\left|x\right|_{\mathcal{A}_1}\right) := \bar{a}(\tau)V_1(\xp,\tau)$ is positive definite with respect to $\mathcal{A}_1$. Therefore, by an application of \cite[Theorem~3.18]{65}, every complete solution to \eqref{eqn:H1} converges to $\mathcal{A}_1$ in \eqref{eqn:setA_UGA}. The arguments above involving the Lyapunov theorem in \cite[Theorem~3.18]{65} yield UGpAS of $\mathcal{A}_1$ for $\HS_1$ in  \eqref{eqn:H1}. Since by Proposition \ref{prop:Existence}, each maximal solution to $\HS_1$ is complete, then $\mathcal{A}_1$ is UGAS for $\HS_1$. 
\end{proof}} 
%
%
\subsection{Uniform Global Asymptotic Stability of $\mathcal{A}$ for $\HS$}
	\label{sec:UGAS}
	
	The hybrid closed-loop algorithm $\HS$ satisfies the hybrid basic conditions by Lemma \ref{lemma:HBC}, satisfying the first assumption of \IfConf{\cite[Theorem~A.3]{dhustigs2022unitingNSC}}{Theorem \ref{thm:hybrid Lyapunov theorem}}. 
	Furthermore, $\Pi(C_0) \cup \Pi(D_0) = \reals^{2n}$, $\Pi(C_1) \cup \Pi(D_1) = \reals^{2n}$, and each maximal solution $(t,j) \mapsto x(t,j) = (\xp(t,j), \xlogic(t,j),\tau(t,j))$ to $\HS$ in \eqref{eqn:HS-TimeVarying}-\eqref{eqn:CAndDGradientsNestNSC} is complete and bounded by \IfConf{\cite[Proposition~2.9]{dhustigs2022unitingNSC}}{Proposition \ref{prop:Existence}}. Since by Assumption \ref{ass:LisNSCVX}, $L$ has a unique minimizer $\minSet$, then $\mathcal{A}$, defined via \eqref{eqn:SetOfMinimizersHS-NSCVX}, is compact by construction, and $\mathcal{U} = \reals^{2n} \times \xlogicSpace \times \reals_{\geq 0}$ contains a nonzero open neighborhood of $\mathcal{A}$, satisfying the second assumption of \IfConf{\cite[Theorem~A.3]{dhustigs2022unitingNSC}}{Theorem \ref{thm:hybrid Lyapunov theorem}}. 
	
	To prove \IfConf{attractivity\footnote{For a definition of attractivity, see \cite[Definition~3.1]{220}.}}{attractivity} of $\mathcal{A}$, we proceed by contradiction. Suppose there exists a complete solution $x$ to $\HS$ such that $\lim\limits_{t+j \rightarrow \infty}\left|x(t,j)\right|_{\mathcal{A}} \neq 0$. 
	Since \IfConf{\cite[Proposition~2.9]{dhustigs2022unitingNSC}}{Proposition \ref{prop:Existence}} guarantees completeness of maximal solutions, 
	we have the following cases:
	\begin{enumerate}[label={\alph*)},leftmargin=*]
		\item There exists $(t',j') \in \dom x$ such that $x(t,j) \in C_1 \setminus D_1$ for all $(t,j) \in \dom x, t+j \geq t' + j'$;
		\item There exists $(t',j') \in \dom x$ such that $x(t,j) \in C_0 \setminus (\mathcal{A} \cup D_0)$ for all $(t,j) \in \dom x, t+j \geq t' + j'$;
		\item There exists $(t',j') \in \dom x$ such that $x(t,j) \in D$ for all $(t,j) \in \dom x, t+j \geq t' + j'$.
	\end{enumerate}\IfConf{\vspace{-0.1cm}}{}

	Case a) contradicts the fact that, by Proposition \ref{prop:UGASNest}, the set $\mathcal{A}_1$, defined via \eqref{eqn:setA_UGA}, is UGAS for $\HS_1$. Such UGAS of $\mathcal{A}_1$\IfConf{}{, guaranteed by Proposition \ref{prop:UGASNest},} implies there exist $\tilde{c}_1 \in (0,\tilde{c}_{1,0})$ and $d_1 \in (0,d_{1,0})$ such that the state $\xp$ reaches $\left(\{\minSet\} + \tilde{c}_1 \ball \right) \times \left(\{0\} + d_1 \ball\right) \subset \T_{1,0}$ at some finite flow time $t \geq 0$ or as $t \rightarrow \infty$. In turn, due to the construction of $C_1$ and $D_1$ in \eqref{eqn:CAndDGradientsNestNSC}, with $\T_{1,0}$ defined via \eqref{eqn:T10}, the solution $x$ must reach $D_1$ at some $(t,j) \in \dom x, t+j \geq t' + j'$. Therefore, case a) does not happen. 
	
	Case b) contradicts the fact that, by Proposition \ref{prop:GAS-HBF}, $\{\minSet\} \times \{0\}$ is UGAS for $\HS_0$. In fact, $\lim\limits_{t+j \rightarrow \infty} \left|x(t,j)\right|_{\mathcal{A}} = 0$, and since $\mathcal{A} \subset C_0$, case b) does not happen.
	
	Case c) contradicts the fact that, due to the construction of $\T_{1,0}$ in \eqref{eqn:T10} and $\T_{0,1}$ in \eqref{eqn:T01}, we have $G(D) \cap D := \left(\left(\T_{0,1} \times \{1\} \times \{0\} \right) \cup \left(\T_{1,0} \times \{0\} \times \{0\} \right)\right)$\\$\cap \left(\left(\T_{0,1} \times \{0\} \times \{0\} \right) \cup \left(\T_{1,0}\times \{1\} \times \reals_{\geq 0} \right) \right) = \emptyset$
	where \IfConf{}{$G(D)$ is defined via \eqref{eqn:G(D)} and} $D$ is defined in \eqref{eqn:CAndDGradientsNestNSC}. Such an equality holds since $\T_{1,0} \cap \T_{0,1} = \emptyset$; see the end of Section \ref{sec:DesignT10}. Therefore, case c) does not happen.

	Therefore, cases a)-c) do not happen, and each maximal and complete solution $x = (\xp, \xlogic, \tau)$ to $\HS$ with $\tau(0,0) = 0$ converges to $\mathcal{A}$. Consequently, by the construction of $C$ and $D$ in \eqref{eqn:CAndDGradientsNestNSC}, the UGAS of $\mathcal{A}_1$ (defined via \eqref{eqn:setA_UGA}) for $\HS_1$ established in Proposition \ref{prop:UGASNest}, the UGAS of $\{\minSet\} \times \{0\}$ for $\HS_0$ established in Proposition \ref{prop:GAS-HBF}, and since each maximal solution to $\HS$ is complete by \IfConf{\cite[Proposition~2.9]{dhustigs2022unitingNSC}}{Proposition \ref{prop:Existence}}, the set $\mathcal{A}$ is UGAS for $\HS$. 

	To show that each maximal and complete solution $x$ to $\HS$ jumps no more than twice, we proceed by contradiction. Without loss of generality, suppose there exists a maximal and complete solution that jumps three times. We have the following possible cases: \IfConf{i) The solution first jumps at a point in $D_0$, then jumps at a point in $D_1$, and then jumps at a point in $D_0$; or ii) The solution first jumps at a point in $D_1$, then jumps at a point in $D_0$, and then jumps at a point in $D_1$.}{
	\begin{enumerate}[label={\roman*)},leftmargin=*]
		\item The solution first jumps at a point in $D_0$, then jumps at a point in $D_1$, and then jumps at a point in $D_0$; or
		\item The solution first jumps at a point in $D_1$, then jumps at a point in $D_0$, and then jumps at a point in $D_1$.
	\end{enumerate}}
	Case i) does not hold since, once the jump in $D_1$ occurs, the solution $x$ is in $(\T_{1,0} \times \{0\} \times \{0\}) \subset C_0$. Due to the construction of $\T_{1,0}$ in \eqref{eqn:T10} and $\T_{0,1}$ in \eqref{eqn:T01} such that $\T_{1,0} \cap \T_{0,1} = \emptyset$, as described in the contradiction of case c) above, and due to the UGAS of $\minSet \times \{0\}$ for $\HS_0$ by Proposition \ref{prop:GAS-HBF}, the solution $x$ will never return to $D_0$. Therefore, case i) does not happen. Case ii) leads to a contradiction for the same reason, and in this case, once the first jump in $D_1$ occurs, no more jumps happen. Therefore, since cases i)-ii) do not happen, each maximal and complete solution $x$ to $\HS$ with $\tau(0,0)=0$ has no more than two jumps.
	
	\subsection{Convergence Rate of $\HS$}
	\label{sec:ConvRateH}
	Finally, we prove the hybrid convergence rate of $\HS$. Letting $\zeta > 0$ and letting $M > 0$ come from Assumption \ref{ass:Lipschitz}, then by Proposition \ref{prop:UpperBoundofV}, since $L$ satisfies Assumptions \ref{ass:LisNSCVX} and \ref{ass:Lipschitz}, each maximal solution $t \mapsto (\xp(t), \tau(t))$ to the closed-loop algorithm $\HS_1$ with $\tau(0,0) = 0$ satisfies \eqref{eqn:BoundOnL}, for all $t \geq 0$, where $c$ is defined below \eqref{eqn:BoundOnL}. By Proposition \ref{prop:HBFConvergenceRate}, since $L$ satisfies Assumptions \ref{ass:LisNSCVX} and \ref{ass:QuadraticGrowth}, then, given $\gamma > 0$ and $\lambda > 0$, for each $m \in (0,1)$ such that $\psi := \frac{m\alpha\gamma}{\lambda} > 0$ and $\nu := \psi (\psi - \lambda) < 0$,
 	each maximal solution $t \mapsto \xp(t)$ to the closed-loop algorithm $\HS_0$ satisfies \eqref{eqn:ConvergenceRateHBFNSCVX} for all $t \ \in \dom \xp \ (= \reals_{\geq 0})$. 	
	Since maximal solutions $(t,j) \mapsto x(t,j) = (\xp(t,j),\xlogic(t,j),\tau(t,j))$ to $\HS$ starting from $C_1$ are guaranteed to jump no more than once, as implied by the contradiction in cases i)-ii) above, then the domain of each maximal solution $x$ to $\HS$ starting from $C_1$ is $\cup_{j=0}^{1} (I^j,j)$, with $I^0$ of the form $[t_0,t_1]$ and with $I^1$ of the form $[t_1,\infty)$. Therefore, given $\zeta > 0$, $\lambda > 0$, $\gamma > 0$, $c_{1,0} \in (0,c_0)$, $\varepsilon_{1,0} \in (0,\varepsilon_0)$, $\alpha > 0$ from Assumption \ref{ass:QuadraticGrowth}, and $M >0$ from Assumption \ref{ass:Lipschitz}, due to the construction of $\mathcal{U}_0$, $\T_{1,0}$, and $\T_{0,1}$ in \eqref{eqn:U0}, \eqref{eqn:T10}, and \eqref{eqn:T01}, with $\tilde{c}_{1,0} \in (0,\tilde{c}_0)$ and $d_{1,0} \in (0,d_0)$ defined via \eqref{eqn:UTilde0SetEquations} and \eqref{eqn:TTilde10SetEquations}, and due to the individual convergence rates of $\HS_1$ and $\HS_0$, each maximal solution $(t,j) \mapsto x(t,j) = (\xp(t,j),\xlogic(t,j),\tau(t,j))$ to the hybrid closed-loop algorithm $\HS$ that starts in $C_1$, such that $\tau(0,0)=0$, satisfies \eqref{eqn:UnitingConvergenceRateNSCHS1} for each $t \in I^0$ at which $q(t,0)$ is equal to 1 and $t \geq 0$, and satisfies item \ref{item:3} of Theorem~\ref{thm:GASNestNSC} for each $t \in I^1$ at which $q(t,1)$ is equal to 0.
\IfConf{
\begin{rem}
	Although it is outside the scope of this paper, a potential approach to discretizing the hybrid closed-loop algorithm $\HS$ in \eqref{eqn:HS-TimeVarying}-\eqref{eqn:CAndDGradientsNestNSC} can be found in \cite{40}. Such a discretization approach, which is designed for hybrid systems and which has assumptions that are satisfied by forward Euler and p-stage Runge-Kutta consistent methods, for example, would yield results similar to Theorem \ref{thm:GASNestNSC}, Proposition \ref{prop:GAS-HBF}, Lemma \ref{lemma:IntermediateDotV}, and Propositions \ref{prop:HBFConvergenceRate}, \ref{prop:ConvergenceNSCVXNesterov}, \ref{prop:UpperBoundofV}, and \ref{prop:UGASNest}.
\end{rem}}{}

\IfConf{}{	
\section{Extensions}
Some possible extensions of Theorem \ref{thm:GASNestNSC}, Proposition \ref{prop:GAS-HBF}, Lemma \ref{lemma:IntermediateDotV}, and Propositions \ref{prop:HBFConvergenceRate}, \ref{prop:ConvergenceNSCVXNesterov}, \ref{prop:UpperBoundofV}, and \ref{prop:UGASNest} are as follows. A potential approach to discretizing the hybrid closed-loop algorithm $\HS$ in \eqref{eqn:HS-TimeVarying}-\eqref{eqn:CAndDGradientsNestNSC} can be found in \cite{40}. Such a discretization approach, which is designed for hybrid systems and which has assumptions that are satisfied by forward Euler and p-stage Runge-Kutta consistent methods, for example, would yield results similar to Theorem \ref{thm:GASNestNSC}, Proposition \ref{prop:GAS-HBF}, Lemma \ref{lemma:IntermediateDotV}, and Propositions \ref{prop:HBFConvergenceRate}, \ref{prop:ConvergenceNSCVXNesterov}, \ref{prop:UpperBoundofV}, and \ref{prop:UGASNest}. 

It is possible to extend Theorem \ref{thm:GASNestNSC}, Proposition \ref{prop:GAS-HBF}, Lemma \ref{lemma:IntermediateDotV}, and Propositions \ref{prop:HBFConvergenceRate}, \ref{prop:ConvergenceNSCVXNesterov}, \ref{prop:UpperBoundofV}, and \ref{prop:UGASNest} to include $\mathcal{C}^1$, convex objective functions $L$ with a compact and connected set of minimizers. Such an extension could be achieved via
the use of Clarke's generalized derivative (see \cite{clarke1990optimization}). Additionally, Clarke's generalized derivative could be utilized to extend the analysis of the hybrid closed-loop algorithm to include nonsmooth convex objective functions $L$ with a compact and connected set of minimizers.}
\section{Conclusion}
\IfConf{

\vspace{-0.3cm}

}{}
We presented an algorithm, designed using hybrid system tools, that unites Nesterov's accelerated algorithm and the heavy ball algorithm to ensure fast convergence and UGAS of the unique minimizer for $\mathcal{C}^1$, convex objective functions $L$. The hybrid convergence rate is $\frac{1}{(t+2)^2}$ globally and exponential locally. In simulation, we showed performance improvement not only over the individual heavy ball and Nesterov algorithms, but also over the HAND-1 algorithm in \cite{poveda2019inducing}. In the process, we proved the existence of solutions for the individual heavy ball and Nesterov algorithms, and we extended the convergence rate results for Nesterov's algorithm in \cite{muehlebach2019dynamical} to functions $L$ with generic $\xp_1^*$, $L^*$, and $\zeta > 0$. Additionally, we established UGAS of the minimizer for Nesterov's algorithm, when $L$ is $\mathcal{C}^1$, convex, and has a unique minimizer. Future work will extend the uniting algorithm to a general framework, allowing the local and global algorithms to be any accelerated gradient algorithm. We will also extend the uniting algorithm to learning applications. 
	
\IfConf{}{
\appendix
\appendixpage
\renewcommand{\thesection}{\Alph{section}}
\section{Proof of Lemma \ref{lemma:HBC}}
\label{sec:ProofHBCLemma}

The objective function $L$ is $\mathcal{C}^1$, convex, and has a single minimizer by Assumption \ref{ass:LisNSCVX}. 
Therefore, since $\nabla L$ is continuous, the following hold: the set $\mathcal{U}_0$, defined via \eqref{eqn:U0}, is closed since it is a sublevel set of the continuous function $V_0$; due to $\bar{a}$ in \eqref{eqn:BarA} being continuous, the set $\T_{1,0}$, defined via \eqref{eqn:T10}, is closed since it is a sublevel set of the continuous function $V_1$; the set $\T_{0,1}$, defined via \eqref{eqn:T01}, is closed since it is the closed complement of a set. Therefore, since the sets $\mathcal{U}_0$, $\T_{1,0}$, and $\T_{0,1}$ are closed, then the sets $D_0$, $D_1$, $C_0$, and $C_1$ are closed. Since $C$ and $D$ are both finite unions of finite and closed sets, then $C$ and $D$ are also closed. 

Since $\bar{d}$ and $\bar{\beta}$, defined via \eqref{eqn:dBarBetaBar}, are continuous, and since by Assumption \ref{ass:LisNSCVX}, $L$ is $\mathcal{C}^1$, then $\hp_{\xlogic}$ in \eqref{eqn:H0H1NSCNesterovHBF} and $\hcq$ in \eqref{eqn:StaticStateFeedbackLawsNSC} are continuous. In turn, the map $\xp \mapsto$ $F_P(\xp,\kappa_{\xlogic}(\hp_{\xlogic}(\xp,\tau),\tau))$ is also continuous since $F_P$ in \eqref{eqn:HBFplant-dynamicsTR} is a $\mathcal{C}^1$ function of $\hcq$ and $\hp_{\xlogic}$. Therefore, $x \mapsto F(x)$ is continuous. The map $\g$ satisfies \ref{item:A3} by construction since it is continuous.
\section{Proof of Proposition \ref{prop:Existence}}
\label{sec:ProofExistence}
Since Assumptions \ref{ass:LisNSCVX}, \ref{ass:QuadraticGrowth}, and \ref{ass:Lipschitz} hold, then $\HS$ satisfies the hybrid basic conditions by Lemma \ref{lemma:HBC}. With $\tilde{c}_0 > 0$ and $d_0 > 0$ defined via \eqref{eqn:UTilde0SetEquations}, since $L$ is $\mathcal{C}^1$, convex, has a single minimizer by Assumption \ref{ass:LisNSCVX}, and has quadratic growth away from $\minSet$ by Assumption \ref{ass:QuadraticGrowth}, from the arguments below \eqref{eqn:c0SublevelSet}, every $\xp \in \mathcal{U}_0$ belongs to the $c_0$-sublevel set of $V_0$; recall that $\mathcal{U}_0$ is defined in in \eqref{eqn:U0} and that $V_0$ is defined via \eqref{eqn:LyapunovHBF}. Additionally, since by Assumption \ref{ass:QuadraticGrowth} $L$ has quadratic growth away from $\minSet$, then $\T_{0,1}$ in \eqref{eqn:T01}, defines the closed complement of a sublevel set of $V_0$ with level equal to $c_0$. Therefore, due to the definitions of $\mathcal{U}_0$ in \eqref{eqn:U0} and $\T_{0,1}$ in \eqref{eqn:T01}, $\Pi(C_0) \cup \Pi(D_0) = \reals^{2n}$. Furthermore, since $\T_{1,0}$ is defined via \eqref{eqn:T10}, and since by the definitions of $C_1$ and $D_1$ in \eqref{eqn:CAndDGradientsNestNSC}, $C_1$ is the closed complement of $D_1$, then $\Pi(C_1) \cup \Pi(D_1) = \reals^{2n}$.

Due to the definitions of $C_0$, $D_0$, $C_1$, and $D_1$ in \eqref{eqn:CAndDGradientsNestNSC}, $\mathcal{U}_0$ in \eqref{eqn:U0}, $\T_{1,0}$ in \eqref{eqn:T10}, and $\T_{0,1}$ in \eqref{eqn:T01}, then $C \setminus D$ is equal to $\mathrm{int}(C)$. Hence, for each point $x \in C \setminus D$, the tangent cone to $C$ at $x$ is 
\begin{equation} \label{eqn:TangentCone}
	T_C(x) := \begin{cases}  \reals^{2n} \times \{0\} \times \{0\} & \ \ \text{if } x \in C_0 \setminus D_0,\\
		\reals^{2n} \times \{1\} \times \reals_{\geq 0} & \ \ \text{if } x \in C_1 \setminus D_1.
	\end{cases}
\end{equation}
Therefore, $F(x) \cap T_C(x) \neq \emptyset$, satisfying (VC) of for each point $x \in C \setminus D$, and nontrivial solutions exist for every initial point in $\left(C_0 \cup C_1\right) \cup \left(D_0 \cup D_1\right)$, where $\Pi(C_0) \cup \Pi(D_0) = \reals^{2n}$ and $\Pi(C_1) \cup \Pi(D_1) = \reals^{2n}$. To prove that item (c) of Proposition \ref{prop:SolnExistence} does not hold, we need to show that $G(D) \subset C \cup D$. With $D$ defined in \eqref{eqn:CAndDGradientsNestNSC}, $G(D)  = \left(\T_{0,1} \times \{1\} \times \{0\} \right) \cup \left(\T_{1,0} \times \{0\} \times \{0\} \right)$.
\begin{equation} \label{eqn:G(D)}
		G(D)  = \left(\T_{0,1} \times \{1\} \times \{0\} \right) \cup \left(\T_{1,0} \times \{0\} \times \{0\} \right).
\end{equation}
Notice that $\T_{1,0} \times \{0\} \times \{0\} \subset C_0$ and $\T_{0,1} \times \{1\} \times \{0\} \subset C_1$. Therefore, $G(D) \subset C$; hence $G(D) \subset C \cup D$. Therefore, item (c) of Proposition \ref{prop:SolnExistence} does not hold. Then it remains to prove that item (b) does not happen. 

To this end, we show first that $\HS_0$, defined via \eqref{eqn:H0}, has no finite time escape\footnote{{\em Finite escape time} describes when there exists a solution $t \mapsto x(t)$ to a continuous-time nonlinear system that satisfies $\lim\limits_{t \nearrow t_e} \left|x(t)\right| = \infty$ for some finite time $t_e$.}, and has unique and bounded solutions. Since $L$ is $\mathcal{C}^1$ by Assumption \ref{ass:LisNSCVX}, and $\nabla L$ is Lipschitz continuous by Assumption \ref{ass:Lipschitz}, then $\hp_0$ in \eqref{eqn:H0H1NSCNesterovHBF} and $\hczero$ in \eqref{eqn:StaticStateFeedbackLawLocal} are Lipschitz continuous, which, since $F_P$ is a $\mathcal{C}^1$ function of $\hp_0$ and $\hczero$, means the map $\xp \mapsto F_P(\xp,\hczero(\hp_0(\xp)))$ is also Lipschitz continuous. Therefore, by \cite[Theorem~3.2]{khalil2002nonlinear}, $\dot{\xp} = F_P(\xp,\hczero(\hp_0(\xp)))$ has no finite time escape and each maximal solution to $\HS_0$ is unique. To show that each maximal solution to $\HS_0$ is bounded, we use the Lyapunov function in \eqref{eqn:LyapunovHBF}, defined for each $\xp \in \reals^{2n}$. Then, solutions to $\dot{\xp} = F_P(\xp,\hczero(\hp_0(\xp)))$ starting from any $c_V$-sublevel set $W := \defset{\xp \in \reals^{2n}}{V_0(\xp) \leq c_V}$, $c_V \geq 0$, remains in such a set for all time since $V_0$ satisfies 
\begin{equation}\label{eqn:VdotHBF}
	\dot{V}_0(\xp)
	\!=\!
	\langle \nabla V_0(\xp),\fp(\xp,\hczero(\hp_0(\xp))) \rangle 
	\!=\! 
	-\lambda \left|\xp_2\right|^2 \leq 0
\end{equation}
for each $\xp \in \reals^{2n}$, since $\lambda$ is positive. Then, to show that $V_0$ is radially unbounded, we derive class-$\mathcal{K}_{\infty}$ functions $\alpha_1$ and $\alpha_2$ such that\footnote{Since $L$ has quadratic growth away from $\xp_1^*$ by Assumption \ref{ass:QuadraticGrowth}, then the choice of $\alpha_1$ comes from lower bounding $L(\xp_1) - L^*$ in $V_0$ via Assumption \ref{ass:QuadraticGrowth}. The choice of $\alpha_2$ comes from the following: since $L$ is $\mathcal{C}^1$, convex, and has a single minimizer by Assumption \ref{ass:LisNSCVX}, then the expression $L(\xp_1) - L^*$ in $V_0$ is upper bounded using the definition of convexity in Footnote \ref{foot:Convexity}, by the same process that $L(\xp_1) - L^*$ is upper bounded in \eqref{eqn:NearOptimalityC}, to get \eqref{eqn:c0SublevelSet}. Then, $\left|\nabla L(\xp_1)\right|$ in \eqref{eqn:c0SublevelSet} is upper bounded via Assumption \ref{ass:Lipschitz} with $u_1 = \xp_1^*$ and $w_1 = \xp_1$.}, for all $\xp \in \reals^{2n}$, with $\xp^* := (\xp_1^*,0)$,
\begin{align} \label{eqn:RadiallyUnboundedV0TR}	
		\alpha_1(\left|\xp - \xp^*\right|) := \min \left\{\!\alpha \gamma,\frac{1}{2}\! \right\} \left|\xp - \xp^*\right|^2 & \leq V_0(\xp) \nonumber\\
		&\leq \alpha_2(\left|\xp - \xp^*\right|) := \left(\!M \gamma + \frac{1}{2}\!\right) \left|\xp - \xp^*\right|^2.
\end{align}
Since $L$ has quadratic growth away from $\xp_1^*$ with constant $\alpha > 0$ by Assumption \ref{ass:QuadraticGrowth}, then the choice of $\alpha_1$ comes from lower bounding $V_0$ as follows
	\begin{align}
		V_0(\xp) = \gamma \left(L(\xp_1) - L^*\right) + \frac{1}{2} \left|\xp_2\right|^2 & \geq \alpha \gamma \left|\xp_1 - \xp_1^*\right|^2 + \frac{1}{2} \left|\xp_2\right|^2 \nonumber\\
		& \geq \min \left\{\!\alpha \gamma,\frac{1}{2}\! \right\} \left|\xp - \xp^*\right|^2 = \alpha_1(\left|\xp - \xp^*\right|)
	\end{align}
	for all $\xp \in \reals^{2n}$.
	The choice of $\alpha_2$ comes from the following. Since $L$ is $\mathcal{C}^1$, convex, has a single minimizer by Assumption \ref{ass:LisNSCVX}, and since $\nabla L$ is Lipschitz continuous with constant $M > 0$ by Assumption \ref{ass:Lipschitz}, we upper bound $V_0$ in the following manner, using the definition of convexity in Footnote \ref{foot:Convexity} to get \eqref{eqn:c0SublevelSet}, and then using the Lipschitz bound in Assumption \ref{ass:Lipschitz} with $u_1 = \xp_1^*$ and $w_1 = \xp_1$ to upper bound \eqref{eqn:c0SublevelSet}, yielding
	\begin{align} \label{eqn:UpperBoundV0}
		V_0(\xp) = \gamma \left(L(\xp_1) - L^*\right) + \frac{1}{2}\left|\xp_2\right|^2 & \leq \gamma \left|\nabla L(\xp_1)\right| \left| \xp_1 - \xp_1^* \right| + \frac{1}{2}\left|\xp_2\right|^2 \nonumber\\
		& \leq M \gamma \left|\xp_1 - \xp_1^*\right|^2 + \frac{1}{2}\left|\xp_2\right|^2 \nonumber\\
		& \leq \left(\!M \gamma + \frac{1}{2}\!\right) \left|\xp - \xp^*\right|^2 = \alpha_2(\left|\xp_1 - \xp_1^*\right|)
	\end{align}
	for all $\xp \in \reals^{2n}$. Since \eqref{eqn:RadiallyUnboundedV0TR} is satisfied for $V_0$ in \eqref{eqn:LyapunovHBF} for all $\xp \in \reals^{2n}$, then $V_0$ is radially unbounded (in $\xp$, relative to $\{\minSet\} \times \{0\}$). Therefore, $W$ is compact and, due to \eqref{eqn:VdotHBF}, forward invariant for $\HS_1$, that is, any nontrivial solution starting in the subset $W$ is complete and stays in $W$. Therefore, each maximal solution to $\HS_0$, defined via \eqref{eqn:H0}, is bounded. 

Next, we show that $\HS_1$ in \eqref{eqn:H1} has no finite time escape from $\reals^{2n} \times \reals_{\geq 0}$, and has unique solutions. Since $\bar{d}$ and $\bar{\beta}$, defined via \eqref{eqn:dBarBetaBar}, are continuous, and since by Assumption \ref{ass:LisNSCVX}, $L$ is $\mathcal{C}^1$, then $\hp_1$ in \eqref{eqn:H0H1NSCNesterovHBF} and $\hcone$ in \eqref{eqn:StaticStateFeedbackLawsNSC} are also continuous. Furthermore, since by Assumption \ref{ass:Lipschitz} $\nabla L$ is Lipschitz continuous, then $\hp_1$ in \eqref{eqn:H0H1NSCNesterovHBF} and $\hcone$ in \eqref{eqn:StaticStateFeedbackLawsNSC} are Lipschitz continuous which, in turn, means the map $\xp \mapsto F_P(\xp,\hcone(\hp_1(\xp,\tau),\tau))$ is Lipschitz continuous. Consequently, since the map $\xp \mapsto F_P(\xp,\hcone(\hp_1(\xp,\tau),\tau))$ is Lipschitz continuous and since the solution component $\tau$ of $\HS_1$ increases linearly, then by \cite[Theorem~3.2]{khalil2002nonlinear}, $\HS_1$ in \eqref{eqn:H1} has no finite escape time from $\reals^{2n} \times \reals_{\geq 0}$ and each maximal solution to $\HS_0$ is unique. Therefore, each maximal solution to $\HS_1$, defined via \eqref{eqn:H1}, is complete and unique. 

Since $\HS_0$ has no finite time escape from $\reals^{2n}$ and $\HS_1$ has no finite time escape from $\reals^{2n} \times \reals_{\geq 0}$, then this means $\dot{x} = F(x)$ has no finite time escape from $C$ for $\HS$, as $\xlogic$ does not change in $C$ and as the state component $\tau$ is bounded in $C$, namely, the state component $\tau$ -- which is always reset to $0$ in $D$ -- increases linearly in $C_1$ and remains at $0$ in $C_0$. Therefore, there is no finite time escape from $C \cup D$, for solutions $x$ to $\HS$. Therefore, item (b) from Proposition \ref{prop:SolnExistence} does not hold.
\section{Proof of Proposition \ref{prop:ConvergenceNSCVXNesterov}}
\label{sec:ProofProp54}
The Lyapunov function $V_1$, defined via \eqref{eqn:LyapunovNesterovNSCVX}, is positive definite with respect to $\mathcal{A}_1$, defined via \eqref{eqn:setA_UGA}, 
since, by Assumption \ref{ass:LisNSCVX}, $L$ is $\mathcal{C}^1$, convex, and has a unique minimizer $\xp_1^*$. 	
Then, letting
\begin{equation} \label{eqn:BarV1}
	\bar{v}_1(\xp,\tau) := \xp_1 + \bar{\beta}(\tau)\xp_2,
\end{equation}
letting $\varphi(\xp,\tau) := \bar{a}(\tau)\left(\bar{a}(\tau)\left( \xp_1 - \xp_1^* \right) + \xp_2\right) +\frac{\zeta^2}{M}\nabla L(\xp_1)$, and since \IfConf{$\nabla V_1(\xp,\tau) $$ = \left[\varphi(\xp,\tau) \ \ \left(\bar{a}(\tau)\left( \xp_1 - \xp_1^* \right) + \xp_2\right) \right.$\\$\left. \diff{\bar{a}(\tau)}{\tau} \left\langle  \xp_1 - \xp_1^*,\left(\bar{a}(\tau)\left( \xp_1 - \xp_1^* \right) + \xp_2\right) \right\rangle \right]$}{$\nabla V_1(\xp,\tau) $\\$ = \left[\varphi(\xp,\tau) \ \ \left(\bar{a}(\tau)\left( \xp_1 - \xp_1^* \right) + \xp_2\right) \quad \diff{\bar{a}(\tau)}{\tau} \left\langle  \xp_1 - \xp_1^*,\left(\bar{a}(\tau)\left( \xp_1 - \xp_1^* \right) + \xp_2\right) \right\rangle \right]$}, we evaluate the derivative of $V_1$, using the map $\xp \mapsto$ $F_P(\xp,\kappa_1(\hp_1(\xp,\tau),\tau))$, where $F_P$ is defined in \eqref{eqn:HBFplant-dynamicsTR}, $\hcone$ is defined via \eqref{eqn:StaticStateFeedbackLawGlobalNSCVX}, and $\hp_1$ is defined in \eqref{eqn:H0H1NSCNesterovHBF}, to yield \IfConf{
	\begin{align} \label{eqn:dotVNSCVX}
		\dot{V}_1(\xp,\tau) = & \left\langle \nabla V_1(\xp,\tau),\matt{F_P(\xp,\kappa_1(\hp_1(\xp,\tau),\tau))\\ 1} \right\rangle \nonumber\\
		= & -\frac{\bar{a}({\tau})\zeta^2}{M}\left\langle \xp_1 - \xp_1^*,\nabla L(\bar{v}_1(\xp,\tau))\right\rangle \nonumber\\
		& + \bar{a}(\tau) \diff{\bar{a}(\tau)}{\tau}\left| \xp_1 - \xp_1^* \right|^2 \nonumber\\
		&+ \left(\bar{a}(\tau) - 2\bar{d}(\tau)\right)\left|\xp_2\right|^2 \nonumber\\
		&+ \left(\bar{a}^2(\tau) - 2\bar{d}(\tau)\bar{a}(\tau) + \diff{\bar{a}(\tau)}{\tau}\right)\left\langle \xp_1 - \xp_1^*,\xp_2 \right\rangle \nonumber\\
		&- \frac{\zeta^2}{M} \left\langle \xp_2, \nabla L(\bar{v}_1(\xp,\tau)) - \nabla L(\xp_1) \right\rangle 
	\end{align}
}{\begin{align} \label{eqn:dotVNSCVX_TR}
		\dot{V}_1(\xp,\tau) = & \left\langle \nabla V_1(\xp,\tau),\matt{F_P(\xp,\kappa_1(\hp_1(\xp,\tau),\tau))\\ 1} \right\rangle \nonumber\\
		= & \left\langle \nabla V_1(\xp,\tau), \matt{\matt{\xp_2 \\ - 2\bar{d}(\tau)\xp_2 - \frac{\zeta^2}{M}\nabla L(\bar{v}_1(\xp,\tau))} \\ 1} \right\rangle \nonumber\\
		= & \bar{a}(\tau) \left\langle \bar{a}(\tau)\left(\xp_1 - \xp_1^*\right) + \xp_2, \xp_2 \right\rangle + \frac{\zeta^2}{M}\left\langle \xp_2,\nabla L(\xp_1) \right\rangle - 2\bar{d}(\tau) \left|\xp_2\right|^2 \nonumber\\	
		&- 2\bar{d}(\tau)\bar{a}(\tau)\left\langle \xp_1 - \xp_1^*,\xp_2 \right\rangle -\frac{\bar{a}(\tau)\zeta^2}{M} \left\langle \xp_1 - \xp_1^*,\nabla L(\bar{v}_1(\xp,\tau)) \right\rangle \nonumber\\ 
		& -\frac{\zeta^2}{M} \left\langle \xp_2,\nabla L(\bar{v}_1(\xp,\tau)) \right\rangle +\bar{a}(\tau)\diff{\bar{a}(\tau)}{\tau}\left|\xp_1 - \xp_1^*\right|^2 +\diff{\bar{a}(\tau)}{\tau}\left\langle \xp_1 - \xp_1^*,\xp_2 \right\rangle \nonumber\\
		= & -\frac{\bar{a}({\tau})\zeta^2}{M}\left\langle \xp_1 - \xp_1^*,\nabla L(\bar{v}_1(\xp,\tau))\right\rangle + \bar{a}(\tau) \diff{\bar{a}(\tau)}{\tau}\left| \xp_1 - \xp_1^* \right|^2 \nonumber\\
		&+ \left(\bar{a}(\tau) - 2\bar{d}(\tau)\right)\left|\xp_2\right|^2 + \left(\bar{a}^2(\tau) - 2\bar{d}(\tau)\bar{a}(\tau) + \diff{\bar{a}(\tau)}{\tau}\right)\left\langle \xp_1 - \xp_1^*,\xp_2 \right\rangle \nonumber\\
		&- \frac{\zeta^2}{M} \left\langle \xp_2, \nabla L(\bar{v}_1(\xp,\tau)) - \nabla L(\xp_1) \right\rangle 
	\end{align}
}
for all $(\xp, \tau) \in \reals^{2n} \times \reals_{\geq 0}$. Since $L$ is $\mathcal{C}^1$, convex, and has a unique minimizer by Assumption \ref{ass:LisNSCVX}, then using the definition of convexity in Footnote \ref{foot:Convexity} with $u_1 = \xp_1^*$ and $w_1 = \bar{v}_1(\xp,\tau)$, where $\bar{v}_1$ is defined via \eqref{eqn:BarV1}, we get
\begin{align} \label{eqn:ConvexityZstarBarV}
	- \left\langle \bar{v}_1(\xp,\tau) - \xp_1^*, \nabla L(\bar{v}_1(\xp,\tau)) \right\rangle & \leq -\left(L(\bar{v}_1(\xp,\tau)) - L^*\right)
\end{align} 
for each $\xp \in \reals^{2n}$ and $\tau \in \reals_{\geq 0}$. Using the definition of convexity in Footnote \ref{foot:Convexity} with $u_1 = \bar{v}_1(\xp,\tau)$, where $\bar{v}_1$ is defined via \eqref{eqn:BarV1},  and $w_1 = \xp_1$ yields
\begin{align} \label{eqn:ConvexityZstarZ}
	\left\langle \nabla L(\xp_1), \bar{\beta}(\tau)\xp_2 \right\rangle & \leq L(\bar{v}_1(\xp,\tau)) - L(\xp_1)
\end{align}
for each $\xp \in \reals^{2n}$ and $\tau \in \reals_{\geq 0}$.
Combining \eqref{eqn:ConvexityZstarBarV} and \eqref{eqn:ConvexityZstarZ} yields\\ $- \left\langle \bar{v}_1(\xp,\tau) - \xp_1^*,\! \nabla L(\bar{v}_1(\xp,\tau)) \right\rangle + \left\langle \nabla L(\xp_1), \bar{\beta}(\tau)\xp_2 \right\rangle \leq \!\! -L(\bar{v}_1(\xp,\tau)) + L(\bar{v}_1(\xp,\tau)) - L(\xp_1) + L^*$.
%
Then, rearranging terms gives, for all $\xp \in \reals^{2n}$ and $\tau \in \reals_{\geq 0}$, 
\begin{align}\label{eqn:BoundOnQAndGradFNSCVX}
	\IfConf{& - \langle \xp_1 - \xp_1^*, \nabla L(\bar{v}_1(\xp,\tau)) \rangle \\
	& \leq - \left(L(\xp_1) - L^*\right)+ \left\langle \bar{\beta}(\tau) \xp_2, \nabla L(\bar{v}_1(\xp,\tau)) - \nabla L(\xp_1) \right\rangle . \nonumber}{& - \langle \xp_1 - \xp_1^*, \nabla L(\bar{v}_1(\xp,\tau)) \rangle \\ 
	& \leq - \left(L(\xp_1) - L^*\right)
	+ \left\langle \bar{\beta}(\tau) \xp_2, \nabla L(\bar{v}_1(\xp,\tau)) - \nabla L(\xp_1) \right\rangle . \nonumber}
\end{align}
Substituting the bound in \eqref{eqn:BoundOnQAndGradFNSCVX} into \IfConf{\eqref{eqn:dotVNSCVX}}{\eqref{eqn:dotVNSCVX_TR}} yields 
\begin{align} \label{eqn:SubsV1Bound}
	\IfConf{\dot{V}_1(\xp,\tau) \leq & - \frac{\bar{a}(\tau)\zeta^2}{M} \left(L(\xp_1) - L^*\right) \nonumber\\
		& + \frac{\bar{a}(\tau)\zeta^2}{M} \left\langle \bar{\beta}(\tau) \xp_2, \nabla L(\bar{v}_1(\xp,\tau)) - \nabla L(\xp_1) \right\rangle  \nonumber\\
		& + \bar{a}(\tau) \diff{\bar{a}(\tau)}{\tau}\left| \xp_1 - \xp_1^* \right|^2 + \left(\bar{a}(\tau) - 2\bar{d}(\tau)\right)\left|\xp_2\right|^2\nonumber\\
		& + \left(\bar{a}^2(\tau) - 2\bar{d}(\tau)\bar{a}(\tau) + \diff{\bar{a}(\tau)}{\tau}\right)\left\langle \xp_1 - \xp_1^*,\xp_2 \right\rangle \nonumber\\
		&- \frac{\zeta^2}{M} \left\langle \xp_2, \nabla L(\bar{v}_1(\xp,\tau)) - \nabla L(\xp_1) \right\rangle}{
	\dot{V}_1(\xp,\tau) \leq & - \frac{\bar{a}(\tau)\zeta^2}{M} \left(\left(L(\xp_1) - L^*\right)
	- \left\langle \bar{\beta}(\tau) \xp_2, \nabla L(\bar{v}_1(\xp,\tau)) - \nabla L(\xp_1) \right\rangle\right)  \nonumber\\
	& + \bar{a}(\tau) \diff{\bar{a}(\tau)}{\tau}\left| \xp_1 - \xp_1^* \right|^2 + \left(\bar{a}(\tau) - 2\bar{d}(\tau)\right)\left|\xp_2\right|^2\nonumber\\
	& + \left(\bar{a}^2(\tau) - 2\bar{d}(\tau)\bar{a}(\tau) + \diff{\bar{a}(\tau)}{\tau}\right)\left\langle \xp_1 - \xp_1^*,\xp_2 \right\rangle \nonumber\\
	&- \frac{\zeta^2}{M} \left\langle \xp_2, \nabla L(\bar{v}_1(\xp,\tau)) - \nabla L(\xp_1) \right\rangle}
\end{align}
for all $(\xp, \tau) \in \reals^{2n} \times \reals_{\geq 0}$. Then, noticing that $\frac{\bar{a}(\tau)}{2}\left|\bar{a}(\tau)\left(\xp_1 - \xp_1^*\right) + \xp_2 \right|^2 =$\IfConf{}{\\}$\frac{\bar{a}^3(\tau)}{2}\left| \xp_1 - \xp_1^* \right|^2$\IfConf{\\}{}$+ \bar{a}^2(\tau)\left\langle \xp_1 - \xp_1^*,\xp_2 \right\rangle + \frac{\bar{a}(\tau)}{2}\left|\xp_2\right|^2$,
adding it to and subtracting it from \eqref{eqn:SubsV1Bound}, and rearranging terms, yields \IfConf{
	\begin{align} \label{eqn:PreVanishBound}
		& \dot{V}_1(\xp,\tau) \\
		& \leq -\bar{a}(\tau)V_1(\xp,\tau) + \left(\frac{\bar{a}^3(\tau)}{2} + \bar{a}(\tau)\diff{\bar{a}(\tau)}{\tau}\right) \left| \xp_1 - \xp_1^* \right|^2\nonumber\\  
		& +\left(\frac{3\bar{a}(\tau)}{2} - 2\bar{d}(\tau)\right)\left|\xp_2\right|^2\nonumber\\
		& + \left(2\bar{a}^2(\tau) - 2\bar{d}(\tau)\bar{a}(\tau) + \diff{\bar{a}(\tau)}{\tau}\right)\left\langle \xp_1 - \xp_1^*,\xp_2 \right\rangle \nonumber\\
		& - \frac{\zeta^2}{M}\left(1 - \bar{\beta}(\tau)\bar{a}(\tau)\right)\left\langle \xp_2, \nabla L(\bar{v}_1(\xp,\tau)) - \nabla L(\xp_1) \right\rangle \nonumber
\end{align}}{
	\begin{align} \label{eqn:PreVanishBoundTR}
		\dot{V}_1(\xp,\tau) \leq & -\bar{a}(\tau)V_1(\xp,\tau) + \bar{a}(\tau) \diff{\bar{a}(\tau)}{\tau}\left| \xp_1 - \xp_1^* \right|^2 + \left(\bar{a}(\tau) - 2\bar{d}(\tau)\right)\left|\xp_2\right|^2 \nonumber\\
		& + \left(\bar{a}^2(\tau) - 2\bar{d}(\tau)\bar{a}(\tau) + \diff{\bar{a}(\tau)}{\tau}\right) \left\langle \xp_1 - \xp_1^*,\xp_2 \right\rangle + \frac{\bar{a}^3(\tau)}{2}\left| \xp_1 - \xp_1^* \right|^2  \nonumber\\
		& + \frac{\bar{a}(\tau)}{2}\left|\xp_2\right|^2 + \bar{a}^2(\tau)\left\langle \xp_1 - \xp_1^*,\xp_2 \right\rangle \nonumber\\
		& - \frac{\zeta^2}{M} \left(1 - \bar{\beta}(\tau)\bar{a}(\tau)\right) \left\langle \xp_2, \nabla L(\bar{v}_1(\xp,\tau)) - \nabla L(\xp_1) \right\rangle \nonumber\\
		\leq & -\bar{a}(\tau)V_1(\xp,\tau) + \left(\frac{\bar{a}^3(\tau)}{2} + \bar{a}(\tau)\diff{\bar{a}(\tau)}{\tau}\right) \left| \xp_1 - \xp_1^* \right|^2 \nonumber\\
		& +\left(\frac{3\bar{a}(\tau)}{2} - 2\bar{d}(\tau)\right)\left|\xp_2\right|^2\nonumber\\
		& + \left(2\bar{a}^2(\tau) - 2\bar{d}(\tau)\bar{a}(\tau) + \diff{\bar{a}(\tau)}{\tau}\right)\left\langle \xp_1 - \xp_1^*,\xp_2 \right\rangle \nonumber\\
		& - \frac{\zeta^2}{M}\left(1 - \bar{\beta}(\tau)\bar{a}(\tau)\right)\left\langle \xp_2, \nabla L(\bar{v}_1(\xp,\tau)) - \nabla L(\xp_1) \right\rangle 
\end{align}}
for all $(\xp, \tau) \in \reals^{2n} \times \reals_{\geq 0}$.
Due to the definitions of the functions $\bar{a}$ and $\bar{d}$, in \eqref{eqn:BarA} and \eqref{eqn:dBarBetaBar}, respectively, the cross term $\left\langle \xp_1 - \xp_1^*,\xp_2 \right\rangle$ vanishes since $2\bar{a}^2(\tau) - 2\bar{d}(\tau)\bar{a}(\tau) + \diff{\bar{a}(\tau)}{\tau} = 2\left(\frac{2}{\tau + 2}\right)^2 - 2\left(\frac{3}{2(\tau+2)}\right)\left(\frac{2}{\tau + 2}\right) - \frac{2}{\left(\tau + 2\right)^2} = 0$.
Moreover, the definitions of the functions $\bar{d}$ and $\bar{a}$ lead to the $\left| \xp_1 - \xp_1^* \right|^2$ and $\left|\xp_2\right|^2$ terms in \IfConf{\eqref{eqn:PreVanishBound}}{\eqref{eqn:PreVanishBoundTR}} vanishing due to $\frac{\bar{a}^3(\tau)}{2} + \bar{a}(\tau)\diff{\bar{a}(\tau)}{\tau} = \frac{\left(\frac{2}{\tau + 2}\right)^3}{2} + \left(\frac{2}{\tau + 2}\right)\left(-\frac{2}{(\tau + 2)^2}\right) = 0$
and $\frac{3\bar{a}(\tau)}{2} - 2\bar{d}(\tau) = \frac{3\left(\frac{2}{\tau + 2}\right)}{2} - 2\left(\frac{3}{2\left(\tau + 2\right)}\right) = 0$.
The bound in \IfConf{\eqref{eqn:PreVanishBound}}{\eqref{eqn:PreVanishBoundTR}} reduces to
\begin{align} \label{eqn:finalDotVNSCVX}
	& \IfConf{\dot{V}_1(\xp,\tau) \\
	& \leq -\bar{a}(\tau)V_1(\xp,\tau) \nonumber\\
	& - \frac{\zeta^2}{M}\left(1 - \bar{\beta}(\tau)\bar{a}(\tau)\right)\left\langle \xp_2, \nabla L(\bar{v}_1(\xp,\tau)) - \nabla L(\xp_1) \right\rangle \nonumber}{
	\dot{V}_1(\xp,\tau) \leq -\bar{a}(\tau)V_1(\xp,\tau) - \frac{\zeta^2}{M}\left(1 - \bar{\beta}(\tau)\bar{a}(\tau)\right)\left\langle \xp_2, \nabla L(\bar{v}_1(\xp,\tau)) - \nabla L(\xp_1) \right\rangle}
\end{align}
for all $(\xp, \tau) \in \reals^{2n} \times \reals_{\geq 0}$.
By Assumption \ref{ass:LisNSCVX}, $L$ is $\mathcal{C}^1$ and convex. By \cite[Theorem~2.1.3]{nesterov2004introductory}, a function $L$ is $\mathcal{C}^1$ and convex if and only if, for each $w_1, u_1 \in \reals^n$, 
\begin{equation}\label{eqn:ConvexityNest}
	\left\langle \nabla L(w_1) - \nabla L(u_1),w_1 - u_1 \right\rangle \geq 0.
\end{equation} 
Then, since $\bar{\beta}(\tau) \geq 0$ for all $t \geq 0$, using the bound in \eqref{eqn:ConvexityNest} 
with $w_1 = \bar{v}_1(\xp,\tau)$, where $\bar{v}_1$ is defined in \eqref{eqn:BarV1}, and $u_1 = \xp_1$, we get, for all $\xp \in \reals^{2n}$ and $\tau \in \reals_{\geq 0}$,
\begin{align} \label{eqn:UpperBoundLastTerm}
	\left\langle \bar{v}_1(\xp,\tau) - \xp_1, \nabla L(\bar{v}_1(\xp,\tau)) - \nabla L(\xp_1) \right\rangle &= \nonumber\\
	\bar{\beta}(\tau)\left\langle \xp_2, \nabla L(\bar{v}_1(\xp,\tau)) - \nabla L(\xp_1) \right\rangle &\geq 0 \nonumber\\
	- \bar{\beta}(\tau)\left\langle \xp_2, \nabla L(\bar{v}_1(\xp,\tau)) - \nabla L(\xp_1) \right\rangle &\leq 0 \nonumber\\
	- \left\langle \xp_2, \nabla L(\bar{v}_1(\xp,\tau)) - \nabla L(\xp_1) \right\rangle & \leq 0
\end{align}
Therefore, since $1 - \bar{\beta}(\tau)\bar{a}(\tau) \geq 0$, due to $\bar{a}$, defined via \eqref{eqn:BarA}, equaling $1$ at $\tau = 0$ and monotonically decreasing toward zero (but being always positive) as $\tau$ tends to $\infty$, and due to $\bar{\beta}$, defined via \eqref{eqn:dBarBetaBar}, equaling $0$ at $\tau = 0$ and monotonically increasing to $1$ as $\tau$ tends to $\infty$, we use \eqref{eqn:UpperBoundLastTerm}
to upper bound the last term of \eqref{eqn:finalDotVNSCVX} as follows:
\begin{equation}
	- \frac{\zeta^2}{M}\left(1 - \bar{\beta}(\tau)\bar{a}(\tau)\right)\left\langle \xp_2, \nabla L(\bar{v}_1(\xp,\tau)) - \nabla L(\xp_1) \right\rangle \leq 0
\end{equation}
This leads to, $\xp \in \reals^{2n}$ and $\tau \in \reals_{\geq 0}$,
\begin{equation} \label{eqn:MJ_dotVNCVX}
	\dot{V}_1(\xp,\tau) \leq -\bar{a}(\tau)V_1(\xp,\tau).
\end{equation}
Applying Gr\"{o}nwall's Inequality to \eqref{eqn:MJ_dotVNCVX}, namely,
\begin{align*}
	V_1(\xp(t),t) & \leq V_1(\xp(0),0) \exp \left(- \int_{0}^{t} \bar{a}(\tau) d\tau\right) \\ 
	& = V_1(\xp(0),0) \exp \left(-2\ln\left(t+2\right) -2\ln\left(2\right)\right) \\
	& = V_1(\xp(0),0)\exp \left(-\ln \left(\frac{t+2}{2}\right)^2\right) \\
	& = V_1(\xp(0),0)\left(\frac{1}{\exp \left(\ln \left(\frac{t+2}{2}\right)^2\right)}\right) \\
	& = \frac{4}{(t+2)^2}V_1(\xp(0),0)
\end{align*}
shows that each maximal solution $t \mapsto (\xp(t),\tau(t))$ to the closed-loop algorithm $\HS_1$, such that $\tau(0) = 0$, satisfies \eqref{eqn:ConvergenceRateNCVX}, for all $t \geq 0$. 
\section{General Results for Hybrid Systems}
\label{sec:GenHybridResults}

The following proposition, from \cite{65}, is used to prove the existence of solutions to the hybrid closed-loop system. 
\IfConf{\begin{prop}[Basic existence of solutions]}{\begin{proposition}(Basic existence of solutions):} \label{prop:SolnExistence}
	Let $\HS = (C,F,D,G)$ satisfy Definition \ref{def:HBCs}. Take an arbitrary $\xi \in C \cup D$. If $\xi \in D$ or\\
	
	\noindent
	(VC) there exists a neighborhood $U$ of $\xi$ such that for every $x \in U \cap C$,	
	$$F(x) \cap T_C(x) \neq \emptyset,$$	
	then there exists a nontrivial solution $x$ to $\HS$ with $x(0,0) = \xi$. If (VC) holds for every $\xi \in C \setminus D$, then there exists a nontrivial solution to $\HS$ from every initial point in $C \cup D$, and every\footnote{The set $\mathcal{S}_{\HS}$ contains all maximal solutions to $\HS$.} $x \in \mathcal{S}_{\HS}$ satisfies exactly one of the following conditions:
	
	\begin{enumerate}[label={(\alph*)}]
		\item \label{item:SolnExistenceItemA} $x$ is complete;
		
		\item \label{item:SolnExistenceItemB} $\mathrm{dom} \ x$ is bounded and the interval $I^J$, where $J = \mathrm{sup}_j \ \mathrm{dom} \ x$, has nonempty interior and $t \mapsto x(t,J)$ is a maximal solution to $\dot{z} \in F(z)$, in fact\\ $\mathrm{lim}_{t \mapsto T} \lvert x(t,J) \rvert = \infty$, where $T = \mathrm{sup}_t \ \mathrm{dom} \ x$;
		
		\item \label{item:SolnExistenceItemC} $x(T,J) \not\in C \cup D$, where $(T,J) = \mathrm{sup} \ \mathrm{dom} \ x$.
	\end{enumerate}
	Furthermore, if $G(D) \subset C \cup D$, then (c) above does not occur. 	
\IfConf{\end{prop}}{\end{proposition}}

The following definition, from \cite[Definition~3.17]{220}, describes the basic properties that a function must satisfy to serve as a Lyapunov function for the hybrid closed-loop algorithm $\HS$. 
\IfConf{\begin{defn}[Lyapunov function candidate]}{\begin{definition}[Lyapunov function candidate]} \label{def:LyapCandidate}
	The sets $\mathcal{U}$, $\mathcal{A} \subset \reals^n$, and the function $V : \dom V \rightarrow \reals$ define a \underline{Lyapunov function candidate} on $\mathcal{U}$ with respect to $\mathcal{A}$ for the hybrid closed-loop system $\HS = \left(C, F, D, G\right)$ if the following conditions hold:
	\begin{enumerate}
		\item $\left(\overline{C} \cup D \cup G(D) \right) \cup \mathcal{U} \subset \dom V$;
		\item $\mathcal{U}$ contains an open neighborhood of\IfConf{\\}{} $\mathcal{A} \cap \left(C \cup D \cup G(D)\right)$;
		\item $V$ is continuous on $\mathcal{U}$ and locally Lipschitz on an open set containing $\overline{C} \cap \mathcal{U}$;
		\item $V$ is positive definite on $\overline{C} \cup D \cup G(D)$ with respect to $\mathcal{A}$. 
	\end{enumerate}
\IfConf{\end{defn}}{\end{definition}}

The following theorem is used to prove the uniform global asymptotic stability of the hybrid closed-loop system, via Lyapunov stability and an invariance principle. 

\IfConf{\begin{thm}[Hybrid Lyapunov theorem]}{\begin{theorem}(Hybrid Lyapunov theorem):} \label{thm:hybrid Lyapunov theorem}
	Given sets $\mathcal{U}, \mathcal{A} \subset \reals^n$ and a function $V : \mathrm{dom} \ V \rightarrow \reals$ defining a Lyapunov candidate on $\mathcal{U}$ with respect to $\mathcal{A}$ for the closed-loop hybrid system $\HS = (C,F,D,G)$, suppose
	
	\begin{itemize}
		\item $\HS$ satisfies the hybrid basic conditions;
		\item $\mathcal{A}$ is compact and $\mathcal{U}$ contains a nonzero open neighborhood of $\mathcal{A}$;
		\item $\dot{V}$ and $\Delta V$ satisfy
		\begin{align}
			& \dot{V}(x) = \max_{\xi \in \f(x)} \langle \nabla V(x),\xi \rangle \leq 0 \qquad \qquad \forall x \in C \cap \mathcal{U} \label{eqn:CTLyapunov}\\
			& \Delta V(x) := \max_{\xi \in \g(x)} V(\xi) - V(x) \leq 0 \quad \quad \ \forall x \in D \cap \mathcal{U} \label{eqn:DTLyapunov}
		\end{align}
	\end{itemize}
	Then $\mathcal{A}$ is stable. Furthermore, $\mathcal{A}$ is attractive and, hence, pre-asymptotically stable if any of the following conditions hold: 
	
	\begin{enumerate}
		\item Strict decrease during flows and jumps:
		\begin{align}
			& \dot{V}(x) < 0 \qquad \qquad \qquad \quad \forall x \in (C \cap \mathcal{U}) \backslash \mathcal{A} \\
			& \Delta V(x) < 0 \qquad \qquad \qquad \; \forall x \in (D \cap \mathcal{U}) \backslash \mathcal{A} 
		\end{align}
		\item Strict decrease during flows and no instantaneous Zeno:
		\begin{enumerate}[label={(\alph*)}]
			\item $\dot{V}(x) < 0$ for each $x \in (C \cap \mathcal{U}) \backslash \mathcal{A}$,
			\item any instantaneous Zeno solution $x$ to $\HS$ where $\mathrm{rge} \ x \subset \mathcal{U}$ converges to $\mathcal{A}$;
		\end{enumerate}
		\item Strict decrease during jumps and no complete continuous solution:
		\begin{enumerate}[label={(\alph*)}]
			\item $\Delta V(x) < 0$ for each $x \in (D \cap \mathcal{U}) \backslash \mathcal{A}$,
			\item any complete continuous solution $x$ to $\HS$ where $\mathrm{rge} \ x \subset \mathcal{U}$ converges to $\mathcal{A}$;
		\end{enumerate}
		\item \label{item:HybridLyapunovWeak} Weak decrease during flows and jumps: for each $\chi \in \mathcal{U}$ with $r := V(\chi) > 0$ there is no complete solution $x$ to $\HS$, $x(0,0) = \chi$ such that
		\begin{equation}
			\mathrm{rge} \ x \subset \defset{x}{V(x) = r} \cap \mathcal{U}
		\end{equation}
		and the set $\mathcal{U}$ is the subset of the basin of pre-attraction.
	\end{enumerate}
\IfConf{\end{thm}}{\end{theorem}}

Observe that, if the set $\mathcal{A}$ is pre-asymptotically stable via Theorem \ref{thm:hybrid Lyapunov theorem} and the Lyapunov function $V$ also has compact sublevel sets, namely, for each $c_V > 0$, $\defset{x}{V(x) \leq c_V}$ is compact, then the origin is {\em globally pre-asymptotically stable}.

The following result is used to show that, when a hybrid closed-loop algorithm $\HS$ has a set $\mathcal{A}$ globally asymptotically stable, then when $\HS$ satisfies the hybrid basic conditions, the set $\mathcal{A}$ is also uniformly globally asymptotically stable\footnote{Uniform global asymptotic stability allows an equivalent characterization involving a class-$\mathcal{KL}$ function \cite{65}.} for $\HS$. 
\IfConf{\begin{thm}[pAS implies $\mathcal{KL}$ pAS]}{\begin{theorem}(Pre-asymptotic stability implies $\mathcal{KL}$ pre-asymptotic stability):}
	\label{thm:GASImpliesUGAS}
	Suppose that the hybrid closed-loop system $\HS$ satisfies the hybrid basic conditions and that a compact set $\mathcal{A}$ is pre-asymptotically stable with basin of pre-attraction $\mathcal{B}^p_{\mathcal{A}}$. Then, $\mathcal{B}^p_{\mathcal{A}}$ is open and $\mathcal{A}$ is $\mathcal{KL}$ pre-asymptotically stable on $\mathcal{B}^p_{\mathcal{A}}$ for $\HS$; namely, there exists a function $\beta \in \mathcal{KL}$ such that 
	\begin{equation}
		\left|x(0,0)\right|_{\mathcal{A}} \leq \beta \left(\left|x(0,0)\right|_{\mathcal{A}}, t+j \right) \quad \forall (t,j) \in \dom x
	\end{equation}
	for each $x \in \mathcal{S}_{\HS}(\mathcal{B}^p_{\mathcal{A}})$.
\IfConf{\end{thm}}{\end{theorem}}

For Proposition \ref{prop:GAS-HBF} and Theorem \ref{thm:HybridInvariancePrinciple} we use the following definition of weak invariance, from \cite{65}.
\IfConf{\begin{defn}[Weak invariance]}{\begin{definition}[Weak invariance]}
	Given a hybrid system $\HS$, a set $S \subset \reals^n$ is said to be
	\begin{itemize}
		\item weakly forward invariant if for every $\xi \in S$ there exists at least one complete $x \in \mathcal{S}_{\HS}(\xi)$ with $\mathrm{rge} \ x \subset S$;
		\item weakly backward invariant if for every $\xi \in S$ and every $T> 0$, there exists at least one $x \in \mathcal{S}_{\HS}(S)$ such that for some $(t^*, j^*) \in \dom x$, $t^* + j^* \geq T$, it is the case that $x(t^*, j^*) = \xi$ and $x(t,j) \in S$ for all $(t,j) \in \dom x$ with $t + j \leq t^* + j^*$;
		\item weakly invariant if it is both weakly forward invariant and weakly backward invariant.
	\end{itemize}
\IfConf{\end{definition}}{\end{definition}}

The following {\em hybrid invariance principle}, from \cite[Theorem~3.23]{220}, is used to establish attractivity when only a ``weak'' Lyapunov function is available -- meaning that the function does not strictly decrease along both flows and jumps of the hybrid system. It is also useful to check where particular solutions of interest converge to. 
\IfConf{\begin{thm}[Hybrid Invariance Principle]}{\begin{theorem}(Hybrid Invariance Principle):} \label{thm:HybridInvariancePrinciple}
	Given a hybrid closed-loop system $\HS = (C,F,D,G)$ with state $x \in \reals^n$ satisfying the hybrid basic conditions, nonempty $\mathcal{U} \subset \reals^n$, and a function $V : \dom V \rightarrow \reals$, suppose that \ref{def:LyapCandidate} is satisfied, and that \eqref{eqn:CTLyapunov} and \eqref{eqn:DTLyapunov} hold. With $X := C \cup D \cup G(D)$, we empoly the following definitions:
	\begin{align}
		V^{-1}(r) & := \defset{x \in X}{V(x) = r\!\!}\\
		\dot{V}^{-1}(0) & := \defset{x \in C}{\dot{V}(x) = 0\!\!}\\
		\Delta V^{-1}(0) & := \defset{x \in D}{\Delta V(x) = 0\!\!}
	\end{align}
	Let $x$ be a precompact solution to $\HS$ with $\overline{\rge x} \subset \mathcal{U}$. Then, for some $r \in V(\mathcal{U} \cap X)$, the following hold:
	
	\begin{enumerate}
		\item The solution $x$ converges to the largest weakly invariant set in 
		\begin{equation}
			V^{-1}(r) \cap \mathcal{U} \cap \left[\dot{V}^{-1}(0) \cup \left(\Delta V^{-1}(0) \cap G\left(\Delta V^{-1}(0)\right) \right) \right];
		\end{equation}
		\item The solution $x$ converges to the largest weakly invariant set in 
		\begin{equation}
			V^{-1}(r) \cap \mathcal{U} \cap \Delta V^{-1}(0) \cap G\left(\Delta V^{-1}(0)\right)
		\end{equation}
		if in addition the solution $X$ is Zeno;
		\item The solution $x$ converges to the largest weakly invariant set in 
		\begin{equation}
			V^{-1}(r) \cap \mathcal{U} \cap \dot{V}^{-1}(0)
		\end{equation}
		if, in addition, the solution $x$ is such that, for some $a > 0$ and some $J \in N$, $t_{j + 1} - t_j > a$ for all $j \geq J$; i.e., the given solution $x$ is such that the elapsed time between consecutive jumps is eventually bounded below by a positive constant $a$.
	\end{enumerate}
\IfConf{\end{thm}}{\end{theorem}}
}

\bibliographystyle{IEEEtran}
\bibliography{dhustigs}
\end{document}